\documentclass[11pt,reqno,English]{amsart}
\usepackage[letterpaper, margin=1in]{geometry}
\usepackage{amsfonts}
\usepackage{amssymb}
\usepackage{amsmath}
\usepackage{mathrsfs}
\usepackage{color}
\usepackage[nospace,noadjust]{cite}
\usepackage{verbatim}
\usepackage{dsfont}
\usepackage{bbm}
\usepackage{enumerate,booktabs,mathtools}
\usepackage{graphicx}
\usepackage{slashed} 
\usepackage[titletoc,title]{appendix}
\usepackage{wrapfig}
\usepackage{marginnote}
\usepackage{bm}
\usepackage{enumitem}

\usepackage{nameref}
\usepackage{varioref}
\usepackage{hyperref}
\usepackage{cleveref}
\usepackage{url}

\numberwithin{equation}{section}

\usepackage{algorithm}     % Añade algoritmos y pseudocódigo
\usepackage{algpseudocode}
\renewcommand{\v}{%
    \ensuremath{v}}

\newcommand{\be}{\begin{equation}}
\newcommand{\ee}{\end{equation}}
\newcommand{\bp}{\begin{proof}}
\newcommand{\ep}{\end{proof}}
\newcommand{\bel}{\begin{equation}\label}
\newcommand{\eeq}{\end{equation}}
\newcommand{\bea}{\begin{eqnarray}}
\newcommand{\eea}{\end{eqnarray}}
\newcommand{\bee}{\begin{eqnarray*}}
\newcommand{\eee}{\end{eqnarray*}}
\newcommand{\ben}{\begin{enumerate}}
\newcommand{\een}{\end{enumerate}}

\newcommand{\R}{\mathbb{R}}

\newtheorem{thm}{Theorem}[section]

\newtheorem{conj}[thm]{Conjecture}

\newtheorem*{prop*}{{\bf Proposition}}

\theoremstyle{definition}
\newtheorem{defn}[thm]{Definition}

\numberwithin{ej}{section}

\theoremstyle{remark}

\begin{document}

\title[Interaction of Jamitons]{Interaction of jamitons  in second-order macroscopic traffic models}

\author[Raimund B\"urger]{Raimund B\"urger}
\address{CI${}^{\mathrm{2}}$MA and Departamento de Ingenier\'{\i}a Matem\'atica, Facultad de Ciencias F\'{i}sicas y Matem\'{a}ticas, 
 Universidad de Concepci\'{o}n, Casilla 160-C, Concepci\'{o}n, Chile.}
\email{rburger@ing-mat.udec.cl}
\thanks{R.B. was partially spported by ANID (Chile) through projects Fondecyt 1210610;
  Anillo ANID/ACT210030;  Centro de Modelamiento Matem\'{a}tico (CMM), project FB210005 of BASAL funds for Centers of Excellence; and CRHIAM, 
   projects ANID/FONDAP/15130015 and ANID/FONDAP/1523A0001. }

\author[Claudio Mu\~noz]{Claudio Mu\~noz}
\address{Departamento de Ingenier\'{\i}a Matem\'atica and Centro
de Modelamiento Matem\'atico (UMI 2807 CNRS), Universidad de Chile, Casilla
170 Correo 3, Santiago, Chile.}
\email{cmunoz@dim.uchile.cl}
\thanks{C.M. was partially funded by Chilean research grants FONDECYT 1231250, Basal CMM FB210005. He was additionally supported by Grant PID2022-137228OB-I00 (DISCOLBE).}

\author[Sebasti\'an Tapia]{Sebasti\'an Tapia}
\address{Departamento de Ingenier\'{\i}a Matem\'atica, Universidad de Chile, Casilla
170 Correo 3, Santiago, Chile.}
\email{sebastian.tapia.st@gmail.com}
\thanks{S. T. was supported by Chilean research grants FONDECYT 1231250 and Basal CMM FB210005. He also deeply thanks INRIA Lille and professor Andr\'e de Laire for their support during a stage visit in early 2024, where part of this work was completed.}


\begin{abstract}
Jamitons are self-sustained traveling wave solutions that arise in certain second-order macroscopic models of vehicular traffic. 
  A necessary condition for a jamiton to appear is that the local traffic density 
   breaks the so-called sub-characteristic condition. This condition states that 
    the characteristic velocity of the corresponding first-order Lighthill-Whitham-Richards (LWR) 
     model formed with the same desired speed function is enclosed by the characteristic speeds of the 
      corresponding second-order model. The phenomenon of collision of jamitons 
       in second-order models of traffic flow is studied analytically and numerically for 
        the particular case of the second-order Aw-Rascle-Zhang (ARZ) traffic model 
        [A.\ Aw, M.\ Rascle,  {\it SIAM J.\ Appl.\ Math.} 
       {\bf  60} (2000) 916--938; H.\ M.\ Zhang,  {\it Transp.\ Res.\ B} {\bf   36} (2002)  275--290].
         A compatibility condition is first defined to select jamitons that can collide each other. The collision of jamitons produces a new jamiton with a velocity different from the initial ones. It is observed that the exit velocities smooth out the velocity of the test jamiton and the initial velocities of the  jamitons that collide. Other properties such as the amplitude of the exit jamitons, lengths, and maximum density are also explored.
          In the cases of the amplitude and maximum exit density  it turns out  that over a wide range of sonic densities, 
           the exit values exceed or equal the input values. On the other hand, the resulting jamiton has a greater length than the incoming ones. Finally, the behavior for various driver reaction times is explored. It is  obtained  that some properties do not depend on that  time, 
            such as the amplitude, exit velocity, or maximum density, while the exit length does depend on driver reaction time. 
\end{abstract}

\maketitle 

%\tableofcontents

%% TABLA DE CONTENIDOS - ÍNDICE
%\templateIndex
%
%% CONFIGURACIONES FINALES
%\templateFinalcfg

% CONFIGURACIONES FINALES
%\templateFinalcfg

% ======================= INICIO DEL DOCUMENTO =======================

\section{Introduction}

\subsection{Scope} 

\emph{Jamitons} are self-sustained traveling wave solutions that arise in second-order models of vehicular traffic. Here we understand as a second-order trafffic model as pair of 
 one-dimensional balance equations that represent the analogues 
  of conservation of mass and linear momentum in continuum mechanics. The notion of ``jamiton'', an apparent amalgamation 
 of traffic ``jam'' and     ``soliton'', was coined by Flynn et al.\ in  \cite{presion}, and studied in 
  the context of the fundamental diagram of traffic flow by Seibold et al.\ \cite{jamitones}.  
   The term  ``jamiton'' is, however, not yet universally used 
    in the traffic modelling literature (cf., e.g., \cite{Garavello-Piccoli, Treiber-Kesting}). 
  One widely studied  second-order model  that describes the formation of jamitons  is the inhomogeneous Aw-Rascle-Zhang (ARZ)  model 
   \cite{resurrection,Z} given by  
    \begin{align}\label{eq:ARZ0}
    \begin{split}
   \partial_t \rho  +  \rho_t \partial_x (\rho u)  &= 0,\\ 
  \partial_t \bigl(u + h(\rho) \bigr)  + u \partial_x \bigl(u + h(\rho) \bigr)  &= \dfrac{U(\rho)-u}{\tau},
    \end{split}
\end{align}
where~$t$ is time, $x$~is spatial position, $\rho= \rho(x,t)$ is the sought density of vehicles
 (assumed to take values between zero and some maximal density~$\rho_{\max}$), 
   $u=u(x,t)$ is the unknown velocity, $h=h( \rho)$ is a strictly increasing hesitation function, $U=U( \rho)$ is the desired speed that is 
  assumed to be a given decreasing function, and~$\tau>0$ is a relaxation time. 
It is the purpose of the present work to undertake a mathematical study of the {\em collision} of jamitons for \eqref{eq:ARZ0}, elucidating the behavior of a jamiton following its collision with another jamiton, studying its size and  output velocity. To the knowledge of the authors and known works, this type of study has not been undertaken  previously.

\subsection{Related work} From a mathematical standpoint, traffic models have been studied from at least three different perspectives in recent years, each distinguished by its own particular approach, 
 namely   microscopic models, cellular automata models,    and macroscopic models. These approaches do not operate independently but are 
    consistent in both their formulations and  the phenomena encompassed. For instance, the simplest Nagel-Schreckenberg cellular automaton model \cite{Na-Sch} turns out to be a particular case of the numerical discretization of the macroscopic Lighthill-Whitham-Richards 
  (LWR) model \cite{LW,R}.  Similarly, the  discrete cell transmission model by Hilliges and Weidlich \cite{hw} represents a monotone numerical scheme for  the LWR model that can 
   even be extended to the multiclass case \cite{family,bcv23}. 
  We mainly focus here on macroscopic traffic models that  describe vehicular traffic  by a two-phase continuum approach.  Introductions to this class of models include \cite{Garavello-Piccoli,Treiber-Kesting,Jinbook2021}. 
  Macroscopic or continuous models   no longer follow the behavior of individual vehicles but rather focus on vehicle density and the velocity field (the velocity present at a particular time  and location).  They are also related to microscopic or cellular automaton models. For example, in \cite{micro-macro}, equivalences between continuous and microscopic models are established, while in \cite{cellular-macro}, equivalences with cellular automata are discussed. 
  
  This type of models presents several advantages for the study of vehicular traffic. For instance, they are preferred in terms of safety and data privacy as they do not access individual vehicle data directly. They also offer good precision and estimation at a large scale for noisy or sparse data \cite{gps_macro, gps_macro_2}. Moreover, they can be extended to routes with multiple lanes or to control problems \cite{control}. Summarizing, macroscopic models usually give rise to partial differential equations (PDEs) or systems of PDEs, and many features of traffic flow can be derived from their qualitative properties as well as by applying numerical methods, 
   see monographs on hyperbolic conservation laws and their numerical treatment (cf., e.g., \cite{shocks,volumen_finito,toro2009,hestbook,kuzminbook}).

The authors' interest in  collisions of traffic  flow interactions  that should  have a potentially complicated structure
 is in part motivated by the so-called conjectured fractal exit velocity that appears in dispersive collisions. 
 An important simplified model of these phenomena  is   the well-known $\phi^4$ equation 
 \begin{align}\label{eq:phi_4}
    \partial_{tt} \phi  - \partial_{xx}  \phi  - \phi + \phi^3=  0, \quad \phi(x,t)\in\mathbb R, \quad (x,t)\in\mathbb R^2.
\end{align}
 whose solutions may exhibit these collisions. 
  (Equation~\eqref{eq:phi_4} is   also known under various other names, such as ``$u$-4 model'' \cite{ablo11} or ``cubic Klein-Gordon equation'' 
   \cite{schneid17}.)  
  Equation~\eqref{eq:phi_4}  admits a family of traveling wave solutions, called {\em kinks}, that are given by
\begin{align}\label{phi4_kink} 
    \phi(x,t) = \phi_{\mathrm{K}} (x-vt) = \tanh\left( \dfrac{x - x_0 - vt}{\sqrt{2(1-v^2)}}\right),  
\end{align} 
where $v \in (-1,1)$ is the (fixed) input velocity. On the other hand, the antikink solution corresponds to the same wave but traveling in the opposite direction, 
 $ \phi_{\bar{\mathrm{K}}}(x,t) = - \phi_{\mathrm{K}} (x+ vt)$.   From numerical simulations, one observes that equation \eqref{eq:phi_4} 
 presents a fractal behavior for the output velocities after the collision, based on the input velocities. 
The phenomenon of the output velocities in the $\phi^4$ equation has been extensively studied in the literature, see e.g.\ \cite{phi4, phi_4_1, phi_4_2, phi_4_3, phi_4_4, phi_4_5}, where 
 the phenomenon known as the ``multiple-bounce" resonance effect is described mathematically.

\subsection{Outline of the paper}

In this paper we study, from a numerical point of view, collisions of jamitons, seeking to represent real-life jam collisions. This work is organized as follows. First of all, in Section \ref{cap:macro} we  present the macroscopic traffic models to be explored. Their properties, conditions, and 
 relation to  the fundamental diagram are  discussed. In Section~\ref{cap:jam}, the mathematical construction of jamitons is presented, including their definition, derivation from the traffic model, and several properties that they exhibit, in particular stability and asymptotic stability. In Section~\ref{cap:ARZ}
   a numerical scheme  for the simulation of  a jamiton is outlined. The scheme is  validated by comparison with  theoretical jamitons and results observed in the literature will be replicated, such as the emergence of jamitons on a long, essentially infinite, route. Finally, in Section~\ref{cap:colision}, the procedure for 
  the simulation of the collision of  jamitons is described. This includes the initial configuration, conditions that jamitons must meet to collide, and the selection of jamitons. 
   The chosen jamitons are caused to  collide, and possible properties arising from the collision will be studied. Section~\ref{cap:conclusion} presents conclusions and future work.

\section{Macroscopic models}\label{cap:macro}

\subsection{Preliminaries} 

Macroscopic models describe vehicular traffic  as two continuous phases (vehicles and the void space). 
 The    spatio-temporal evolution of vehicle density is modeled through scalar conservation  PDEs. The main variables 
  are usually the density $\rho = \rho(x, t)$, which represents the number of vehicles per unit length and will always is assumed nonnegative, the velocity field $u=u(x, t)$  that represents the
 local car  velocity at  position~$x$ at time~$t$, and the vehicle flow rate $Q = \rho u$, which indicates the number of vehicles passing a fixed point per unit of time. (Some important modifications to the previous setting are formulated within kinetic models. We will not consider these formulations here, but for a detailed introduction and results, including the collision phenomenon, see \cite{Bellouquid2012,Albi2019,Delitala2017}.) With these variables it is possible to deduce the principle of conservation of mass in integral form, 
 \begin{align}  \label{intform} 
  \frac{\mathrm{d}}{\mathrm{d} t} \int_a^b \rho(x,t) \, \mathrm{d} x = \rho(a,t) u(a,t) - \rho(b,t) u (b,t) \quad 
   \text{for any interval $[a,b] \subset \mathbb{R}$ and $t \geq 0$,} 
 \end{align} 
 which indicates that vehicles are neither created nor destroyed, and from \eqref{intform}  we deduce the conservation PDE in differential form 
\begin{equation}\label{eq:continuidad}
   \partial_t   \rho  +  \partial_x (\rho u) = 0. 
\end{equation}
  The model~\eqref{eq:continuidad}  is still incomplete 
 since  there are two unknowns ($\rho$ and~$u$) and only one equation. The well-known Lighthill-Whitham-Richards (LWR) model \cite{LW,R} 
  provides the required model closure via the kinematic assumption 
  $u = u(\rho)$. In contrast, second-order models such as the ARZ model \eqref{eq:ARZ0}   or the    Payne-Whitham (PW) model  \cite{P,W} (to be outlined in Section~\ref{subsec:PW}) 
    add a second equation for the velocity field that makes both variables~$\rho$ and~$u$ independent. Before continuing, it is necessary to introduce the concept of fundamental diagram  for both first- and second-order traffic models. The  fundamental diagram 
     condenses the main properties to be encapsulated in macroscopic models, such as vehicular congestion or free flow at low densities.

\subsection{Fundamental diagram}
\begin{figure}[t]
\centering
   \begin{tabular}{cc}(a) & (b) \\ 
   \includegraphics[width=.49\linewidth]{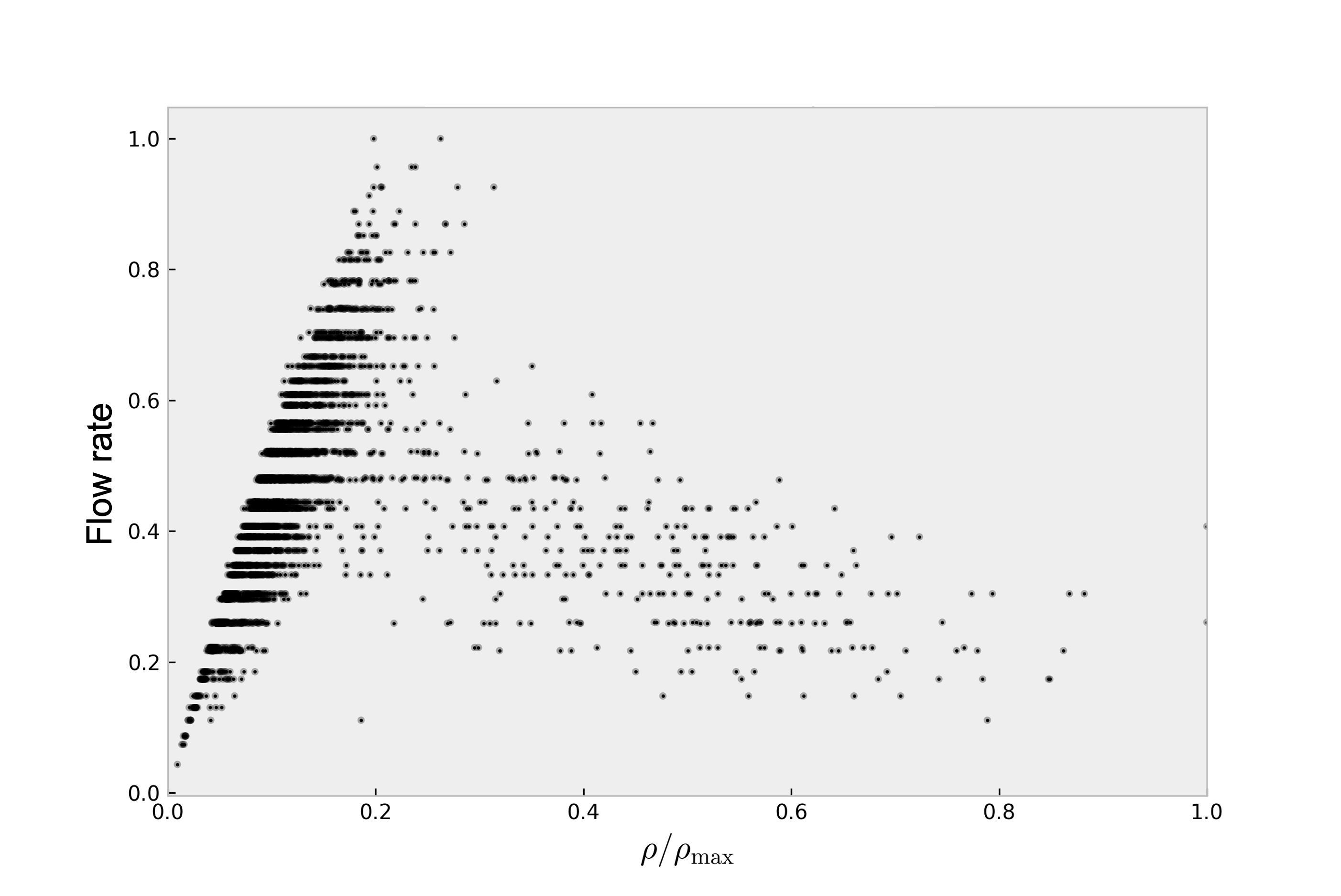} & 
   \includegraphics[width=.49\linewidth]{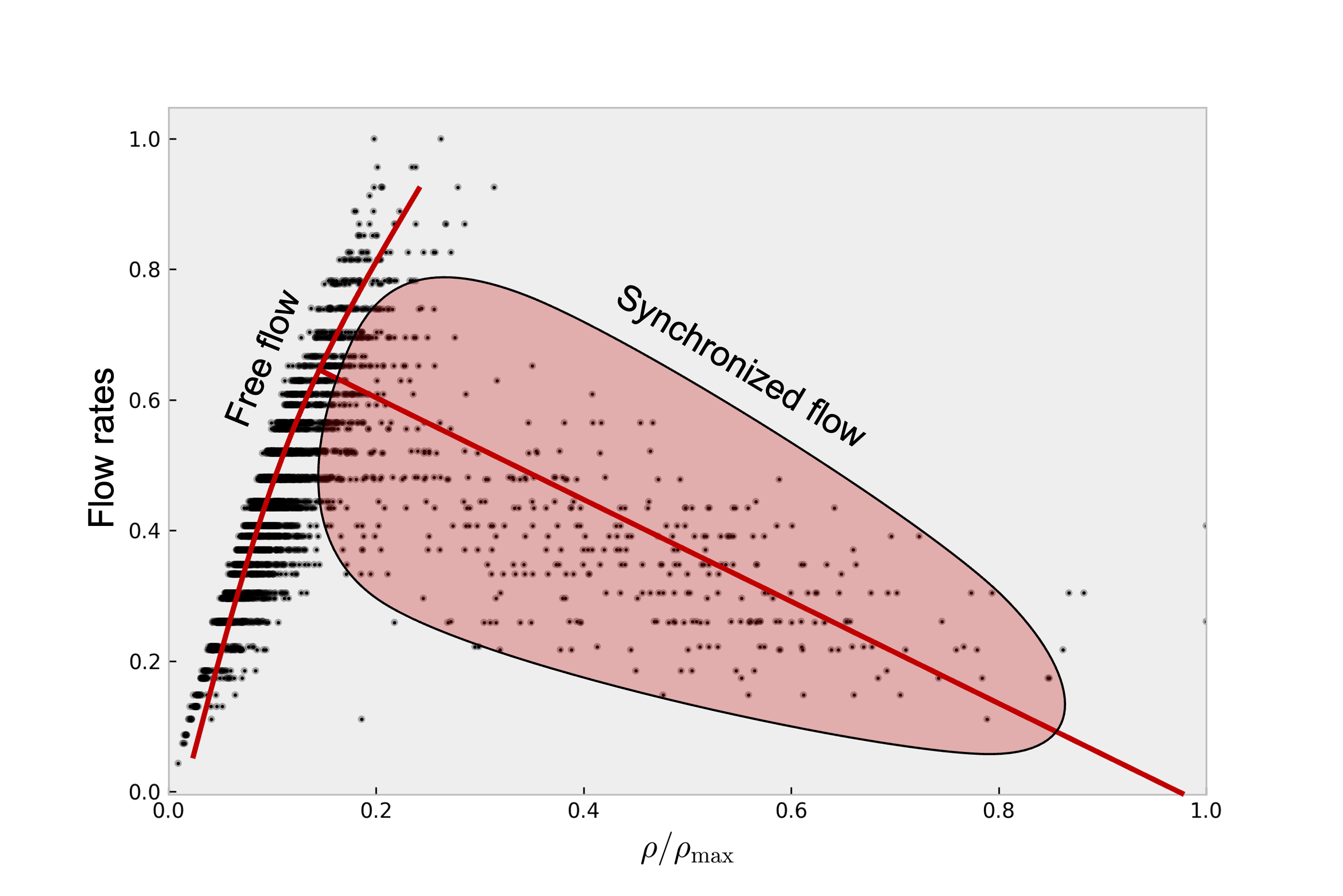} \end{tabular} 
   \caption{Example of an empirical  fundamental diagram with its phases: (a) fundamental diagram, (b) decomposition of the flow in two phases.}
   \label{img:DFs}
\end{figure}

The fundamental diagram is the plot of empirical measures of vehicle flow rate versus density. An example of a fundamental diagram obtained from a 2003 Minnesota Department of Transportation \cite{dataset} database of a highway sensor is shown in Figure~\ref{img:DFs}(a). In macroscopic models the fundamental diagram 
is approximated  by a unimodal function (i.e., a function that  has one extremum) $Q=Q(\rho)$. Then $Q$ is given by $Q(\rho) = \rho U(\rho)$, where $U(\rho)$ is known as the desired speed, 
 which is a decreasing  function of~$\rho$ (higher vehicle density implies less freedom of movement for vehicles). In other words, 
  the decreasing behavior of~$U( \rho)$ describes drivers' attitude to reduce speed with increasing  local traffic density~$\rho$.

\begin{figure}[t]
   \centering
   \begin{tabular}{cc} (a) & (b) \\ 
   \includegraphics[width=.49\linewidth]{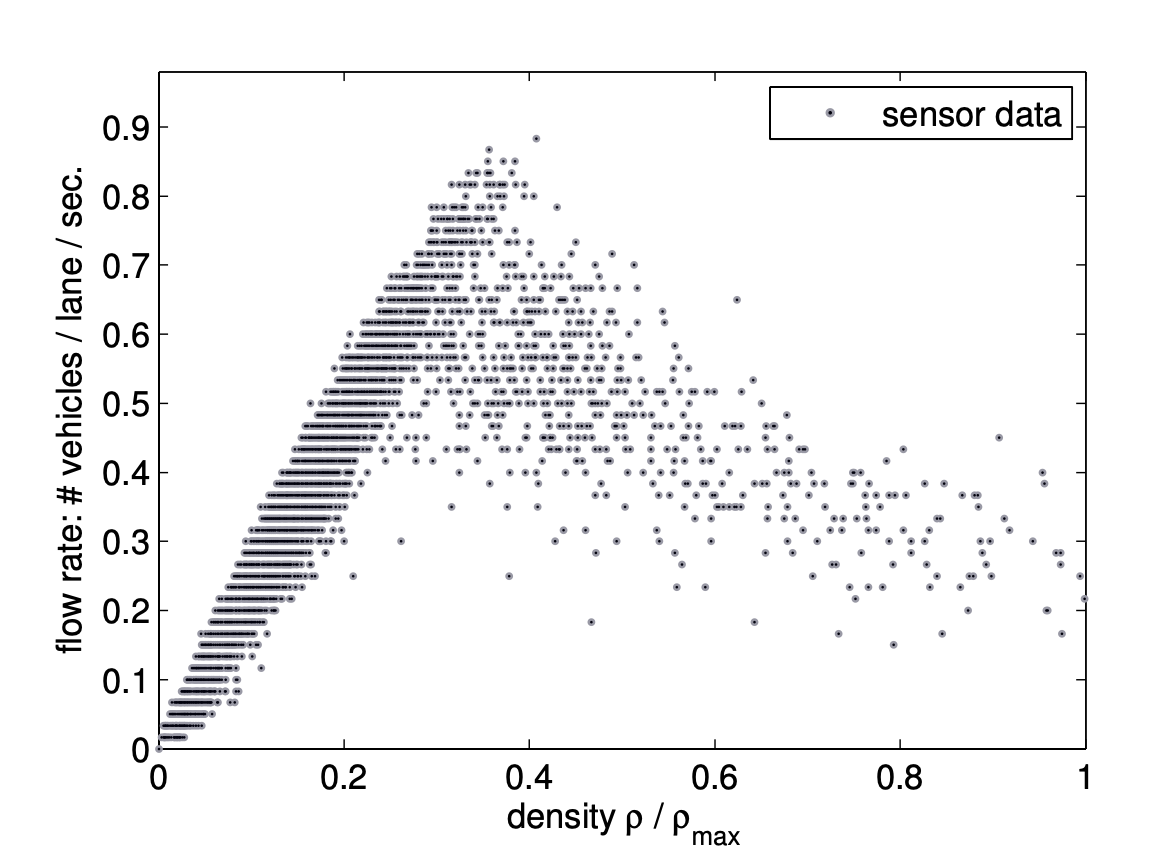} &  
   \includegraphics[width=.49\linewidth]{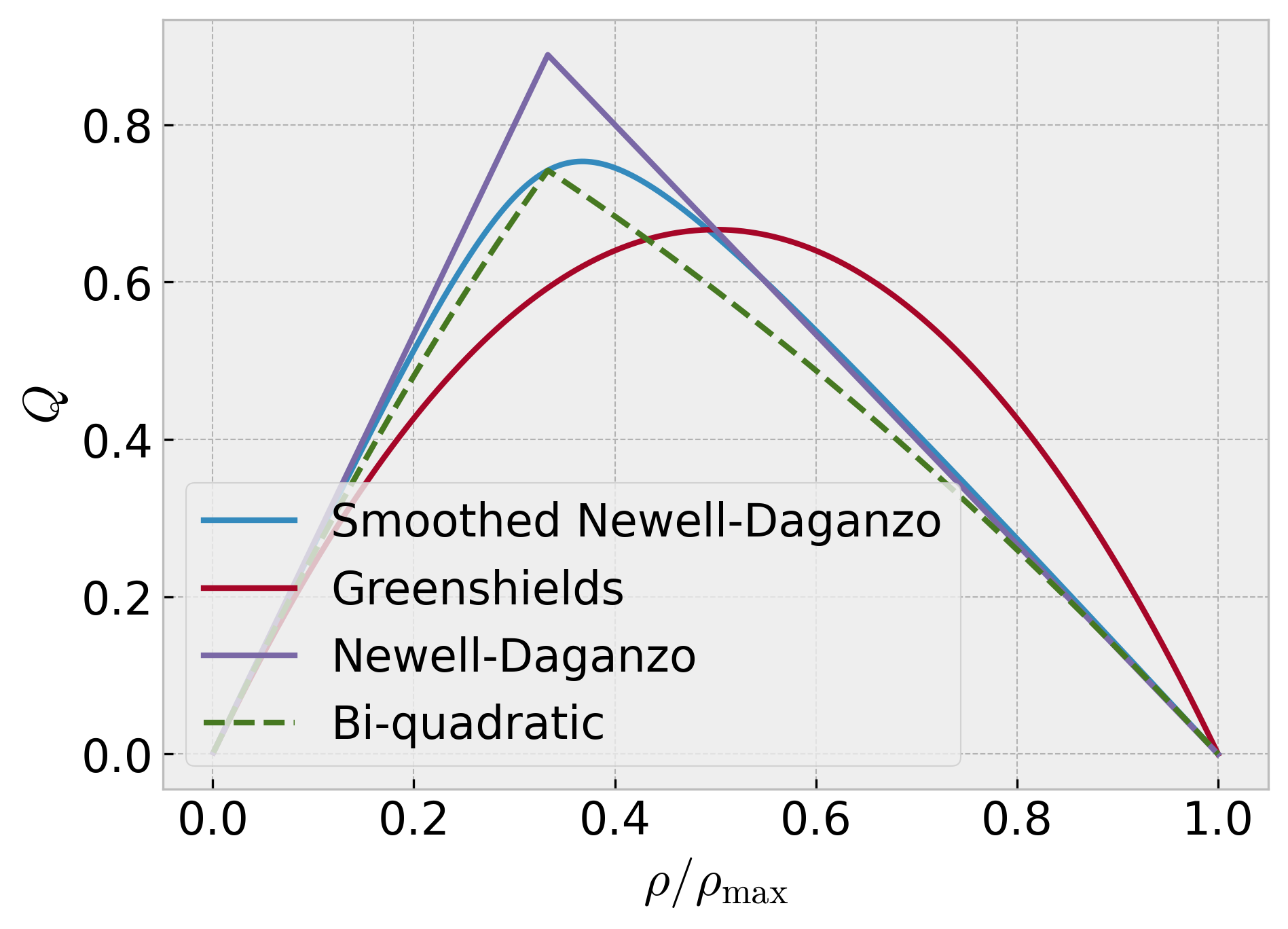}
   \end{tabular} 
   \caption{Fundamental diagram used in \cite{jamitones} and comparison of some flows: (a)  fundamental diagram, (b) comparison between several different flows.}
   \label{img:DFs_aprox}
\end{figure}

The first  fundamental diagram (FD) was measured in the  1930s by Greenshields \cite{primer_DF} and led to the quadratic expression  $Q(\rho) = u_{\max} \rho ( 1- \rho/ \rho_{\max})$, 
corresponding to 
\begin{align} \label{U-lin} 
 U ( \rho) = u_{\max} ( 1- \rho/ \rho_{\max}), 
 \end{align} 
 where $u_{\max}$ is a maximal velocity (corresponding to a free highway) and $\rho_{\max}$ is a maximal density of cars, usualled associated with a bumper-to-bumper situation.  Subsequent studies revealed that a quadratic approximation would not be adequate and numerous alternative  functions have been proposed, such as piecewise linear functions \cite{Newell, Daganzo} or bi-quadratic ones \cite{problema_riemann}. Usually, the fundamental diagram is classified into two or three phases, depending on the theory. The two-phase theory divides the FD into free-flow and synchronized-flow phases, as shown in Figure~\ref{img:DFs}(b). The three-phase theory identifies flow for low density (free flow), medium density (synchronized flow), or high density (moving jams) \cite{df_teoria}. This subdivision into different regimes 
  defines various portions of the function $Q=Q(\rho)$. An example of a function $Q(\rho)$ can be the Newell-Daganzo flow \cite{Newell, Daganzo}, which corresponds to a piecewise linear flow given by
\begin{equation}\label{eq:Newell-Daganzo}
Q(\rho)=
\begin{cases} 
 \dfrac{Q_{\max}\rho}{\rho_{\mathrm{c}}}& \text{if $0 \leq \rho<\rho_{\mathrm{c}}$}, \\
 Q_{\max} \left( \dfrac{\rho_{\max}-\rho}{\rho_{\max}-\rho_{\mathrm{c}}} \right) & \text{if $\rho_{\mathrm{c}} \leq \rho<\rho_{\max}$},
\end{cases}
\end{equation}
where $Q_{\max} = u_{\max}\rho_{\mathrm{c}}$ and $\rho_{\mathrm{c}}$ is the so-called critical  density at which  the flow is maximal. Another example is the function  
\begin{equation}\label{eq:funcion_Q}
    Q(\rho) =  c\left(g(0) + (g(1) - g(0))\frac{\rho}{\rho_{\max}} - g\left(\frac{\rho}{\rho_{\max}}\right)\right) \quad \text{with} \quad   
    g(y) \coloneqq  \sqrt{1+\left(\frac{y-b}{\lambda}\right)^2}  
\end{equation}
proposed in \cite{jamitones}, which is a smoothed version of the Newell-Daganzo flow \eqref{eq:Newell-Daganzo} and whose parameters ($c,b,\lambda$) can be adjusted by using least squares with a dataset to approximate any fundamental diagram \cite{parametros}. In Figure~\ref{img:DFs_aprox}(b)   the  various  flows mentioned 
 are compared, where we set  $\rho_{\mathrm{c}} \coloneqq  \rho_{\max}/3$.

\subsection{LWR model}
The LWR  model (or kinematic traffic model)  \cite{LW,R}     is based on  a fundamental relationship between $\rho$ and $u$. Let $u(\rho) = U(\rho)$, where $U(\rho)$ is defined based on the fundamental diagram. A first choice for $U$ could be
 \eqref{U-lin}. That equation interpolates  linearly between   the maximum speed~$u_{\max}$ 
   at which an individual vehicle can travel and the velocity zero. 
     This choice leads to a quadratic approximation of the fundamental diagram, as seen in Figure~\ref{img:DFs_aprox}(b). However, other alternatives can be defined based on flows $Q(\rho)$ that better approximate the fundamental diagram, with the relationship $U(\rho)=Q(\rho)/\rho$ (remember that $\rho>0$). The LWR model is therefore given by the scalar, nonlinear conservation law 
\begin{equation}\label{eq:LWR}
   \partial_t  \rho + \partial_x Q(\rho)  = 0.
\end{equation}
Although this model  adequately describes  vehicular traffic at low densities, it cannot  reproduce the phenomenon of  ``phantom congestion'' \cite{taco_fantasma}, which consists in  the appearance of traffic jams without the influence of external factors and occurs due to the propagation of small disturbances in traffic. 
 This latter problem is critical because traffic jams arise precisely from this phenomenon, which  is the main focus of this work.

\subsection{Payne-Whitham (PW) model} \label{subsec:PW} 
The Payne-Whitham (PW) model \cite{P,W} was  the first second-order model proposed. It is given by the   system  
\begin{align}\label{eq:PW}
\begin{split}
\partial_t \rho +  \partial_x (\rho u) = 0,\qquad  \partial_t u + u \partial_x u + \dfrac{1}{\rho} \partial_x p(\rho) = \dfrac{U(\rho)-u}{\tau},
\end{split}
\end{align}
where $ p(\rho) $ is a strictly increasing function and (the ``traffic pressure'')   that models preventive driving, and $ \tau $ is the \emph{relaxation time} representing the time it takes for drivers to adjust their speed to the desired speed. Possible choices include  the so-called ``regular pressure'' $p(\rho) \coloneqq  \beta \rho^\gamma$ with coefficients $\beta , \gamma >0$ 
or  the more complicated so-called ``singular pressure'' 
\begin{equation}\label{eq:presion}
    p(\rho) = -\beta\left(\dfrac{\rho}{\rho_{\max}} + \ln\left(1 - \dfrac{\rho}{\rho_{\max}}\right)\right). 
\end{equation}
 In \cite{jamitones}, it is proved that the PW model has solutions governed by self-sustained nonlinear traveling waves (jamitons). However, 
  the PW model has been criticized  in  \cite{Requiem} because, under certain conditions, it admits solutions with negative flows and velocities (vehicles 
   moving in the  direction of decreasing~$x$). For this reason  other second-order traffic models are preferred, although they maintain the relaxation term in the second equation.

\subsection{The inhomogeneous ARZ model}
 The inhomogeneous ARZ model   \cite{resurrection,Z} given by  
    \eqref{eq:ARZ0} modifies the PW model \eqref{eq:PW}. 
     Instead of a pressure function (such as 
    \eqref{eq:presion})  it involves a so-called  hesitation function $h=h(\rho)$, 
      which is assumed to be strictly increasing and plays a role similar  to pressure in the PW  model. 
       The  model \eqref{eq:ARZ0} is a hyperbolic system of conservation laws with a relaxation term. It had originally been proposed in its homogeneous version, but the addition of the relaxation term in the second equation is precisely the ingredient that allows the appearance of jamitons. The functions $U= U(\rho) $ and $ h= h(\rho) $ are assumed to be  twice differentiable and 
        satisfy the following conditions:
\begin{enumerate}[label=\alph*]
    \item[(a)] The function~$U$ satisfies $U'(\rho) < 0$ and $U''(\rho) > 0$,  hence  $ U(\rho)$ is decreasing and $Q(\rho) $ is concave since 
     $Q''(\rho) = 2U'( \rho) + U''(\rho) >0$. \label{cond:cond_1}
    \item[(b)]  The function~$h$ satisfies $h'(\rho)>0$ and $\mathrm{d}^k (\rho h(\rho))/ \mathrm{d} \rho^k >0$, $k=1,2$. \label{cond:cond_2}
   %  \textcolor{red}{It is not clear that strict convexity implies $\mathrm{d} (\rho h(\rho))/ \mathrm{d} \rho >0$. } 
\end{enumerate}
These two assumptions ensure that the equations are well posed in the presence of jump-type solutions. In \cite{resurrection}, it is proven that the system  \eqref{eq:ARZ0} 
 is indeed hyperbolic. This is done by multiplying  first equation of \eqref{eq:ARZ0} by $h'(\rho)$ and subtracting the result from the second, which yields 
  the  system 
\begin{align}\label{eq:ARZ_lineal} \begin{split} 
        \partial_t  \rho  + \partial_x (\rho u) &= 0,  \\
       \partial_t   u_t +  \bigl(u - \rho h'(\rho) \bigr) \partial_x  u  &= \frac{h'(\rho)}{\tau} \bigl(U(\rho) - u \bigr).
    \end{split}
\end{align}
We define
\begin{equation}\label{2p10}
    \boldsymbol{w} \coloneqq   \begin{pmatrix} \rho \\ u \end{pmatrix}, \quad \boldsymbol{L}( \boldsymbol{w}) \coloneqq  \begin{pmatrix} 
        u & \rho \\ 0 & u - \rho h'(\rho)
    \end{pmatrix}, \quad 
    \boldsymbol{f} ( \boldsymbol{w}) \coloneqq    \begin{pmatrix} 0 \\ \dfrac{h'(\rho)}{\tau}(U(\rho) - u) \end{pmatrix} 
\end{equation}
to write the linearized ARZ system \eqref{eq:ARZ_lineal} as
\begin{equation}\label{eq:ARZ_lin}
  \partial_t  \boldsymbol{w}    + \boldsymbol{L}( \boldsymbol{w}) \partial_x  \boldsymbol{w}   = \boldsymbol{f}( \boldsymbol{w}) .
\end{equation}
The eigenvalues of the matrix $  \boldsymbol{L}( \boldsymbol{w})$  are given by
\begin{equation}\label{eq:vel_carac}
    \lambda_1 = u-\rho h'(\rho), \quad \lambda_2 = u,
\end{equation}
which, according to assumption \eqref{cond:cond_2}, satisfy $ \lambda_1 < \lambda_2 $ if $ \rho > 0 $. Therefore, the system \eqref{eq:ARZ0} is strictly hyperbolic. The eigenvalues of the system are known as characteristic speeds and describe conditions for the existence of jamitons. Depending on the variables defined, the ARZ model \eqref{eq:ARZ_lin} may be formulated 
 in two other ways, namely either in conservative or in Lagrangian form.  
 
  The   conservative form is achieved if we define  the variable $y  \coloneqq \rho(u + h(\rho))$,   
   which implies $u = y/\rho - h(\rho)$. Substituting this into the first equation of \eqref{eq:ARZ0} yields 
\[%begin{equation}\label{eq:cons_1}
   \partial_t  \rho +  \partial_x \bigl(y - \rho h(\rho) \bigr)  = 0.
\]%end{equation}
For the second relationship, the first equation of  \eqref{eq:ARZ0} is multiplied by $ (u+h(\rho)) $ and added to the second  equation of  \eqref{eq:ARZ0} multiplied by $ \rho $, 
 which results in the system in Eulerian variables %\textcolor{red}{OK?} 
\begin{align*}
        \partial_t  \rho +  \partial_x  (\rho u)& = 0, \\ 
       \partial_t    \bigl(\rho(u+h(\rho)) \bigr)  +  \partial_x \bigl(u\rho(u+h(\rho)) \bigr)  & =  \frac{\rho}{\tau}\bigl(U(\rho)-u \bigr). 
\end{align*}
Written in conservative variables $\boldsymbol{Q} \coloneqq (\rho, y)^{\mathrm{T}}$, this system  reads  
\begin{align}\label{eq:ARZ_cons}
    \begin{split}
    \partial_t \rho + \partial_x \bigl(y - \rho h(\rho) \bigr) &= 0, \\
\partial_t    y +  \partial_x \left(\dfrac{y^2}{\rho} - yh(\rho)\right) &= \frac{1}{\tau}\bigl(\rho \left(U(\rho)+h(\rho)\right)-y\bigr).
    \end{split}
\end{align}
Defining 
\begin{align} \label{eq:ARZ_cons_ec_prel}
    \boldsymbol{F} (\boldsymbol{Q}) \coloneqq  \begin{pmatrix} y-\rho h(\rho)\\[3pt] \dfrac{y^2}{\rho} - yh(\rho)\end{pmatrix}, \quad 
      \boldsymbol{S} (\boldsymbol{Q}) \coloneqq  \begin{pmatrix}
    0 \\ \dfrac{1}{\tau}\left(\rho \left(U(\rho)+h(\rho)\right)-y\right)
    \end{pmatrix}, 
\end{align} 
we may rewrite the system \eqref{eq:ARZ_cons} as
\begin{equation}\label{eq:ARZ_cons_ec}
   \partial_t  \boldsymbol{Q}  + \partial_x  \boldsymbol{F} (\boldsymbol{Q})  =  \boldsymbol{S} (\boldsymbol{Q}).
\end{equation}
This form  will be essential for  numerical simulations in forthcoming sections, as well as being useful for calculating the characteristic velocities \eqref{eq:vel_carac} through the Jacobian matrix~$\boldsymbol{J}_{\boldsymbol{F}}$ of~$\boldsymbol{F}$. 

To write the ARZ model \eqref{eq:ARZ0} in  Lagrangian form  \cite{lagrangeano_1, lagrangeano_2} we define  the variables $v(\sigma, t)$ and $u(\sigma, t)$, where 
 $v  \coloneqq 1/\rho$ and $\sigma$ satisfies 
 %is the number of vehicles \textcolor{red}{double check that $\sigma$ is correctly defined here. 
 % Shouldn't ``number'' be referred to some time or length interval, or number of vehicles counted form some fixed position?}. These quantities   satisfy  
\[%begin{equation}
     \mathrm{d}\sigma = \rho \, \mathrm{d}x - \rho u\,  \mathrm{d}t,\quad  \sigma(t=0)=\sigma_0 ~  \hbox{given}.
\]%end{equation}
In terms these variables, the model \eqref{eq:ARZ0} can be written as 
\begin{align}\label{eq:ARZ_lagr}
    \begin{split}
        \partial_t  v    &=  \partial_{\sigma} u ,\\
       \partial_t  \bigl(u + \hat{h}(v) \bigr)  &= \dfrac{\hat{U}(v) - u}{\tau},
    \end{split}
\end{align}
where $\hat{h}(v) \coloneqq  h(1/v)$ and $\hat{U}(v)\coloneqq  U(1/v)$   satisfy
\[%begin{equation}\label{eq:hipotesis_lagr}
    \dfrac{ \mathrm{d}\hat{U}}{ \mathrm{d} v} > 0, \quad \dfrac{ \mathrm{d} ^2\hat{U}}{ \mathrm{d} v^2} < 0,\quad \dfrac{ \mathrm{d}\hat{h}}{ \mathrm{d}v} < 0, \quad \text{ and }\quad  \dfrac{ \mathrm{d}^2\hat{h}}{ \mathrm{d}v^2} > 0.
\]%end{equation}
%In Appendix \ref{An:Lag}, the derivation of the system is presented in detail. 
This formulation is useful for  the theoretical construction of jamitons, as well as for  the Smoothed Particle Hydrodynamics (SPH) method (see \cite{SPH-1, SPH-2}) for numerical 
 simulations (as in \cite{presion}, but the SPH method is  not used herein). Both forms have associated Rankine-Hugoniot (jump) conditions. If $[\![ u ]\!] \coloneqq  u^+ - u^-$ 
 denotes  value of the jump of a variable~$u$ across a  discontinuity and $s$ is the jump propagation velocity, 
  then the jump condition for \eqref{eq:ARZ_cons_ec} is  $[\![ \boldsymbol{F}(\boldsymbol{Q}) ]\!] = s [\![ \boldsymbol{Q} ]\!]$, 
   which for Eulerian variables means that  
\begin{align}\label{eq:Rankine-hugoniot-cons}
    \begin{split}
        [\![ y-\rho h(\rho) ]\!]  &= s[\![ \rho ]\!] , \\
         [\![ y^2/\rho - yh(\rho) ]\!]  &= s[\![ y ]\!] 
    \end{split}
\end{align}
or equivalently, if we recall  that $y = \rho(u+h(\rho))$ and use \eqref{eq:Rankine-hugoniot-cons},  
\begin{align}\label{eq:Rankine-hugoniot-eu}
    \begin{split}
        [\![ \rho u ]\!]  &= s [\![ \rho ]\!], \\
        [\![ \rho u^2 + \rho u h(\rho) ]\!] &= s [\![ \rho(u+h(\rho))]\!].
    \end{split}
\end{align}
In  Lagrangian form,  conditions \eqref{eq:Rankine-hugoniot-eu} are replaced by
\begin{align}\label{eq:Rankine-hugoniot-lagr}
    \begin{split}
        [\![  u ]\!] &= m[\![  v ]\!], \\
        [\![  u ]\!] &=- [\![  \hat{h}(v) ]\!],
    \end{split}
\end{align}
where $-m$ is the jump propagation velocity but in Lagrangian variables (in the Eulerian scheme, $m$ represents the traffic flow passing through a jump). In Section 4
 the relation between both quantities $s$ and $m$ is  established.

\subsection{Sub-characteristic condition}

The macroscopic models share a type of solution called a steady state where $\rho(x, t) = \tilde{\rho} = \text{const.}$ and $u(x, t) = \tilde{u} = \text{const.}$. 
 This   means   that the traffic flows uniformly: vehicles always keep the same distance from each other and move at their desired speed $\tilde{u} = U(\tilde{\rho})$. On the other hand, the relaxation term introduced in second-order models causes the traffic velocity $u$ to converge to $U(\rho)$ as $\tau \to 0^+$. In this case, the solutions of 
  the ARZ and PW models are dominated by the continuity equation \eqref{eq:continuidad} with $u=U(\rho)$, which precisely corresponds to the LWR model. Also note that the characteristic velocity 
   of the LWR model is given by 
   \begin{align} \label{murhodef} 
   \mu(\rho) = Q'(\rho) = \rho U'(\rho) + U(\rho). 
   \end{align} 
   Let $\rho(x, t) = \tilde{\rho} \geq0$ and $u(x, t) =\tilde{u}= U(\tilde{\rho})$ denote a base state for some second-order model. We say that the base state $(\tilde{\rho}, \tilde{u})$ is \emph{linearly stable} if any small perturbation decays over time. On the other hand, the base state $(\tilde{\rho}, \tilde{u})$ is said to satisfy the  \emph{sub-characteristic condition} (SCC) if the characteristic velocity $\mu = \mu(\tilde{\rho})$
    defined by \eqref{murhodef}  for $\rho=\tilde{\rho}$  lies between the two characteristic velocities $\lambda_1 = \lambda_1({\tilde{\rho}, \tilde{u}})$, $\lambda_2 = \lambda_2({\tilde{\rho}, \tilde{u}})$ with $\lambda_1<\lambda_2$ of the second-order model, that is:
    \begin{equation}\label{SCC}
        \lambda_1 < \mu < \lambda_2.
    \end{equation}
The previous concepts may be considered different, but actually Whitham's theorem proven in \cite{CSC_demo} relates both:
\begin{thm}[Whitham \cite{CSC_demo}]\label{teo:Whitham}
    The base state $(\tilde{\rho}, \tilde{u})$ is linearly stable if and only if it satisfies the sub-characteristic condition \eqref{SCC}.
\end{thm}

Since second-order models are capable of reproducing traffic instabilities, Whitham's  theorem is crucial for distinguishing between stable and unstable base states more easily through the SCC. Nonetheless, instabilities can be directly studied by calculating the growth factor for perturbations of the form $\mathrm{e}^{\mathrm{i}kx}$, as is done in \cite{stabilidad_jam}, 
 where the  same SCC inequality \eqref{SCC} is obtained.

For the  ARZ model, 
  the characteristic speeds are given by  \eqref{eq:vel_carac}, so its SCC   \eqref{SCC} corresponds to 
\begin{equation*}
    \tilde{u} - \tilde{\rho} h'(\tilde{\rho}) < \tilde{\rho}U'(\tilde{\rho}) + U(\tilde{\rho}) < \tilde{u},
\end{equation*}
but $\tilde{u} = U(\tilde{\rho})$,  so we obtain 
\begin{equation*}
    - \tilde{\rho} h'(\tilde{\rho}) < \tilde{\rho}U'(\tilde{\rho}) < 0.
\end{equation*}
By assumption \eqref{cond:cond_1}, the second inequality always holds. Thus, it follows that for the ARZ model the base state $(\tilde{\rho}, \tilde{u})$ is linearly stable if and only if
$ - \tilde{\rho} h'(\tilde{\rho}) < \tilde{\rho}U'(\tilde{\rho})$ or equivalently, 
\[%begin{equation}\label{eq:CSC_ARZ}
     h'(\tilde{\rho}) + U'(\tilde{\rho}) >0. 
\]%end{equation}

\subsection{Some remarks on Whitham's theorem} 

Theorem \ref{teo:Whitham} has been investigated and extended in several works \cite{whitham_1, whitham_2, whitham_3, CSC_demo, W}, with a primary focus on the SCC and 
  on even  more general $N$-equation hyperbolic systems than traffic models.   These results include the  following findings. 
 We first mention that  Theorem~\ref{teo:Whitham} can be extended to relate the SCC to the stability of solutions, in the sense 
  that a  steady-state solution with small perturbations converges to the initial steady state as $t \to \infty$. Furthermore, 
  solutions to second-order models that satisfy the SCC everywhere converge to LWR solutions as $\tau \to 0$. More  
   generally, if $\tau$ is not small enough, solutions to second-order models converge to LWR solutions but with an additional nonlinear viscosity term of magnitude $\mathcal{O}(\tau)$. 
    Finally, dedicated analyses of the degenerating cases of the SCC (corresponding to the situation when 
     $\lambda_1 = \mu < \lambda_2$ or $\lambda_1 < \mu = \lambda_2$) are available  \cite{deg_1, deg_2}.
     
An interesting question that has been studied in \cite{jamitones, presion, stabilidad_jam} is the following: how do the solutions close to constant base states evolve in second-order models when 
 the SCC is {\em not}  satisfied? In the mentioned works and in numerical studies \cite{presion, lagrangeano_1}, it is shown that these solutions converge to regimes dominated by jamitons, which will be the main focus of this work.

\section{Jamitons}\label{cap:jam}

In this section we describe in detail the construction of 
 jamitons as exact solutions to the ARZ model. We closely follow the approach by Ramadan  \cite{Ramadan}. 

\subsection{Model functions}\label{cap:jam_1}

\begin{figure}[t]
   \centering
   \begin{tabular}{cc} 
   (a) & (b) \\ 
   \includegraphics[width=.45\linewidth]{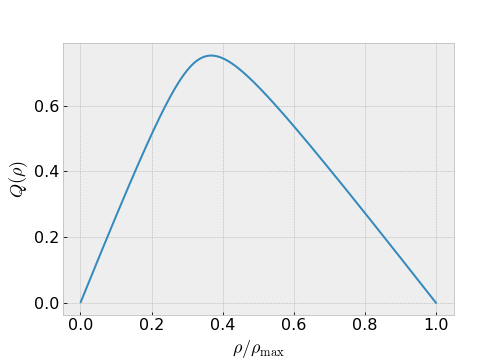} 
   &     \includegraphics[width=.45\linewidth]{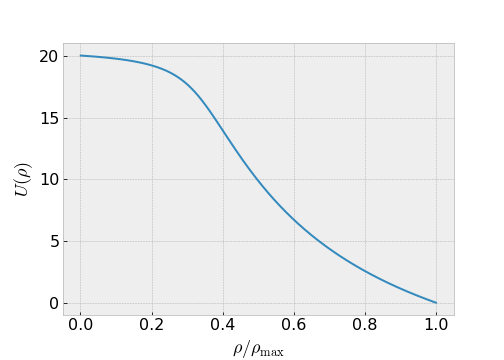} \\
   (c) & (d)   \\ 
   \includegraphics[width=.45\linewidth]{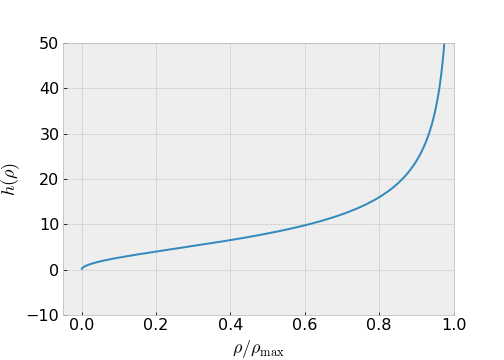} & 
    \includegraphics[width=.45\linewidth]{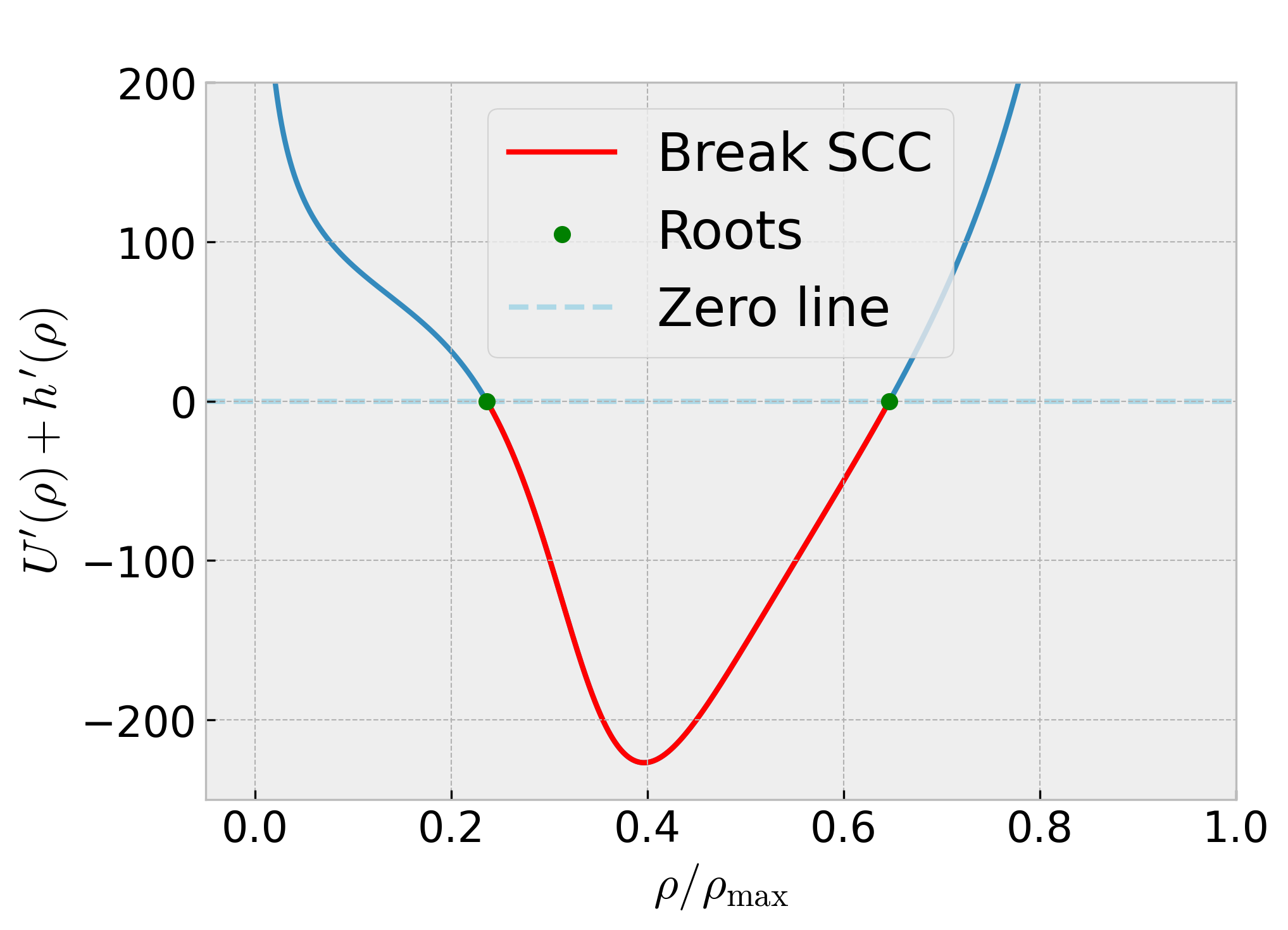}
    \end{tabular} 
   \caption{Plots of functions defined  in section \ref{cap:jam_1}: (a) flow~$Q(\rho)$, (b)   desired velocity~$U(\rho)$, 
    (c) hesitation function~$h(\rho)$, (d) plot of $U'(\rho) + h'(\rho)$ illustrating the interval of violation of the SCC for the corresponding ARZ model.}
   \label{img:funciones}
\end{figure}

The  results of the analysis are  illustrated in a series  of plots  that are based on the following 
specific functions   for the ARZ model \eqref{eq:ARZ0}.  We assume that the function $Q= Q( \rho)$ is given 
 by \eqref{eq:funcion_Q} with  parameters $b=1/3$, $c = 0.078u_{\max}\rho_{\max}$, $\lambda = 1/10$,
   and $u_{\max}=20 \, \mathrm{m} \mathrm{s}^{-1}$, chosen to fit real data \cite{jamitones}. 
The value of $\rho_{\max}$ is based on   assuming a  length of $5 \, \mathrm{m}$   per vehicle plus an extra $50\%$ of safety distance between them. 
 Then, at $\rho_{\max}$ there will be exactly one vehicle each $7.5$ meters, and therefore $\rho_{\max} =(1/7.5) \, \mathrm{m}^{-1}$.  
    The desired velocity is given by $U(\rho) = Q(\rho)/\rho$, 
 and  the hesitation function is  
 \begin{align*} 
 h(\rho) = \beta\left(\frac{\rho}{\rho_{\max} - \rho}\right)^\gamma \quad \text{with $\beta = 8$ and $\gamma = 1/2$}.
 \end{align*} 
 The relaxation time $\tau$ for testing and obtaining general results will be $\tau = 5$ seconds, but in Section~\ref{cap:colision}, effects for other values of~$\tau$ will be studied. 
  The plots of the functions are presented in Figure~\ref{img:funciones}.

\subsection{Traveling-wave analysis}
Jamitons are defined as self-sustained traveling waves that arise in second-order models. These can be constructed following the Zel'dovich-von Neumann-Döring (ZND) theory \cite{ZND} since they have the same mathematical structure as detonation waves. A   traveling-wave  approach leads  to obtain an expression for  jamitons. The construction is based 
  the Lagrangian formulation \eqref{eq:ARZ_lagr}, but the final expressions and simulations will be translated to Eulerian variables for better physical visualization. Let 
  \begin{align*} 
  \chi \coloneqq  \frac{mt + \sigma}{\tau}. 
  \end{align*} 
  The goal is to find traveling-wave solutions $v(\sigma, t) = v(\chi)$ and $u(\sigma, t) = u(\chi)$ of \eqref{eq:ARZ_lagr}. Substituting these expressions into \eqref{eq:ARZ_lagr} yields the system of equations
\begin{align}
    \frac{m}{\tau}v'(\chi) - \frac{1}{\tau}u'(\chi) &= 0,\label{eq:ARZ_ansatz_1}\\
    \frac{m}{\tau}u'(\chi) + \hat{h}'(v(\chi))\frac{m}{\tau}v'(\chi) &= \frac{\hat{U}(v(\chi)) - u(\chi)}{\tau}.\label{eq:ARZ_ansatz_2}
\end{align}
Directly integrating equation \eqref{eq:ARZ_ansatz_1}, one obtains 
\begin{equation}\label{eq:rel_v_u}
    mv - u = -s,
\end{equation}
where $s$ is an integration constant. In Lagrangian variables, $-m$ corresponds to the propagation speed of the jamiton, while $s$ corresponds to the flow of vehicle mass passing through the wave. In Eulerian variables, both constants interchange their meanings ($s$ corresponds to speed and $-m$ to mass flow). Equation \eqref{eq:rel_v_u} implies that 
\begin{equation}\label{eq:u_funcion_v}
    u = s + mv,
\end{equation}
and replacing \eqref{eq:u_funcion_v} in \eqref{eq:ARZ_ansatz_2} yields the ODE of the jamiton:
\begin{equation}\label{eq:EDO_jam}
    v'(\chi) = \dfrac{w(v(\chi))}{r'(v(\chi))}, \quad \text{where $w(\v) = \hat{U}(\v) - (m\v + s)$ and $r(\v) = m\hat{h}(\v) + m^2\v$. }
\end{equation}
Since $\hat{h}'(\v) < 0$ and $\hat{h}''(\v) >0$, the function  $r'(\v)$  has  at most one root $\v_s$, i.e., 
 $r'(v_s) = 0$ or equivalently, $\hat{h}'(\v_s) = -m$. The ODE \eqref{eq:EDO_jam} can then be integrated 
 across  $\v_s$ if $w(\v)$ has one root at $\v_s$ as well, that is $w( v_s) =0$ or equivalently, 
\begin{equation}\label{eq:Chapman_jouguet}
    \hat{U}(\v_s) = m\v_s + s.
\end{equation}
Equation \eqref{eq:Chapman_jouguet} is known as the \emph{Chapman-Jouguet condition}  in ZND theory, and the point $\v_s$ is called the \emph{sonic point}. Then, smooth traveling wave solutions of the ODE \eqref{eq:EDO_jam} can be parametrized by $\v_s$, the sonic point, where the constants $m$ and $s$ are given by 
\begin{equation}\label{eq:m-s}
    m = -\hat{h}'(\v_s) \quad \text{ and } \quad s = \hat{U}(\v_s) - m\v_s.
\end{equation}
Jumps moving at velocity $-m$ can be added to the smooth profile found using the Rankine-Hugoniot conditions \eqref{eq:Rankine-hugoniot-lagr}. The first condition in \eqref{eq:Rankine-hugoniot-lagr}, 
 $[\![ u]\!] = m [\![  v   ]\!]   $,  states that the quantity $m\v - u$ is conserved across a jump, which corresponds to condition \eqref{eq:u_funcion_v} for the smooth part of the traveling wave. Combining both conditions in \eqref{eq:Rankine-hugoniot-lagr} and multiplying 
 the result by~$m$ yields 
\begin{align*} 
 [\![ r (\v) ]\!] =    [\![ m^2\v + \hat{h}(\v) ]\!] = 0;
\end{align*} 
in other words, $r(\v)$ is conserved across a jump. Then, by integrating \eqref{eq:EDO_jam} one can incorporate a jump at some value $\v^-$ that jumps to a value $\v^+$ such that $r(\v^-) = r(\v^+)$ and continue integrating from there. Moreover, a necessary condition for   the jump  to satisfy  the Lax entropy conditions is that   the specific volume must decrease  through the jump, i.e.,  
 $\v^+ < \v_s < \v^-$, so the {\em smooth} portion of    the jamiton profile $\v(\chi)$ must be {\em increasing}. 
 Applying  L'Hospital's rule to the solution of  \eqref{eq:EDO_jam} at the point $\v_s$ and keeping in mind that $\v$ must be increasing, we obtain  that 
\[%begin{equation}
    \dfrac{\hat{U}'(\v_s) + \hat{h}'(\v_s)}{-\hat{h}'(\v_s)\hat{h}''(\v_s)} >0,
\]%end{equation}
that is, $\hat{U}'(\v_s) + \hat{h}'(\v_s) >0$ and therefore the SCC  does not hold at $\v_s$, so jamitons with jumps exist if and only if the SCC is violated. Consequently, since $w''(\v) = \hat{U}''(\v) < 0$ always and $w'(\v) > 0$ for $\v < \v_s$, the function $w$ may have a second root $\v_M$ only for $\v > \v_s$,  %, \textcolor{red}{I find it difficult to follow why  ``the function $w(\v)$ will have a second root $\v_M$ for $\v > \v_s$'', could you insert some intermediate steps?} 
and since \eqref{eq:EDO_jam} cannot be integrated beyond $\v_M$, this profile corresponds to a maximal jamiton connecting the points $\v_M$ with $\v_R$, where $r(\v_R) = r(\v_M)$. It should be noted that the appearance of a maximal jamiton is purely theoretical as it is infinitely long. The same procedure can be carried out in Eulerian variables for \eqref{eq:ARZ_cons} by defining $\eta \coloneqq   (x - st)/\tau$. 
 If we compare this with the kink solution \eqref{phi4_kink} of the $\phi^4$ equation  \eqref{eq:phi_4}, then  the same profile 
  behavior in Eulerian variables is obtained. In the case of jamitons, however, an explicit formula is not available, but the traveling-wave behavior is preserved.

\begin{figure}[t]
   \centering
   \begin{tabular}{cc} (a) & (b) \\ 
   \includegraphics[width=.45\linewidth]{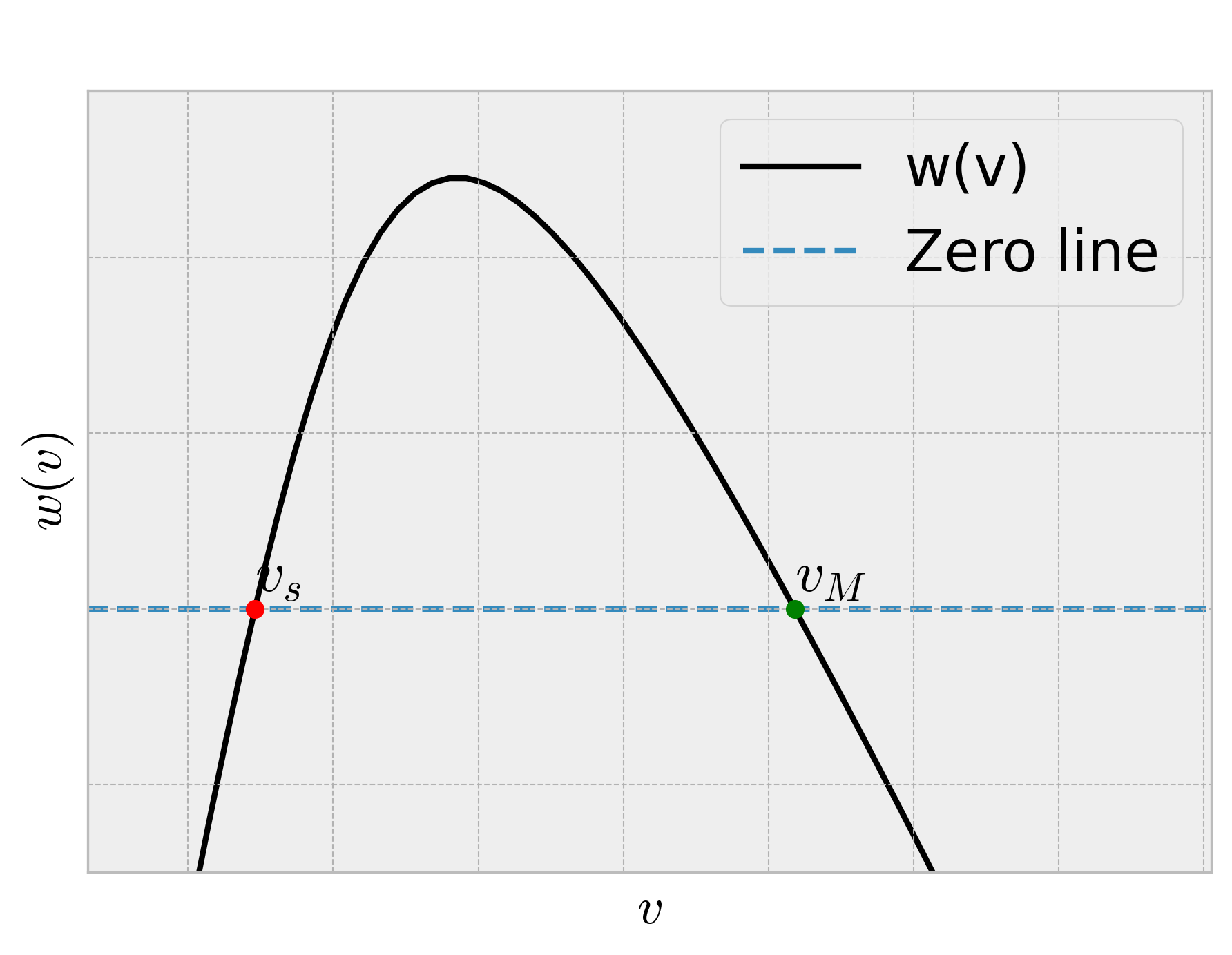} & 
    \includegraphics[width=.45\linewidth]{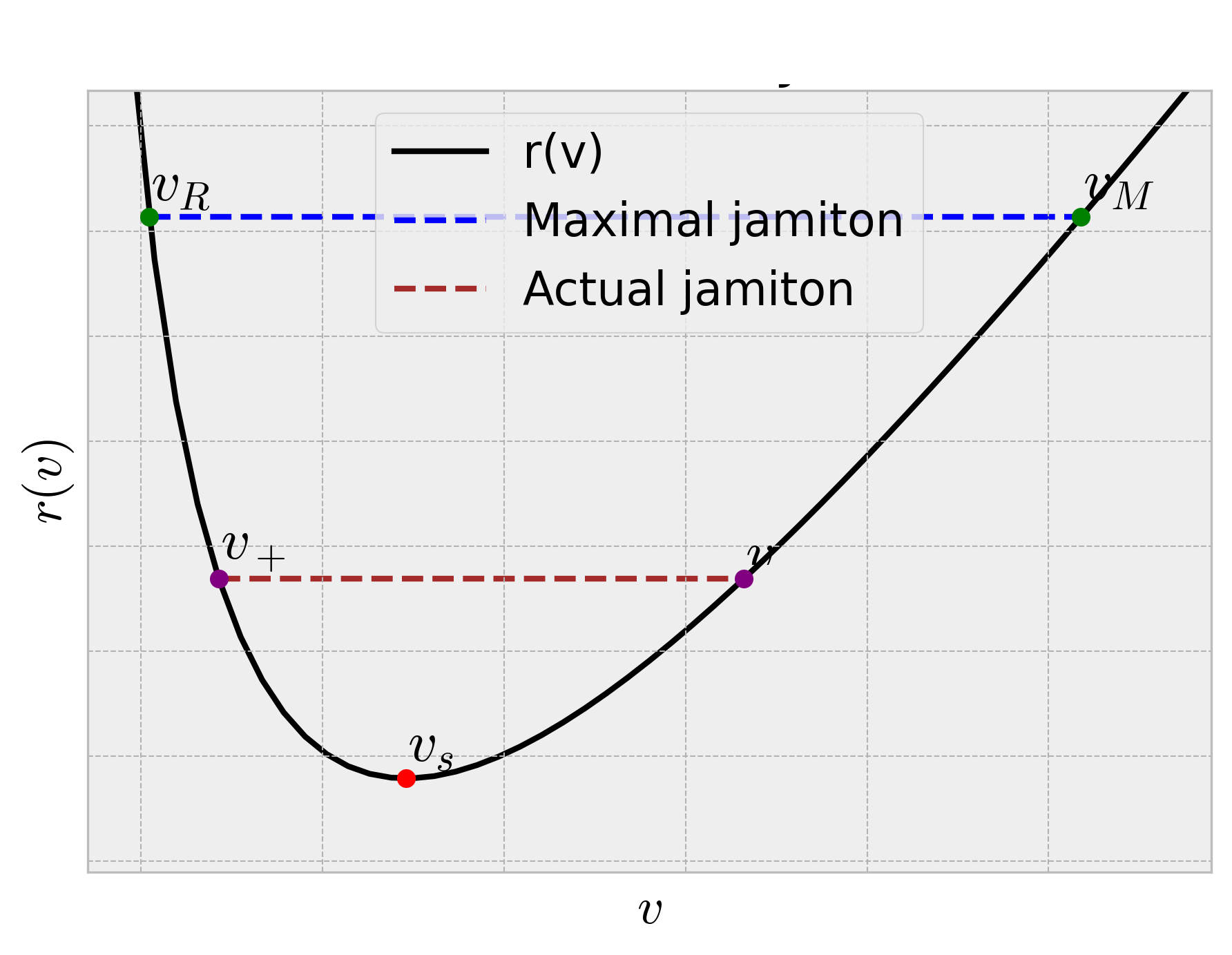}
    \end{tabular} 
   \caption{Functions $w$ and $r$ with a maximal jamiton and an example of standard jamiton: (a)  function $w(\v)$ with a second root $v_M$, 
    (b)   function $r(\v)$ with a maximal jamiton and an example of other jamiton. The point $\v_s$ corresponds to the minimum of~$r$.}
   \label{img:w-y-r}
\end{figure}

\subsection{Construction of jamitons}\label{sec:construccion}
Based on the previously  discussed theoretical aspects, jamitons can be constructed as follows \cite{Ramadan}:
\begin{enumerate}
    \item Choose $\rho_s$ that does not satisfy the SCC condition and define $\v_s \coloneqq  1/\rho_s$.
    \item Define $m \coloneqq  -h'(\v_s)$ and $s  \coloneqq  U(\v_s) - m\v_s$.
    \item Look for $\v_M>\v_s$ such that $w(\v_M)=0$ and define $r_{\min}  \coloneqq r(\v_s)$ and $r_{\max}  \coloneqq r(\v_M)$ (see Figure~\ref{img:w-y-r}).
    \item Search for $\v_R < \v_s$ such that $r(\v_R) = r_{\max}$.
    \item Choose $\v^-$ such that $\v_s < \v^-$ and $r(\v^-) \in (r_{\min}, r_{\max})$.
    \item Choose $\v^+$ such that $\v^+< \v_s$ and $r(\v^-) = r(\v^+)$.
    \item One solves the initial-value problem  of the  ODE \eqref{eq:EDO_jam} from a point $\chi_{\mathrm{start}}$ with initial condition $v(\chi_{\mathrm{start}}) = v^+$ and stops when the value 
     $v^-$ is reached. %Solve the ODE  from $\v^+$ to $\v^-$ to obtain the jamiton equation in the variable~$\chi$. 
    %\textcolor{red}{Is it sufficiently clear what is meant here or should we write that out, e.g. as follows? } 
\end{enumerate}

\begin{figure}[t]
   \centering
   \begin{tabular}{cc} (a) & (b) \\ 
   \includegraphics[width=.47\linewidth]{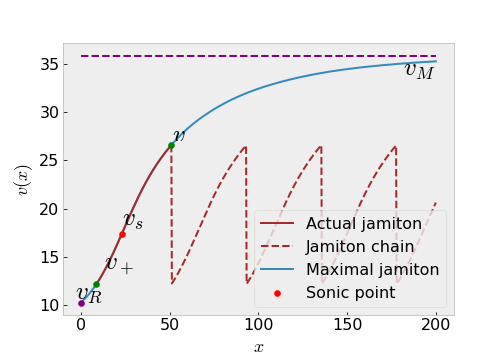} & % requires the graphicx package
    \includegraphics[width=.47\linewidth]{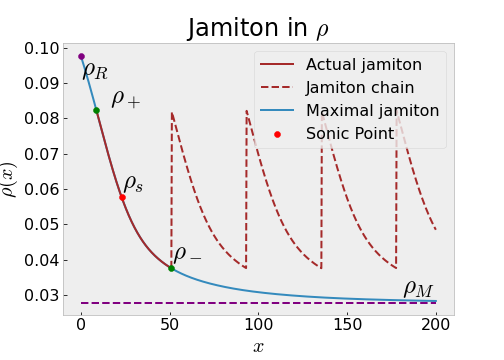} 
    \end{tabular} 
    (c) \\
    \includegraphics[width=.47\linewidth]{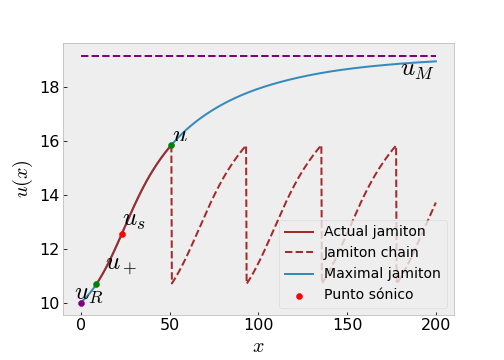}
   \caption{Example of a chain of jamitons, along with a maximal jamiton, 
    shown as a function  (a)  $v=v(x)$, with jamitons between $\v^+$ and $\v^-$, 
    (b)  $\rho = \rho(x)$, with jamitons between $\rho^+$ and $\rho^-$, 
     (c) $u = u(x)$, with jamitons between $u^+$ and $u^-$.}
   \label{img:varios_jamitones}
\end{figure}

\begin{figure}[t]
   \centering
   (a)  \\ 
   \includegraphics[width=\linewidth]{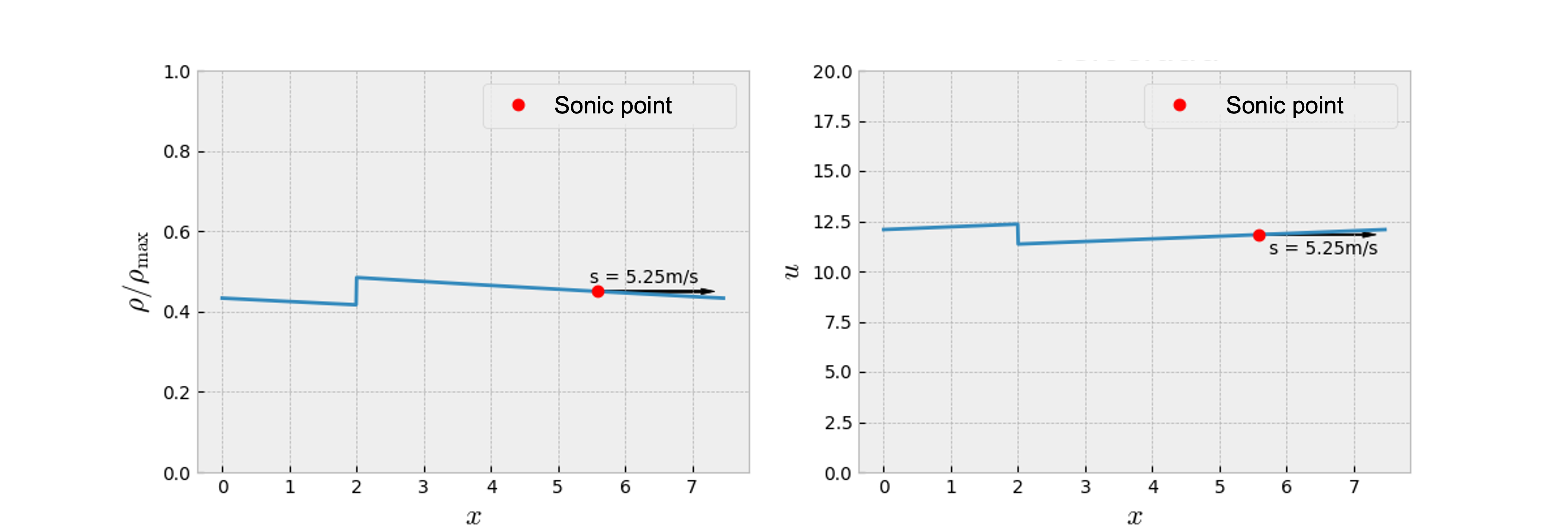} \\
   (b) \\ 
   \includegraphics[width=\linewidth]{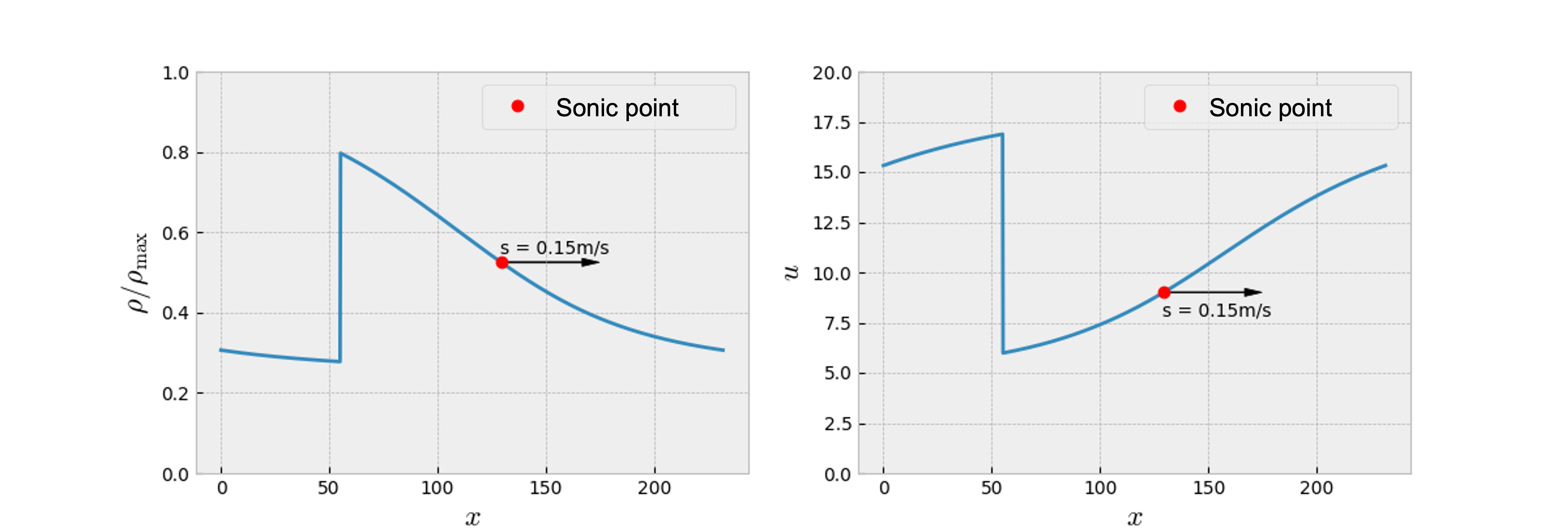} \\
   (c) \\ 
    \includegraphics[width=\linewidth]{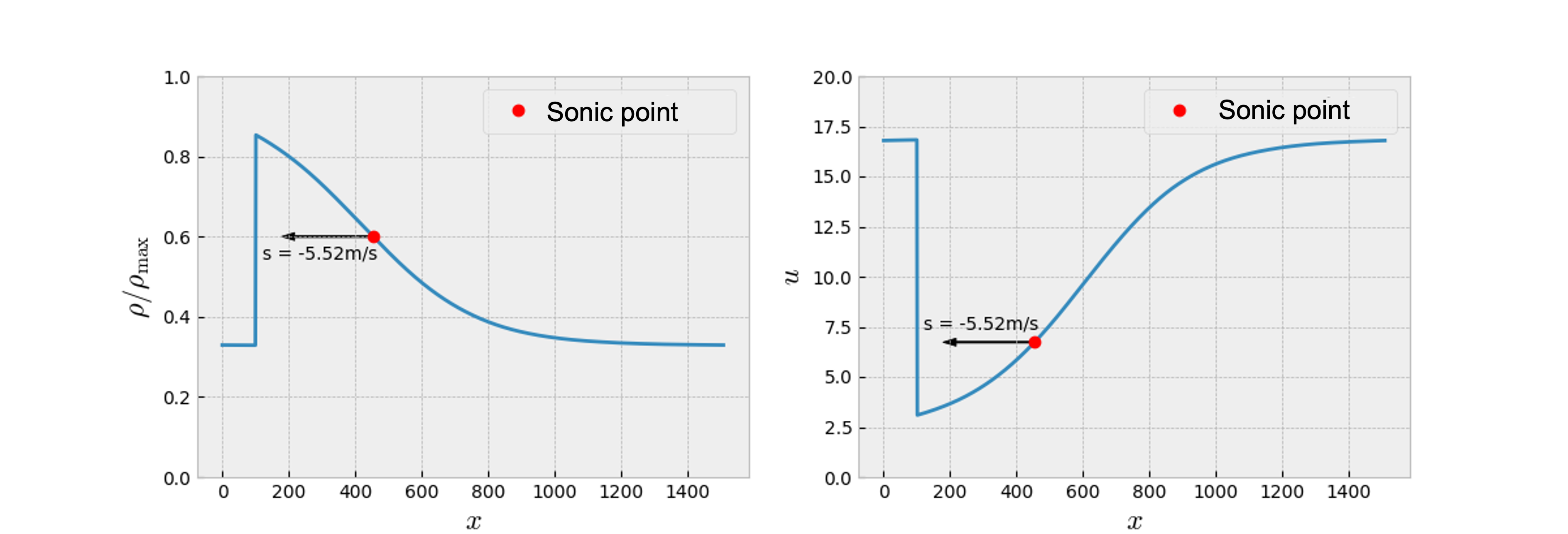}  
       \caption{Jamitons of various lengths: (a) small-length jamiton with $\rho_s = 0.45\rho_{\max}$, with length $L\approx 15.27$ and 
       $N\approx 0.91$ vehicles, (b) medium-length jamiton with  $\rho_s = 0.525\rho_{\max}$,   $L\approx 231.70$ and $N\approx 14.05$
        (notice that this jamiton moves to the left), 
        (c) long  jamiton with $\rho_s = 0.6\rho_{\max}$,   $L\approx 1507.27$ and $N\approx 93.65$, 
        shown in all cases in the variables (left) $\rho$ and (right) $u$.}
   \label{img:jam_grande}
\end{figure}

In Figure~\ref{img:varios_jamitones}, examples of jamitons in the variables~$v$, $\rho$, and~$u$ are presented, where we 
 recall that~$\v = 1/\rho$ and $u$~is obtained from \eqref{eq:u_funcion_v}. In such cases, integration was performed in the variable~$x$, 
  where we used that 
\begin{align*} 
    \mathrm{d} \eta = \frac{1}{\tau}(x - s\, \mathrm{d} t) = \frac{1}{\tau}(\mathrm{d}x + m\v\,\mathrm{d}t -u\, \mathrm{d}t) 
     =\frac{1}{\tau}(\v\,  \mathrm{d}\sigma + m\v\, \mathrm{d}t) = \v\, \mathrm{d} \chi,
\end{align*} 
and the chain rule in the ODE \eqref{eq:EDO_jam}. With this construction, some properties of the jamiton in the coordinate~$x$ can be calculated, such as 
 its length~$L$ given by 
        \begin{equation}\label{eq:Largo}
            L = \tau \int_{\v^+}^{\v^-}\v\dfrac{r'(\v)}{w(\v)} \, \mathrm{d}\v, 
        \end{equation}
  the  total number~$N$ of vehicles \cite{jamitones} in the jamiton given by
 %  \textcolor{red}{Please double check this formula; maybe a step of derivation is appropriate.} 
    \begin{equation}\label{eq:Vehiculos_total}
        N = \tau \int_{\v^+}^{\v^-}\dfrac{r'(\v)}{w(\v)} \, \mathrm{d} \v, 
    \end{equation} 
 and its   amplitude $A$ defined by 
        \begin{equation}\label{eq:Amplitud}
            A = \rho_+ - \rho_-.
        \end{equation} 
In Figure~\ref{img:jam_grande}  jamitons of various  sizes and speeds are shown, along with their lengths \eqref{eq:Largo} and total numbers of vehicles \eqref{eq:Vehiculos_total}. In this case, the number of vehicles may not be natural since it is a continuous model, but an approximate value is obtained.

\subsection{Jamitons in the fundamental diagram}

\begin{figure}[t]
   \centering
   \begin{tabular}{cc} (a) & (b) \\ 
   \includegraphics[width=.45\linewidth]{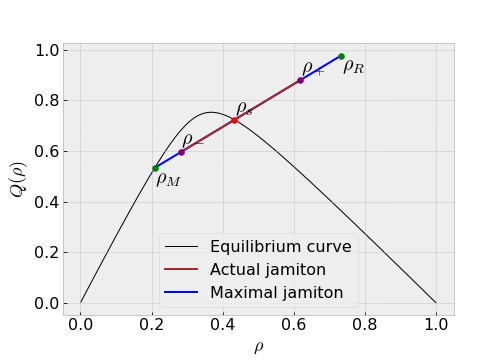} & 
   \includegraphics[width=.45\linewidth]{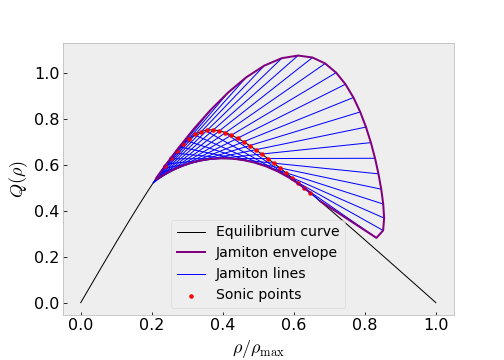} 
   \end{tabular} 
   \caption{Illustration of jamitons in the fundamental diagram: (a) maximal jamiton and a jamiton in the $(\rho,q)$-plane (both lie on the same line passing through $(\rho_s, Q(\rho_s))$), (b) jamiton region composed of several maximal jamitons, along with the jamitonic envelope.}
   \label{img:jam_seg}
\end{figure}

It is instructive to visualize jamitons in the fundamental diagram. Roughly speaking, the idea consists in visualizing 
all pairs $(\rho, q \coloneqq \rho u)$ that arise within the definition  of one jamiton. Clearly, since jamitons are not constant solutions, 
 they cannot be points on the FD. By multiplying the relationship \eqref{eq:u_funcion_v} by~$\rho$, one has
\begin{equation}\label{eq:jam_segmento}
    q = \rho u =  s\rho + m,
\end{equation}
in other words, a jamiton gives rise to a  line segment in the $(\rho, q)$-plane of flow rate versus density, and its slope corresponds to the propagation speed of the jamiton. Furthermore, using the Chapman-Jouguet condition \eqref{eq:Chapman_jouguet} for $\v_s = 1/\rho_s$ and the equalities \eqref{eq:m-s}, we may rewrite \eqref{eq:jam_segmento} as
\begin{align}  \label{eq:jam_segmento2}
    q = s(\rho - \rho_s) + \rho_s U(\rho_s) 
\end{align} 
and evaluating \eqref{eq:jam_segmento2} at $\rho = \rho_s$, we obtain  that $q(\rho_s) = \rho_s U(\rho_s) = Q(\rho_s)$ by definition of $Q$, and therefore, the jamiton segment intersects the equilibrium flow at $\rho_s$. In Figure \ref{img:jam_seg}, an example of a jamiton in the fundamental diagram is depicted, where the endpoints of the segment $(\rho^+, \rho^+u^+)$ and $(\rho^-, \rho^-u^-)$ represent the two states of the jamiton across the jump.

For each $\rho_s$ that violates the SCC, the corresponding maximal jamiton is the line segment connecting the points $(\rho_M, m + s\rho_M)$ and $(\rho_R, m + s\rho_R)$,  where $\rho_M = 1/\v_M$ and $\rho_R = 1/\v_R$ are constructed following the steps in Section~\ref{sec:construccion}. 
 Since $w(\v_M) = 0$ by construction, there holds  $\smash{\hat{U}(\v_M) = m\v_M + s}$,  and multiplying this identity 
  by~$\rho_M$ yields $  m + s\rho_M= Q(\rho_M)$, hence  $(\rho_M, m + s\rho_M)$ also belongs to the equilibrium curve.  Figure~\ref{img:jam_seg} shows some maximal jamitons for various $\rho_s$ that violate  the SCC, along with the jamitonic envelope. This envelope indicates the boundary on which the jamitons must fall, called the jamiton region. The upper envelope corresponds to the curve generated by joining the points $(\rho_R, m+s\rho_R)$. On the other hand, the segments of jamitons also intersect below the equilibrium line, so they also form a lower envelope. This curve is given by points $(\rho^*, q^*)$ that satisfy
\[%begin{equation}
    \rho^*(\rho_s) = -\frac{m'(\rho_s)}{s'(\rho_s)}, \quad  q^*(\rho_s) = m(\rho_s) + s(\rho_s) \rho^*(\rho_s).
\]%end{equation}
An interesting result proven in \cite{jamitones, stabilidad_jam} is that the velocity of a jamiton is determined  by the sonic density~$\rho_s$ only. Furthermore, $s$~decreases with respect to~$\rho_s$. Indeed, from \eqref{eq:m-s} written in Eulerian variables, 
 we deduce that 
\[%begin{equation}
    s(\rho_s) = U(\rho_s) - \rho_s h'(\rho_s), 
\]%end{equation}
hence  
\[%begin{equation}
    s'(\rho_s) = U'(\rho_s) - h'(\rho_s) - \rho_s h''(\rho_s) = U'(\rho_s) + h'(\rho_s) - \bigl(\rho_s h(\rho_s) \bigr)''. 
\]%end{equation}
We have $U'(\rho_s) + h'(\rho_s) <0$  since it is assumed  that $ \rho_s $ violates  the SCC (if not, there would be no jamiton), while  $ (\rho_s h(\rho_s))''>0 $ since it was assumed that $ \rho h(\rho) $ is convex. Thus, $ s'(\rho_s) < 0 $, i.e., $ s $  decreases  with respect to~$\rho_s$.   Figure~\ref{img:jam_seg} illustrates how the slope of each maximal jamiton segment decreases as $\rho_s$~increases.

\subsection{Stability of jamitons}

\begin{figure}[t]
   \centering
   \begin{tabular}{cc} (a) & (b)  \\ 
   \includegraphics[width=.47\linewidth]{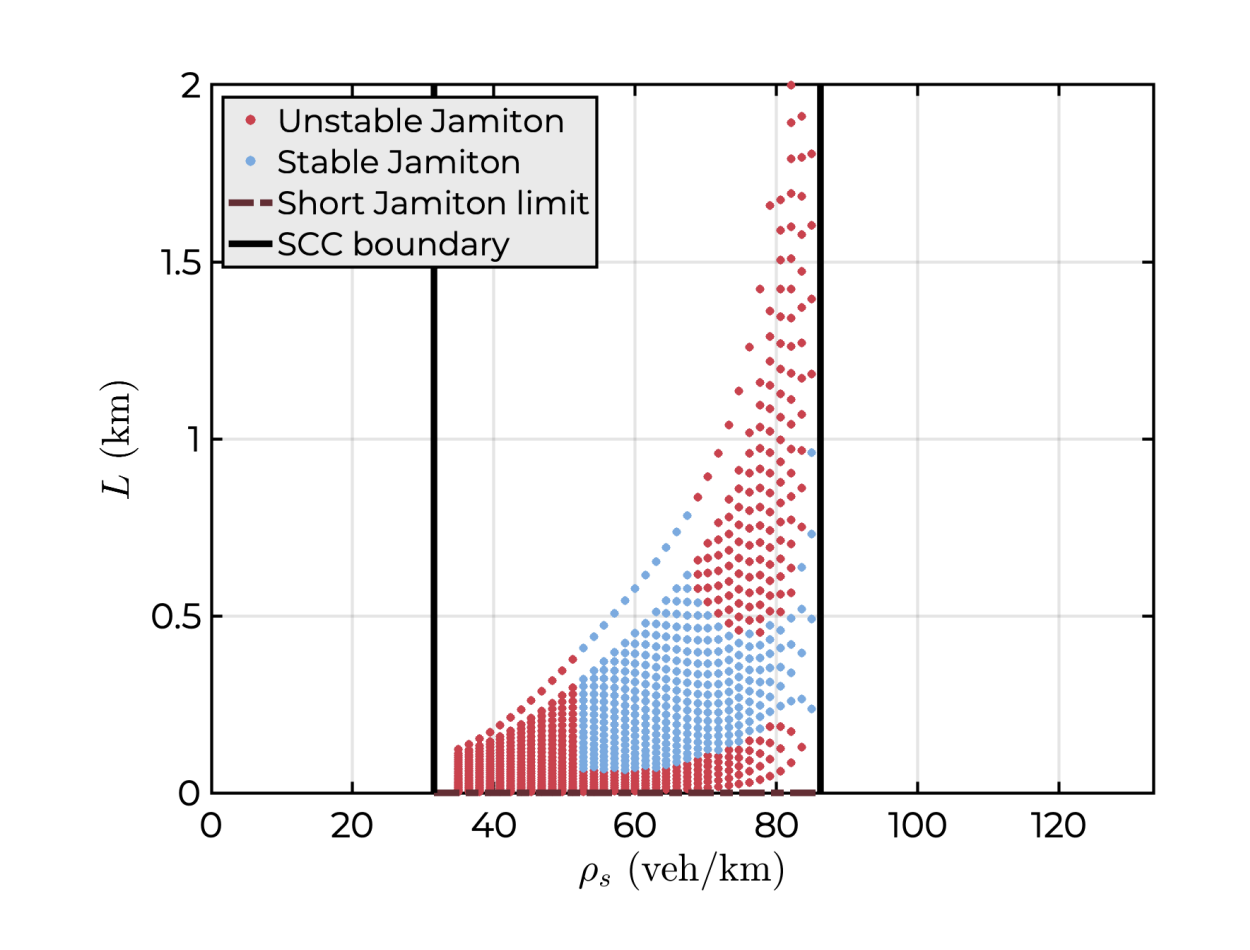} & 
    \includegraphics[width=.47\linewidth]{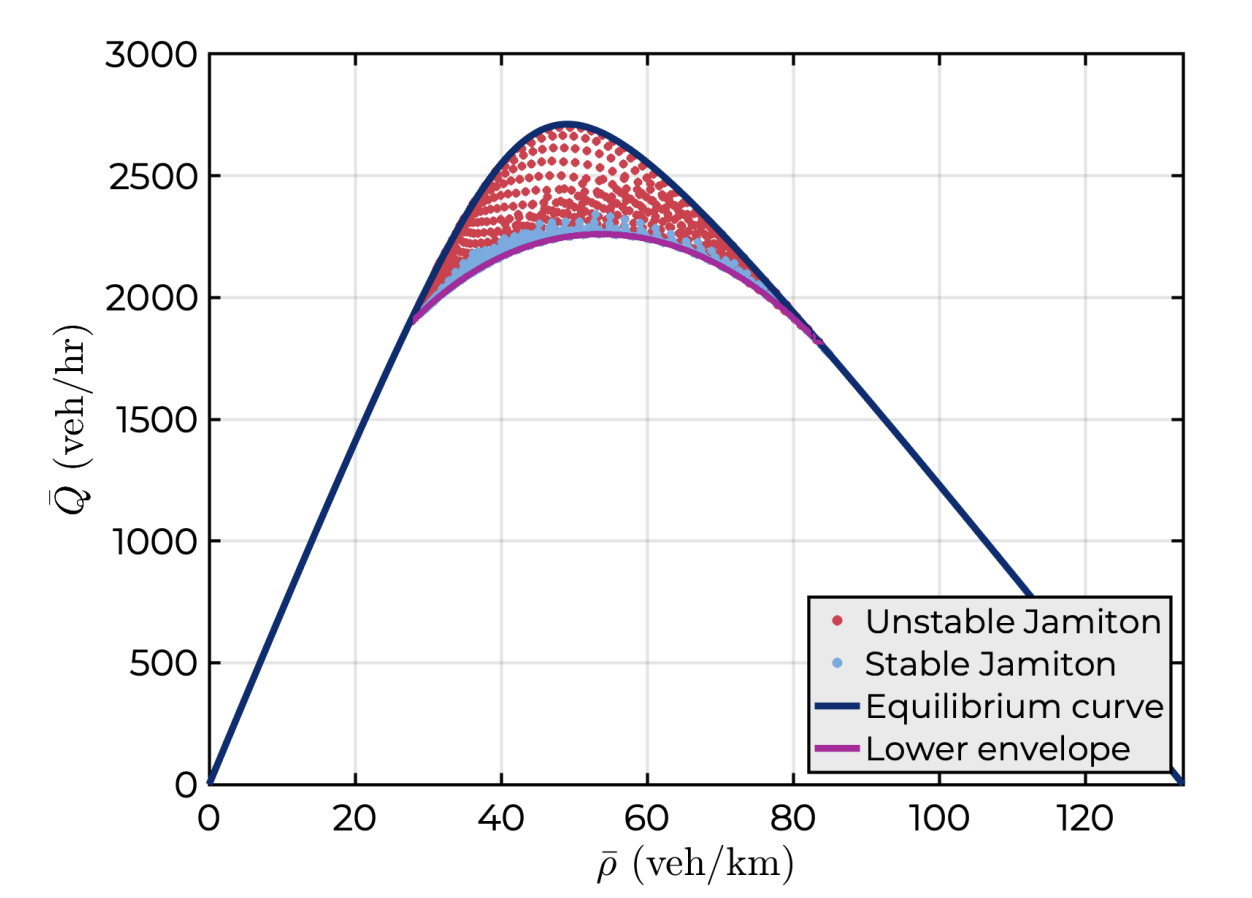}
    \end{tabular} 
   \caption{Stability plots obtained from \cite{stabilidad_jam} Left: Stable and unstable jamitons in the $(\rho_s, L)$ plane. Right: Stable and unstable jamitons in the fundamental diagram.}
   \label{img:estabilidad}
\end{figure}

In \cite{stabilidad_jam}, the dynamical stability of jamitons under small perturbations was  studied numerically. In this case, it was determined that stability depends on the size of the jamiton: very small jamitons are unstable because they collide with new jamitons arising from the perturbations, while larger jamitons are also unstable, but in this case, they split into two or more jamitons. Both cases are visible in the phase plane $ (\rho_s, L) $, as well as in the fundamental diagram, as shown in Figure \ref{img:estabilidad}.
Additionally, in \cite{Ramadan} both concepts of (orbital) stability and asymptotic stability of jamitons are also discussed. While the notion of stability is clearly obtained as the numerical survival of perturbations of jamitons (see also the previously cited work \cite{stabilidad_jam} for a deeper analysis), the notion of asymptotic stability is defined as the existence of exponentially decaying (in time) solutions for the corresponding linearized dynamics, also called sometimes in dispersive PDEs as ``convective stability''. For further details, see \cite[Chapter 5]{Ramadan}.

\section{Simulation of the ARZ model}\label{cap:ARZ}

The numerical method  for the discretization of  the ARZ model  \eqref{eq:ARZ0}
  involves 
 the well-known Harten-Lax-van Leer (HLL) \cite{hllpaper} numerical flux, which represents an  approximate Riemann solver. Since the ARZ model includes a relaxation term that makes it non-homogeneous, it is approximated by a  two-step time-splitting 
   scheme (cf.\  \cite[Chapter~17]{volumen_finito}): first, the homogeneous model is solved  by using 
   a  finite volume scheme (involving the HLL numerical flux), and then, the  differential equation 
\begin{equation}\label{eq:EDO_relajacion}
   \partial_t  \boldsymbol{Q}  = \boldsymbol{S} (\boldsymbol{Q}) 
\end{equation}
is solved, 
where $\boldsymbol{S}$ corresponds to the  relaxation source term in the system. In each time step, 
   the numerical solution of the ODE  \eqref{eq:EDO_relajacion} is applied 
 to  the finite volume solution obtained in the previous step. 
 
 %A  major challenge for implementation  arises when the source term   is, for example, nonlinear or involves higher-order derivatives (diffusion or dispersion terms).

%\subsection{Homogenous model simulation}
%The homogeneous ARZ model corresponds to the system \eqref{eq:ARZ} without the relaxation term:
%\begin{align}\label{eq:ARZ_homogeneo}
%\begin{split}
%    \rho_t + (\rho u)_x &= 0,\\ 
%(u + h(\rho))_t + u(u + h(\rho))_x &= 0.
%\end{split}
%\end{align}
%and written in its conservative form \eqref{eq:ARZ_cons_ec}, \eqref{eq:ARZ_homogeneo} reads
%\begin{equation}\label{eq_Q}
%     \boldsymbol{Q}_t + F(\boldsymbol{Q})_x = 0.
%\end{equation}
%Thus, the finite volume scheme
%\[%begin{equation}\label{eq:esquema_VF}
%    \boldsymbol{Q}_i^{n+1} = \boldsymbol{Q}_i^n - \dfrac{\Delta t}{\Delta x}(\bF_{i+1/2}^n - \bF_{i-1/2}^n)    
%\]%end{equation}
% can be applied with the Godunov flux
%\[%begin{equation}\label{eq:flujo_godunov}
%    \F(\boldsymbol{Q}_{L}, \boldsymbol{Q}_R) = \bF(q^{\downarrow}(\boldsymbol{Q}_{L}, \boldsymbol{Q}_R)) 
%\]%end{equation}
%where $ q^{\downarrow}(Q_l, Q_r) $ denotes the solution to the Riemann problem between the states $ Q_l $ and $ Q_r $. In general, solving the Riemann problem analytically is computationally expensive and leads to undesired behaviors. Therefore, an approximator will be used.

\subsection{HLL solver for a system of  conservation laws}%\label{sec:HLL}

A simple approximator that yields good results and ensures that the numerical flux satisfies the entropy condition \cite{entropia} is the HLL numerical scheme \cite{stabilidad_jam}. To introduce the scheme, we consider the system of conservation laws 
\begin{equation}\label{eq:ARZ_cons_ec_hom}
   \partial_t  \boldsymbol{Q}  + \partial_x  \boldsymbol{F} (\boldsymbol{Q})  =  \boldsymbol{0} 
\end{equation}
  and assume a uniform grid of cells 
$ I_{i,n} \coloneqq [x_{i-1/2}, x_{i+1/2}] \times [t_n, t_{n+1}]$, where it is assumed that 
 $x_{i-1/2} = (i-1/2) \Delta x$ and $t_n = n \Delta t$ for $i \in \mathbb{Z}$ and $n \in  \mathbb{N}_0$. 
  If $\boldsymbol{Q}_i^n$~denotes an approximate value of the cell average of~$\boldsymbol{Q}$ on~$I_{i,n}$ then 
   an explicit   finite volume scheme for \eqref{eq:ARZ_cons_ec_hom} can be written  as the marching formula 
   \begin{align} \label{marching}  
    \boldsymbol{Q}_i^{n+1}   = \boldsymbol{Q}_i^{n} - (\Delta t / \Delta x) \bigl( 
     \boldsymbol{\mathcal{F}}  (  \boldsymbol{Q}_i^{n}, \boldsymbol{Q}_{i+1}^{n} )- 
       \boldsymbol{\mathcal{F}}   (  \boldsymbol{Q}_{i-1}^{n}, \boldsymbol{Q}_{i}^{n} ) \bigr), 
        \end{align} 
        where $  \boldsymbol{\mathcal{F}}  =  \boldsymbol{\mathcal{F}}  (  \boldsymbol{Q}_L, \boldsymbol{Q}_R )$ is a 
         numerical  flux function that among other properties  should be consistent with the 
          exact flux $\boldsymbol{F} = \boldsymbol{F} (\boldsymbol{Q}) $ in the sense that 
           $  \boldsymbol{\mathcal{F}}  (  \boldsymbol{Q}, \boldsymbol{Q} ) = \boldsymbol{F} (\boldsymbol{Q}) $ for all~$\boldsymbol{Q}$. 
       
The HLL scheme  is defined by the particular choice of  $  \boldsymbol{\mathcal{F}}  =  \boldsymbol{\mathcal{F}}  (  \boldsymbol{Q}_L, \boldsymbol{Q}_R )$ 
 given by 
 \begin{align} \label{hllflux} 
 \boldsymbol{\mathcal{F}}^{\mathrm{HLL}}   (  \boldsymbol{Q}_L, \boldsymbol{Q}_R ) 
 = \frac{1}{s_R^+ - s_L^-} \bigl(  s_R^+ \boldsymbol{F}(\boldsymbol{Q}_L) - s_L^- \boldsymbol{F} (\boldsymbol{Q}_R) + s_R^+ s_L^- 
  (\boldsymbol{Q}_R - \boldsymbol{Q}_L) \bigr), 
 \end{align} 
 where $s_L$ and $s_R$ are lower and upper estimates of the characteristic speeds involved in the solution of the Riemann 
  problem for \eqref{eq:ARZ_cons_ec_hom} defined for left and right states~$\boldsymbol{Q}_L$ and~$\boldsymbol{Q}_R$,  and we define 
  $a^+ \coloneqq \max \{ a, 0\}$ and $a^- \coloneqq \min \{ a, 0 \}$.   For a hyperbolic system   \eqref{eq:ARZ_cons_ec_hom} 
   of a general number~$m$ of equations and  unknowns, 
   where the characteristic speeds (eigenvalues of the Jacobian matrix $\boldsymbol{J}_{\boldsymbol{F}} ( \boldsymbol{Q}) 
    = ( \partial F_i ( \boldsymbol{Q}) / \partial \phi_j)_{1 \leq i,j \leq m }$)  are given by $\lambda_1 (  \boldsymbol{Q}) 
     \leq   \lambda_2 (  \boldsymbol{Q})  \leq \cdots \leq  \lambda_m (  \boldsymbol{Q})$,  
      one may choose $s_R = \max \{ \lambda_m  (  \boldsymbol{Q}_L), \lambda_m  (  \boldsymbol{Q}_R) \}$ 
       and  $s_L = \min \{ \lambda_1  (  \boldsymbol{Q}_L), \lambda_1  (  \boldsymbol{Q}_R) \}$.   
       These formulas are utilized for $m=2$ when \eqref{eq:ARZ_cons_ec_hom} represents the 
       the homogeneous version 
of the system of balance equations \eqref{eq:ARZ_cons_ec},  and the corresponding characteristic speeds~$\lambda_1$ and~$\lambda_2$ are given
  by~\eqref{eq:vel_carac}.

 The preceding description is limited to the final form of the HLL scheme. For its motivation based on the simplified solution
  of the above-mentioned Riemann problem through a wave configuration that consists of just two waves separating three 
  states (namely $\boldsymbol{Q}_L$, an intermediate one, and~$\boldsymbol{Q}_R$) we refer to monographs on numerical 
   schemes for conservation laws, for instance  \cite{toro2009,volumen_finito,hestbook,kuzminbook}.

\subsection{Simulation with relaxation term}\label{sec:simulacion_arz}
To incorporate the source term into the model, one simply needs to numerically solve the ODE \eqref{eq:EDO_relajacion}, where the update accounts for the term $\boldsymbol{Q}_i^*$ obtained from simulating the homogeneous scheme. The non-homogeneous ARZ model written in conservative variables 
 reads \eqref{eq:ARZ_cons_ec}, where the source term $\boldsymbol{S} ( \boldsymbol{Q})$ is given by the second equation of 
  \eqref{eq:ARZ_cons_ec_prel}. Then one must solve \eqref{eq:EDO_relajacion} 
and since the relaxation term only appears in the second equation, it suffices to update  the variable~$y$. Equation~\eqref{eq:EDO_relajacion}   is solved using an implicit time-stepping scheme, as it is undesirable to impose additional CFL conditions on spatial and temporal steps, and to take advantage of the absence of a relaxation term for the variable $\rho$, simplifying the computational implementation. 
 Discretizing \eqref{eq:EDO_relajacion}   gives that
\begin{equation}\label{eq:imp_relax}
    \frac{y_i^{n+1} - y_i^*}{\Delta t} = \frac{\rho_i^{n+1}U(\rho_i^{n+1}) + h(\rho_i^{n+1}) - y_i^{n+1}}{\tau}, \qquad \rho_i^{n+1} = \rho_i^*,
\end{equation}
where the second equation holds because there is no source term for $\rho$. Thus, defining $\alpha \coloneqq  \Delta t/ \tau$, 
  we  obtain from  \eqref{eq:imp_relax} the following update scheme for the variable~$y_i^{n+1}$:
\[%begin{equation}
    y_i^{n+1} = \frac{\alpha}{\alpha + 1}(\rho_i^*U(\rho_i^*) + h(\rho_i^*)) + \frac{1}{\alpha + 1}y_i^*,
\]%end{equation}
where $\rho_i^*$ and $y_i^*$ are obtained from \eqref{marching} with $\boldsymbol{Q}_i^{n+1}$ replaced by 
 $\boldsymbol{Q}_i^*=(\rho_i^*, y_i^*)^{\mathrm{T}}$.  

\subsection{Results}
In \cite{MIT_jam}, several simulations of traffic jams on circular roads and the appearance of initially small jamitons, also known as ``jamitinos,''  on an infinite road are presented. To validate the scheme, the results obtained in \cite{MIT_jam} will be emulated, and error tables between the simulation with a jamiton as the initial condition and its theoretical solution will also be presented. In all simulations, periodic boundary conditions (circular road) were considered. The benefit of this type of conditions is the preservation of vehicular mass over time since there are no exits or entrances of vehicles, meaning that the number of vehicles $ N $ given by
\begin{align*} 
    N(t)  \coloneqq  \int_{0}^L \rho(x, t) \, \mathrm{d}x,
\end{align*}  
is constant for all times~$t$.  This provides another way to validate the numerical simulations.

\subsection{Example~1: accuracy of the numerical scheme}

\begin{figure}[t]
   \centering
\begin{tabular}{cc} 
(a) & (b)  \\
   \includegraphics[width=.48\linewidth]{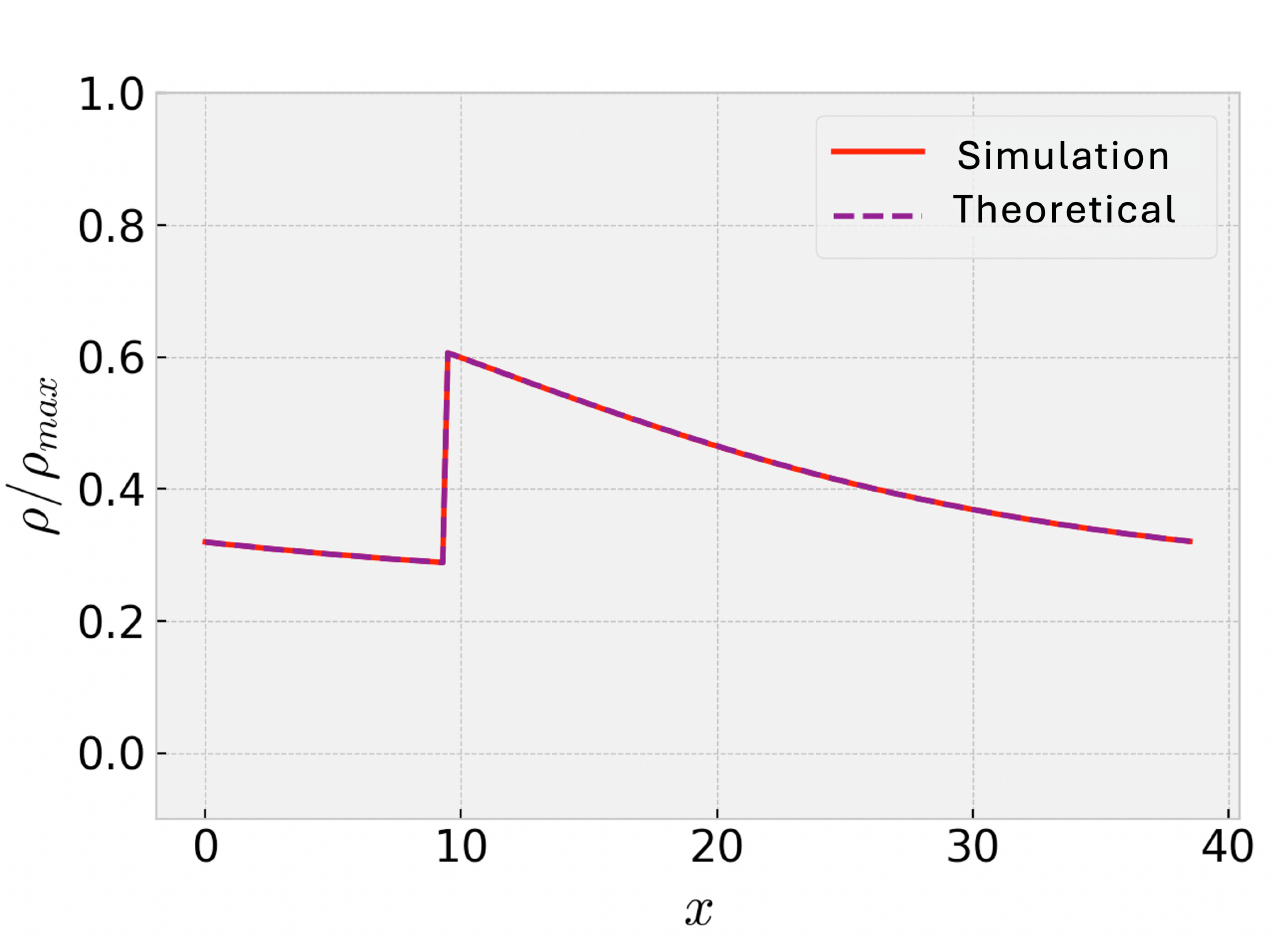} & 
   \includegraphics[width=.48\linewidth]{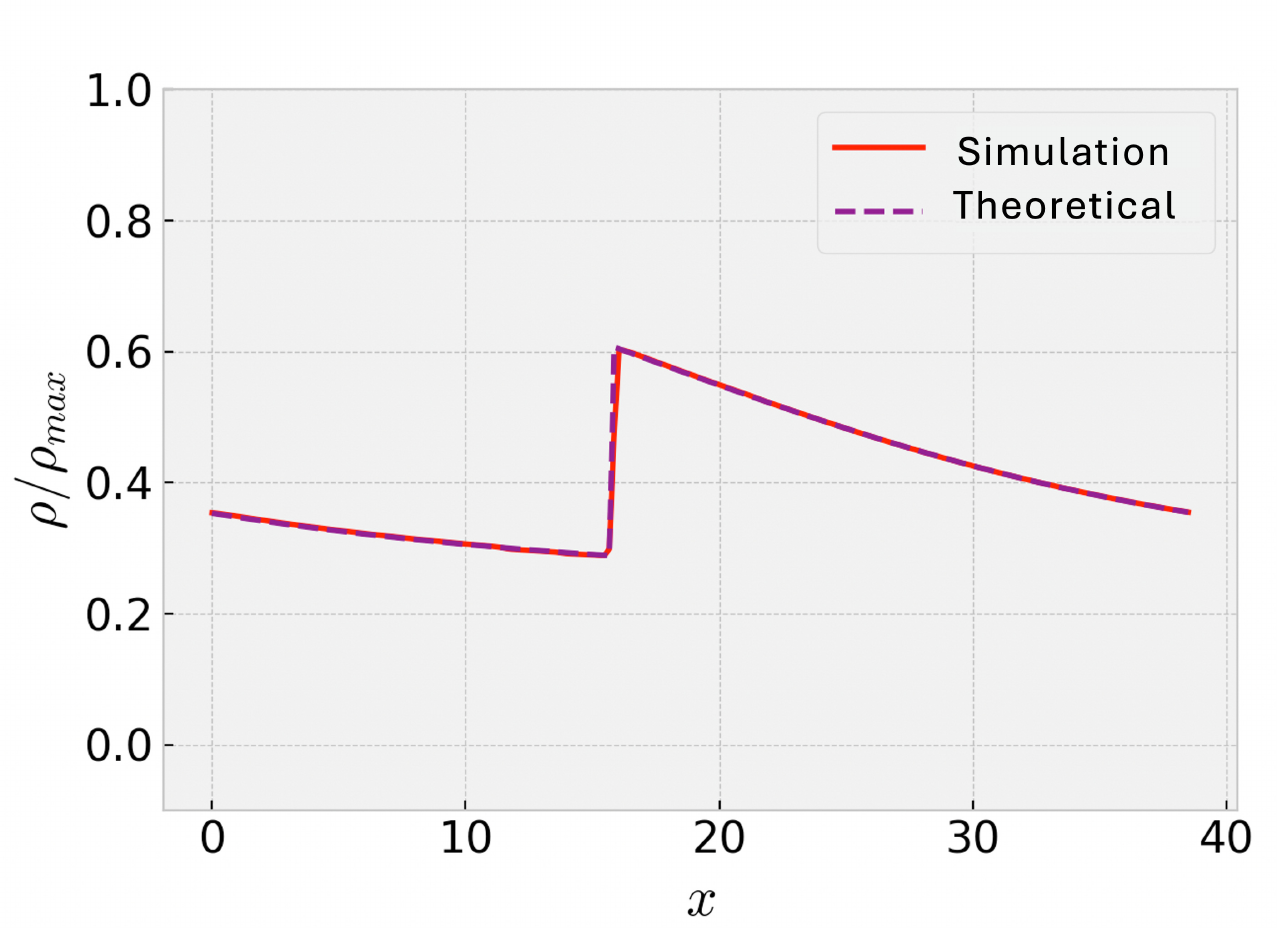}  \\  
   (c) & (d) \\ 
   \includegraphics[width=.48\linewidth]{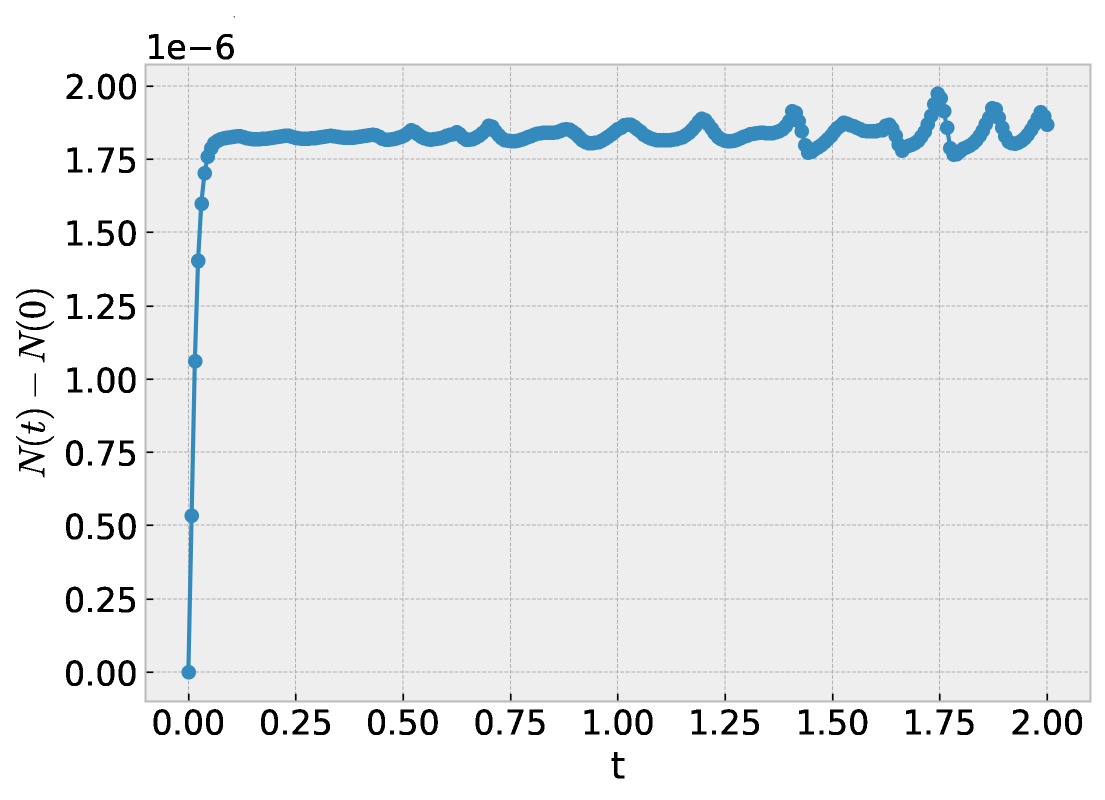} & 
   \includegraphics[width=.48\linewidth]{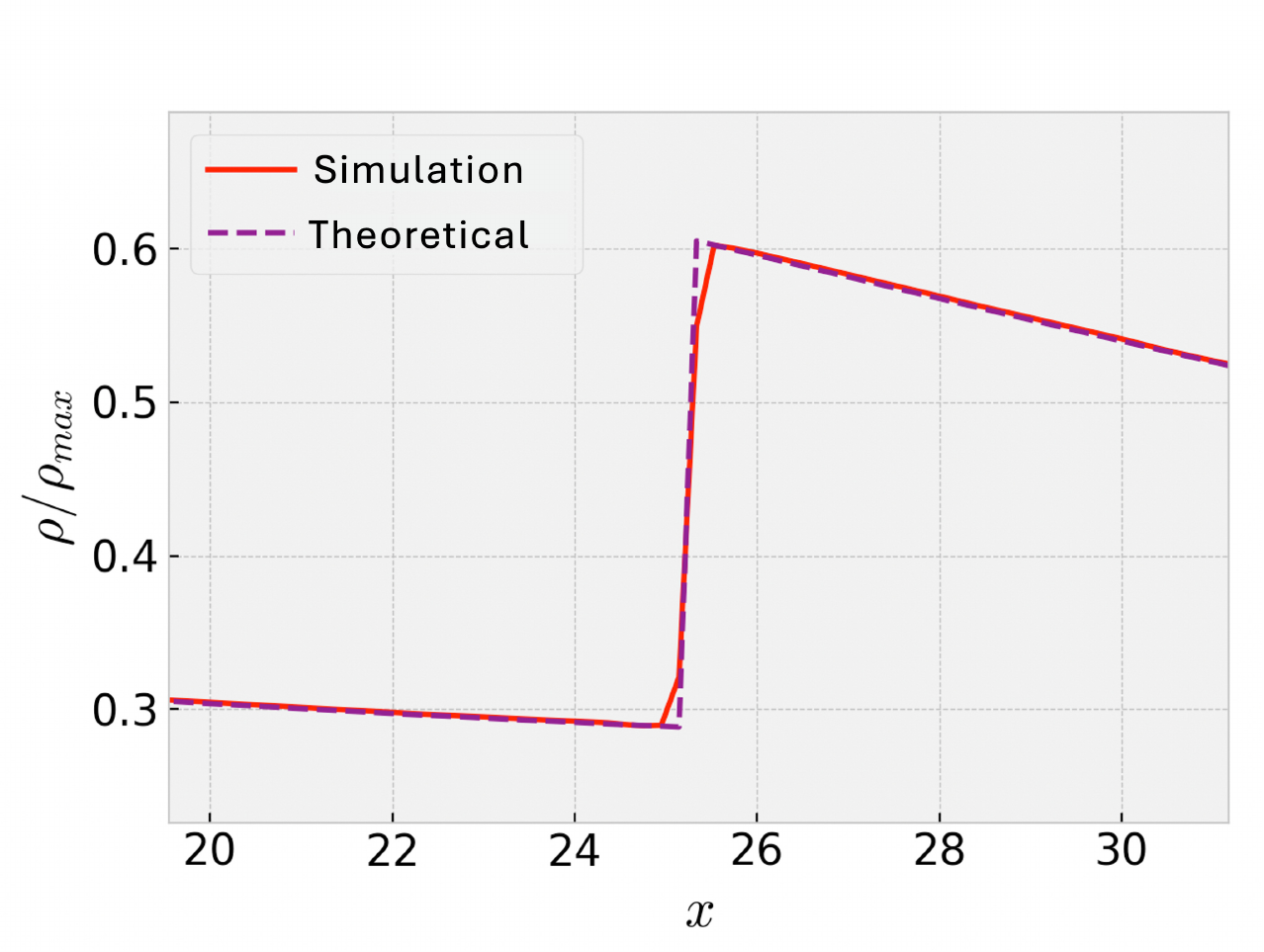} 
\end{tabular} 
   \caption{Example 1 (accuracy test): sample comparison between the numerical simulation and theoretical jamiton with $ \rho_s = 0.433\rho_{\max} $, $ v_- = 26 $, and $\mathcal{N} = 200$:  
   (a) jamiton-like initial condition, (b) comparison of exact and numerical solutions at simulated time $t=1.01$, 
    (c) numerically  simulated number of  vehicles (for the exact solution, the difference $N(t) - N(0)$ should be zero), (d) enlarged view of the  exact and numerical solutions at simulated time $t=2.47$.} 
   \label{comparativa}
\end{figure}

\begin{table}[t]
\caption{Example 1 (accuracy test): errors for various values of~$\mathcal{N}$, $\tau$,   and  simulated times $t_{\mathrm{final}}$.}
\centering  \addtolength{\tabcolsep}{-0.5pt} 
\begin{tabular}{ccccccccccccc} 
\toprule 
& \multicolumn{6}{c}{$t_{\mathrm{final}}=0.5 $} &  \multicolumn{6}{c}{$t_{\mathrm{final}}=2 $}  \\  \cmidrule(lr){2-7}\cmidrule(lr){8-13} 
& \multicolumn{2}{c}{$\tau=1 $} & \multicolumn{2}{c}{$\tau=5 $} &  \multicolumn{2}{c}{$\tau=10 $} & \multicolumn{2}{c}{$\tau=1 $} & \multicolumn{2}{c}{$\tau=5 $} &  \multicolumn{2}{c}{$\tau=10 $}  \\
  \cmidrule(lr){2-3}\cmidrule(lr){4-5}  \cmidrule(lr){6-7}\cmidrule(lr){8-9}  \cmidrule(lr){10-11}\cmidrule(lr){12-13}
$\mathcal{N}$ 
 & $\varepsilon_\rho^{\Delta x}$ & $\varepsilon_u^{\Delta x}$ & $\varepsilon_\rho^{\Delta x}$ & $\varepsilon_u^{\Delta x}$ & $\varepsilon_\rho^{\Delta x}$ & $\varepsilon_u^{\Delta x}$  & 
   $\varepsilon_\rho^{\Delta x}$ & $\varepsilon_u^{\Delta x}$ & $\varepsilon_\rho^{\Delta x}$ & $\varepsilon_u^{\Delta x}$ & $\varepsilon_\rho^{\Delta x}$ & $\varepsilon_u^{\Delta x}$  \\ \midrule 
20   & 5.908      & 2.661        & 2.519         & 1.205           & 1.779      & 0.862 & 19.760     & 9.725   & 6.055      & 2.690   & 3.524      & 1.663 \\ 
40   & 3.079      & 1.445        & 1.992         & 0.874           & 1.187      & 0.526 & 11.521     & 5.495   & 2.837      & 1.251   & 2.589      & 1.226 \\
80   & 1.437      & 0.667        & 0.594         & 0.296           & 0.957      & 0.402 &  7.052      & 3.282   & 1.308      & 0.622   & 1.092      & 0.491 \\
160  & 1.028      & 0.503        & 0.350         & 0.170           & 0.290      & 0.143 & 4.473      & 2.021   & 0.722      & 0.363   & 0.493      & 0.246  \\
320  & 0.649      & 0.304        & 0.203         & 0.111           & 0.179      & 0.088 &  2.834      & 1.235   & 0.337      & 0.187   & 0.217      & 0.127   \\
640  & 0.244      & 0.123        & 0.103         & 0.051           & 0.131      & 0.072 & 1.207      & 0.538   & 0.186      & 0.106   & 0.140      & 0.088  \\
1280 & 0.099      & 0.048        & 0.070         & 0.034           & 0.086      & 0.046  & 0.487      & 0.235   & 0.094      & 0.064   & 0.120      & 0.067  \\
2560 & 0.055      & 0.026        & 0.050         & 0.025           & 0.073      & 0.039 & 0.295      & 0.144   & 0.065      & 0.044   & 0.094      & 0.053  \\ \bottomrule 
\end{tabular}
\label{tab:errores}
\end{table}

The scheme described in Section~\ref{sec:simulacion_arz} was applied to  an initial jamiton-like solution 
  constructed in Section~\ref{cap:jam}. The simulation was run until a final time $t_{\mathrm{final}} $ and compared with the theoretical (exact)  solution given by 
    $\boldsymbol{Q}_{\mathrm{theo}} (x, t) \coloneqq  \boldsymbol{Q}_{\mathrm{jam}}((x-st)/ \tau)$, 
where $ \boldsymbol{Q}_{\mathrm{jam}} $ is obtained from the construction in Section \ref{sec:construccion}. 
 The relative  $L^1$ error 
\begin{align*} 
     \varepsilon^{\Delta x} \coloneqq  100 \frac{\| \boldsymbol{Q}_{\Delta x} - \boldsymbol{Q}_{\mathrm{theo}} \| _{L^1}}{\| \boldsymbol{Q}_{\mathrm{theo}} \|_{L^1}} 
\end{align*} 
  of a numerical solution $\boldsymbol{Q}_{\Delta x}$ computed  with a spatial discretization~$\Delta x$ 
   was calculated  for various grid refinements and values of the relaxation parameter~$\tau$. 
 See Figure~\ref{comparativa} for further details.

In Table \ref{tab:errores}, the errors obtained for $ t_{\mathrm{final}} = 0.5 $ and $ t_{\mathrm{final}} = 2 $, respectively, for three different values of $ \tau $ are shown. These errors correspond to a fixed jamiton with $ \rho_s = 0.433\rho_{\max} $ and $ v_- = 26 $. As expected, as the mesh becomes finer, the solution better approximates the theoretical solution. However, for $ \tau = 1 $, the errors are larger for both final times compared to $ \tau = 5 $ and $ \tau = 10 $. This could indicate that the simulation does not accurately approximate the exact  solution when $ \tau $ is small, especially if one wishes to study the convergence of ARZ to LWR as $ \tau $ tends to 0. This case will not be covered in this work, but in \cite[Chap. 17]{volumen_finito}, an approach and a fractional step update for $ \tau \to 0 $ are proposed.

\subsection{Example~2: emergence of jamitons}

\begin{figure}[t]
   \centering
   \begin{tabular}{cc} 
    (a) & (b) \\ 
   \includegraphics[width=.48\linewidth]{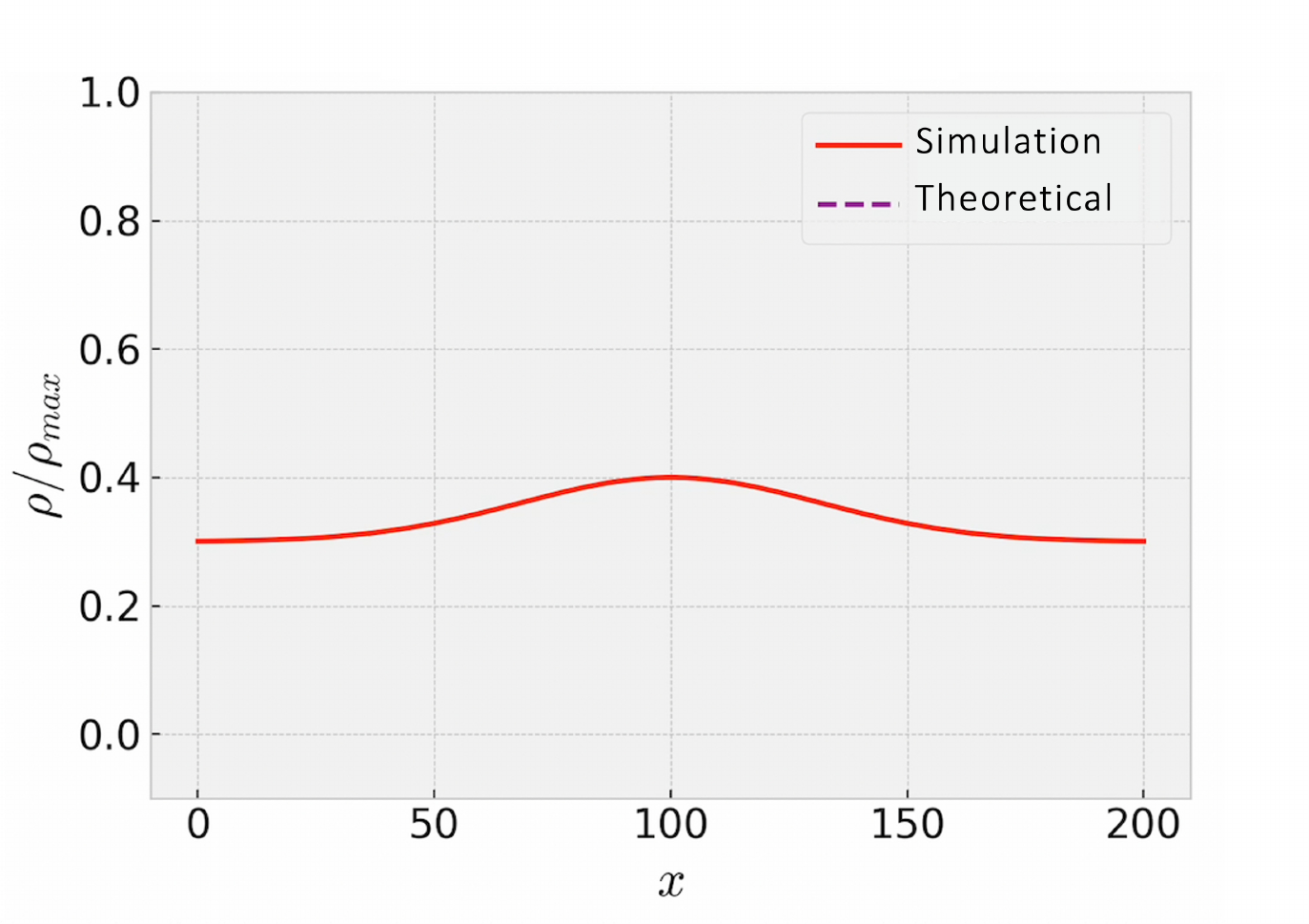} & 
   \includegraphics[width=.48\linewidth]{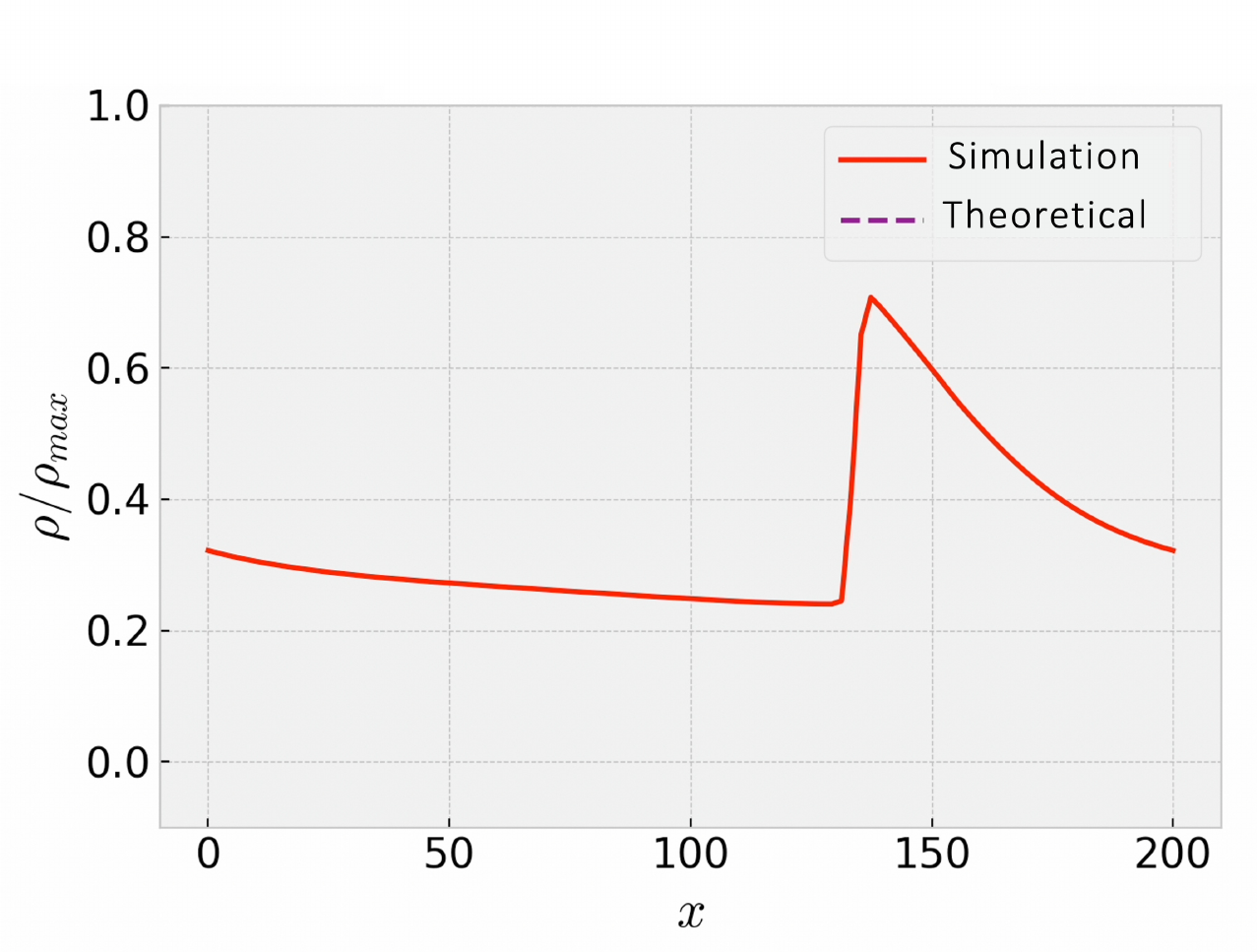}  \\
   (c) & (d) \\ 
 \includegraphics[width=.48\linewidth]{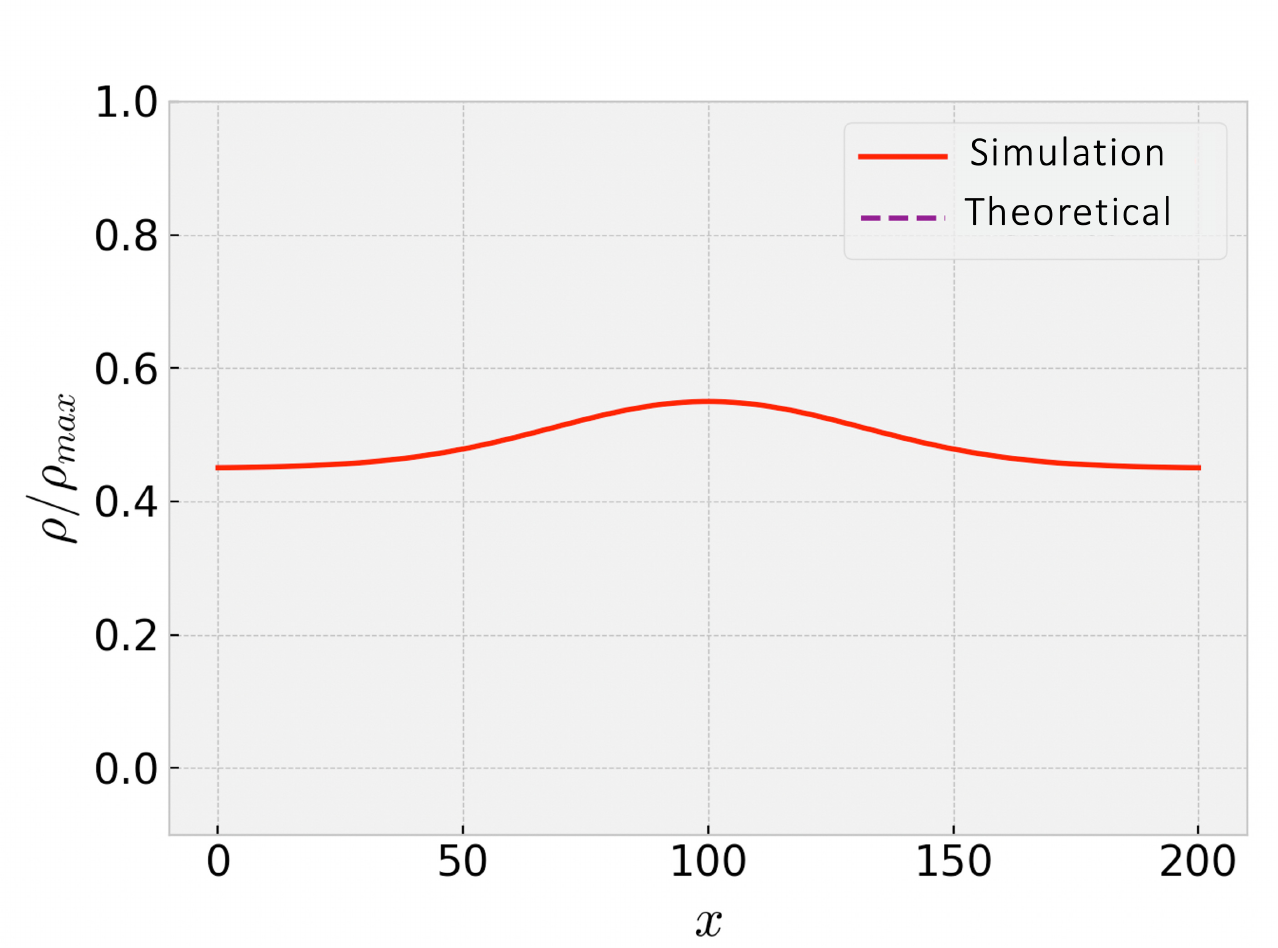} & 
   \includegraphics[width=.48\linewidth]{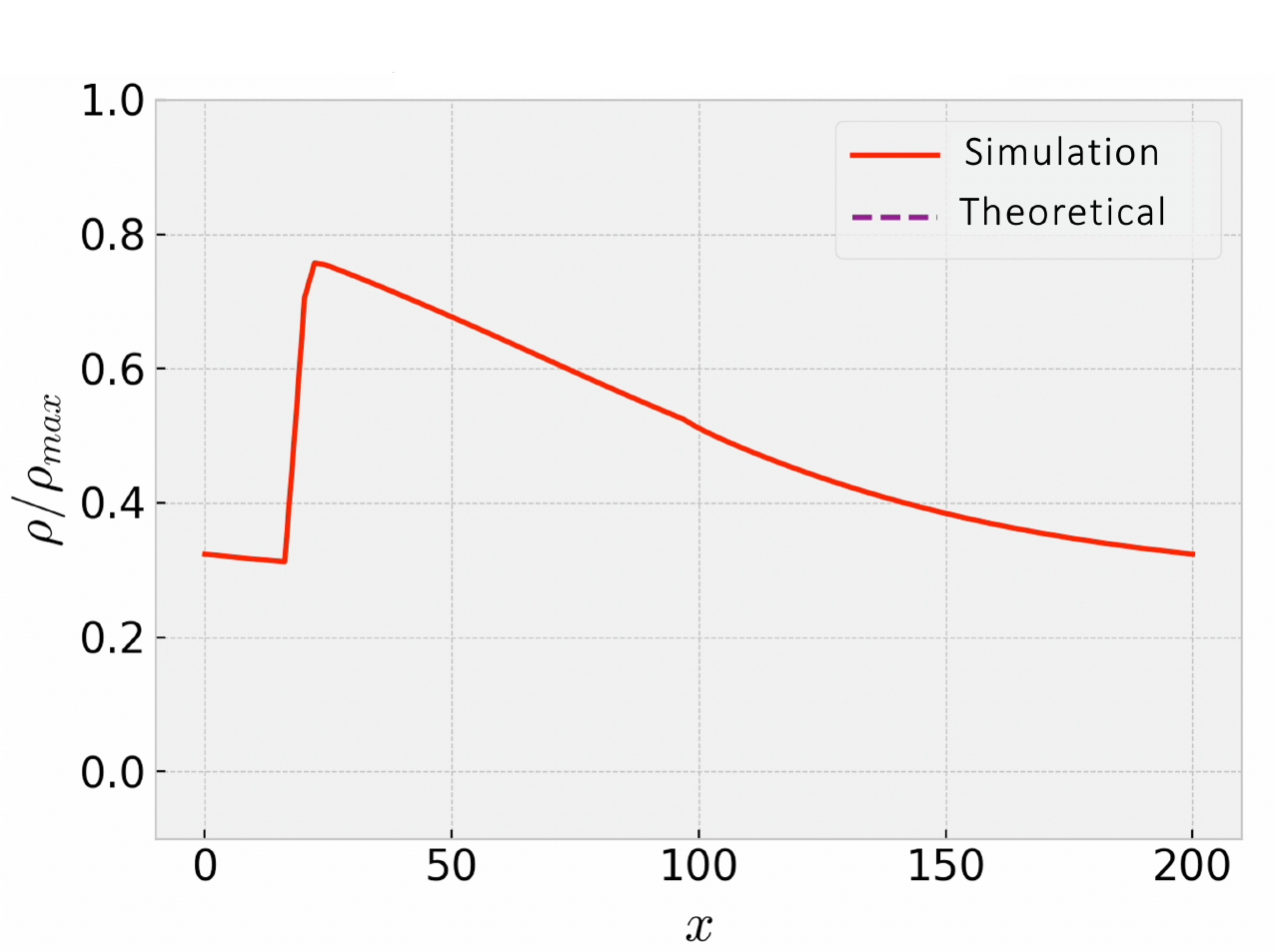} 
   \end{tabular} 
   \caption{Example~2 (emergence of jamitons): examples with (a, b) low initial density, (c, d) high initial density;  (a, c) initial condition, 
    (b, d) formation of a   jamiton at simulated time (b)  $t=13.42$, (d) $t=23.07$.} 
   \label{img:baja_densidad}
\end{figure}

In Figure~\ref{img:baja_densidad} %, \ref{img:alta_den}, and \ref{img:jamitinos}, 
 formations of jamitons are observed with initial conditions that emulate a small perturbation. The size and displacement of the obtained jamitons depend to some extent on the initial density, which in all cases was chosen to violate  the SCC.
In the case of Figures~\ref{img:baja_densidad}(a) and~(b), the jamiton travels with positive velocity to the right, 
 while in  Figures~\ref{img:baja_densidad}(c) and~(d), the resulting jamiton, besides being longer, moves with negative velocity, i.e., towards the left.

 \subsection{Example~3: jamitons on a long circular road}

\begin{figure}[t]
   \centering
   \begin{tabular}{cc} (a) & (b) \\ 
   \includegraphics[width=.48\linewidth]{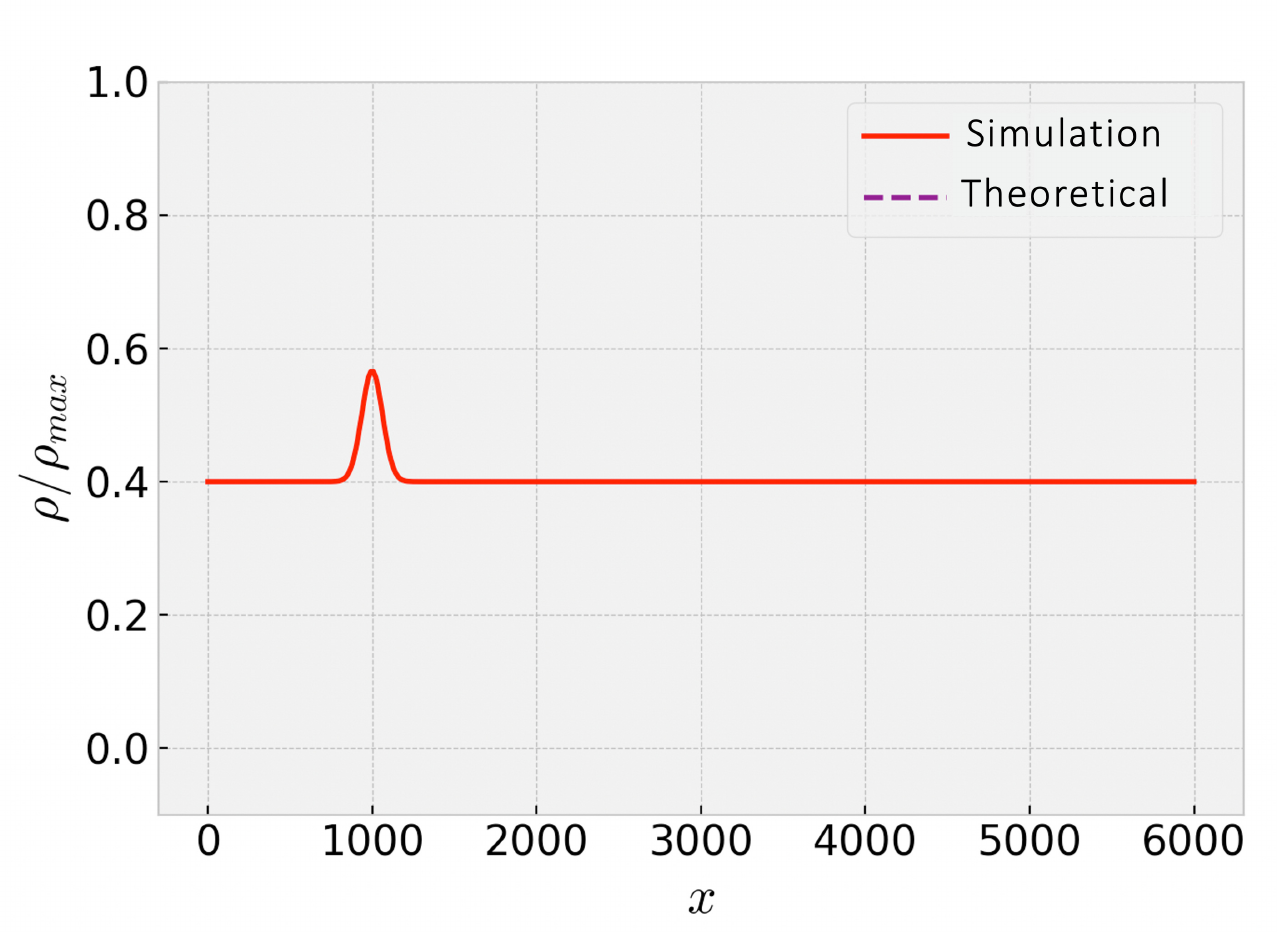} & 
   \includegraphics[width=.48\linewidth]{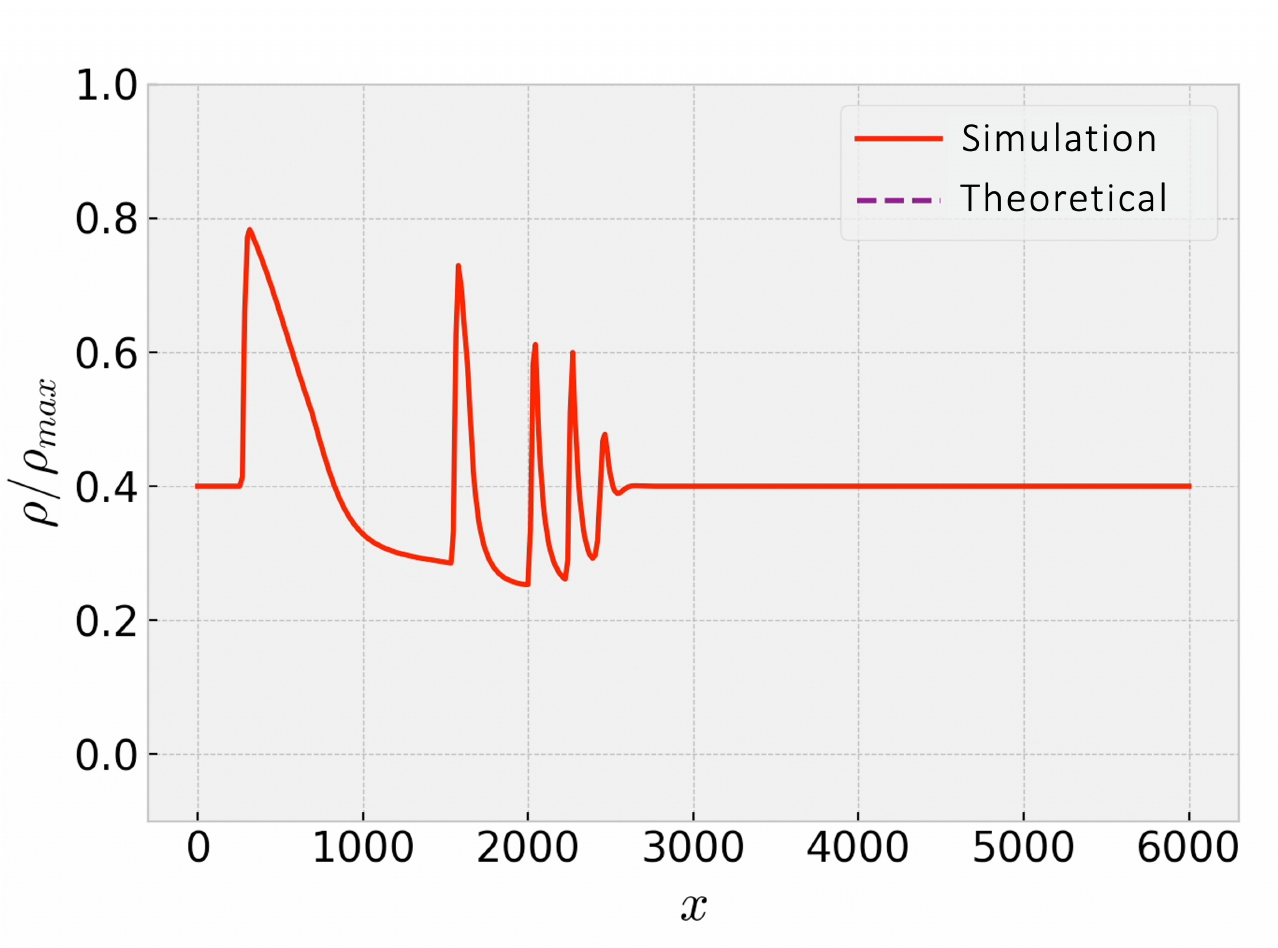}
    \end{tabular} 
    (c) \\
   \includegraphics[width=.48\linewidth]{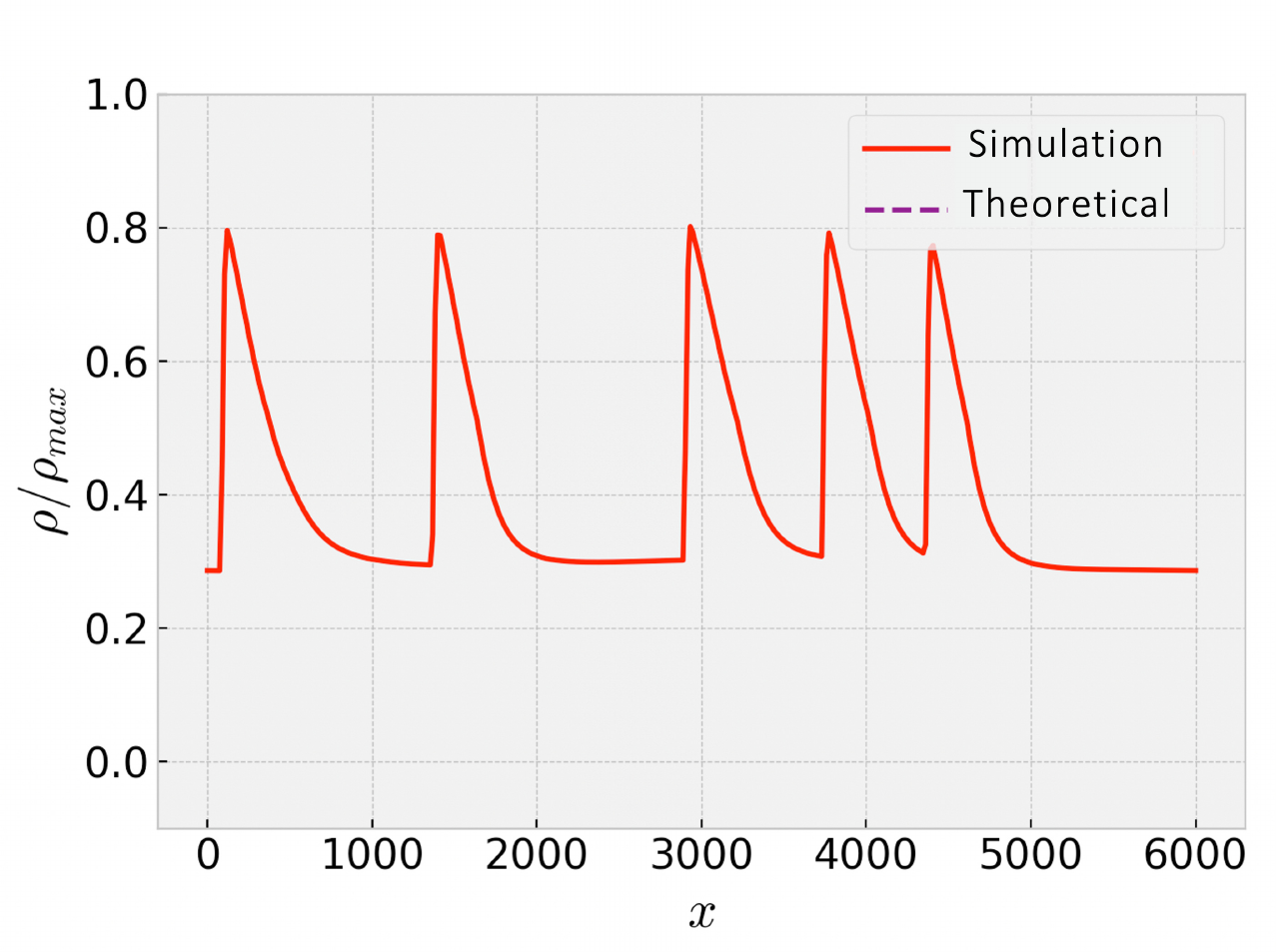}
      \caption{Example 3 (jamitons on a long circular road): (a) initial condition, 
       (b) formation of a large jamiton and several jamitinos at $ t = 117.09 $, (c)     solution is governed by the jamitonic regime at $ t = 1646.79$. }
       \label{img:jamitinos} 
\end{figure}

In Figure~\ref{img:jamitinos}, the formation of  so-called ``jamitinos'', i.e., initially small jamitons, as shown in \cite{MIT_jam}, is observed. The setting is as follows: one starts from a more pronounced Gaussian initial condition on a circular road six kilometers long, essentially approximating an infinite road. After some time, a chain of jamitinos with positive velocity forms from a larger jamiton with negative velocity. Since one is modeling a circular road, naturally the large jamiton collides with some jamitinos. What is observed is that the result of colliding two jamitons is another jamiton, but with different length and velocity compared to those initiated before the collision. This will be examined in more detail in Section \ref{cap:colision}.
In Figure~\ref{img:jamitinos}(c),    the simulation of jamitons was run for a long time. It is observed that the solution is dominated by a jamitonic regime, with several jamitons connected by jumps and of equal amplitude but apparently of different lengths. This situation can occur if each jamiton has a different sonic density $\rho_s$. It can also be noted that all jamitons in Figure~\ref{img:jamitinos}(c), being connected by jumps, share the parameter $\rho_-$. This observation will be important for establishing the collision scheme in Section \ref{cap:colision}.

\subsection{Example 4: approximation of the jamiton propagation velocity}\label{sec:aprox_m_s}

\begin{table}[t]
\caption{Example 4 (approximation of the jamiton propagation velocity): errors in  $s$ and $m$  for various values of~$\mathcal{N}$ and  $\tau$    at simulated times $t_{\mathrm{final}}=2$.}
\centering  \addtolength{\tabcolsep}{-0.5pt} 
\begin{tabular}{ccccccc} 
\toprule 
& \multicolumn{2}{c}{$\tau=1 $} & \multicolumn{2}{c}{$\tau=5 $} &  \multicolumn{2}{c}{$\tau=10 $}  \\
  \cmidrule(lr){2-3}\cmidrule(lr){4-5}  \cmidrule(lr){6-7} 
$\mathcal{N}$ 
  & $\varepsilon^{\Delta x}_s$ & $\varepsilon^{\Delta x}_m$ & $\varepsilon^{\Delta x}_s$ & $\varepsilon^{\Delta x}_m$ & $\varepsilon^{\Delta x}_s$ & $\varepsilon^{\Delta x}_m$ \\  \midrule 
 20   & 0.01623   & 0.03049   & 0.00371   & 0.01826   & 0.03326   & 0.02255   \\
    40   & 0.00312   & 0.01316   & 0.00522   & 0.00253   & 0.01849   & 0.01286   \\
    80   & 0.00410   & 0.00235   & 0.00566   & 0.00125   & 0.01087   & 0.00785   \\
    160  & 0.00429   & 0.00007   & 0.00526   & 0.00272   & 0.00701   & 0.00533   \\
    320  & 0.00423   & 0.00128   & 0.00487   & 0.00328   & 0.00473   & 0.00374   \\
    640  & 0.00171   & 0.00048   & 0.00248   & 0.00170   & 0.00209   & 0.00162   \\
    1280 & 0.00029   & 0.00008   & 0.00233   & 0.00179   & 0.00091   & 0.00070   \\
    2560 & 0.00025   & 0.00006   & 0.00165   & 0.00131    & 0.00063   & 0.00051  
\\ \bottomrule 
\end{tabular}
\label{tab:error_m_s}
\end{table}

An important value to obtain from the simulation of jamitons is the propagation velocity, especially if one wishes to study numerically jamitons and the theoretical parameters of the jamiton are unknown. A first approximation would be to obtain~$\rho_+$ and~$\rho_-$ from the simulation, since  $\rho_+$ and $\rho_+$ are the maximum and minimum of 
 $\rho(x,t)$ for a jamiton-like solution at arbitrary, sufficiently large  time~$t$ (i.e., when the solution has attained the traveling wave regime). With these parameters, the value of~$\rho_s$ can be approximated using the construction of a jamiton. However,  this method  has a high error rate since it depends only on two points of the simulation, and the error of the numerical scheme could propagate in the estimation of~$\rho_s$. Therefore,  another way to approximate~$s$ and~$m$ that has higher precision is desirable.
 In Section~\ref{cap:jam} we demonstrated  that  jamitons    corresponded to line segments    in the fundamental diagram, 
  and that these line  segments   intersect the equilibrium curve at $(\rho_s, Q(\rho_s))$, where  $m$ and $s$ correspond to the intersection and slope of the line passing through the segment, respectively. 
   This result allows us  to approximate the  velocity of propagation of a jamiton by resorting to the $(\rho, q)$ plane. The method for obtaining this was presented in \cite{stabilidad_jam}. It consists 
    in  taking a jamiton-like solution obtained from the simulation $\rho_{\mathrm{jam}}$ and $u_{\mathrm{jam}}$, and plotting the obtained data in the $(\rho, q)$ plane (i.e., plotting $(\rho_{\mathrm{jam}}, \rho_{\mathrm{jam}} u_{\mathrm{jam}})$). The expected result would correspond to a line segment with the parameters of the jamiton. Thus, $m_{\mathrm{jam}}$ and $s_{\mathrm{jam}}$ can be obtained from a linear regression of the mentioned plot.

To observe the effectiveness of the method, a jamiton with $\rho_s = 0.433\rho_{\max}$ and $v_- = 26$ was taken as the initial condition (as in Example~1).
  The relative error between $s$ and $m$ of the jamiton and the one obtained by linear regression was calculated, given by:
\[%begin{equation}
    \varepsilon^{\Delta x}_s = 100\frac{|s_{\Delta x} - s_{\mathrm{theo}} |}{|s_{\mathrm{theo}} |}, \qquad  \varepsilon^{\Delta x}_m = 100\frac{|m_{\Delta x} - m_{\mathrm{theo}} |}{|m_{\mathrm{theo}} |}.
\]%end{equation}
For this particular jamiton, we have $m_{\mathrm{theo}}  \approx 0.356$ and $s_{\mathrm{theo}}  \approx 6.374$. The effectiveness of the method was tested by running the simulation until $t_{\mathrm{final}} = 2$ for $\tau = 1, 5, 10$ with increasingly finer grids. The results are shown in Table \ref{tab:error_m_s}. The effectiveness of the method is evident from Table \ref{tab:error_m_s}, even with relatively coarse grids, and its precision increases as the grid becomes finer. This method will be used in Section \ref{cap:colision} to determine the exit velocities of the jamitons resulting from the collision of jamitons.

\section{Collision of jamitons}\label{cap:colision}

\subsection{Compatibility of jamitons}
\begin{figure}[t]
   \centering
   \includegraphics[width=.5\linewidth]{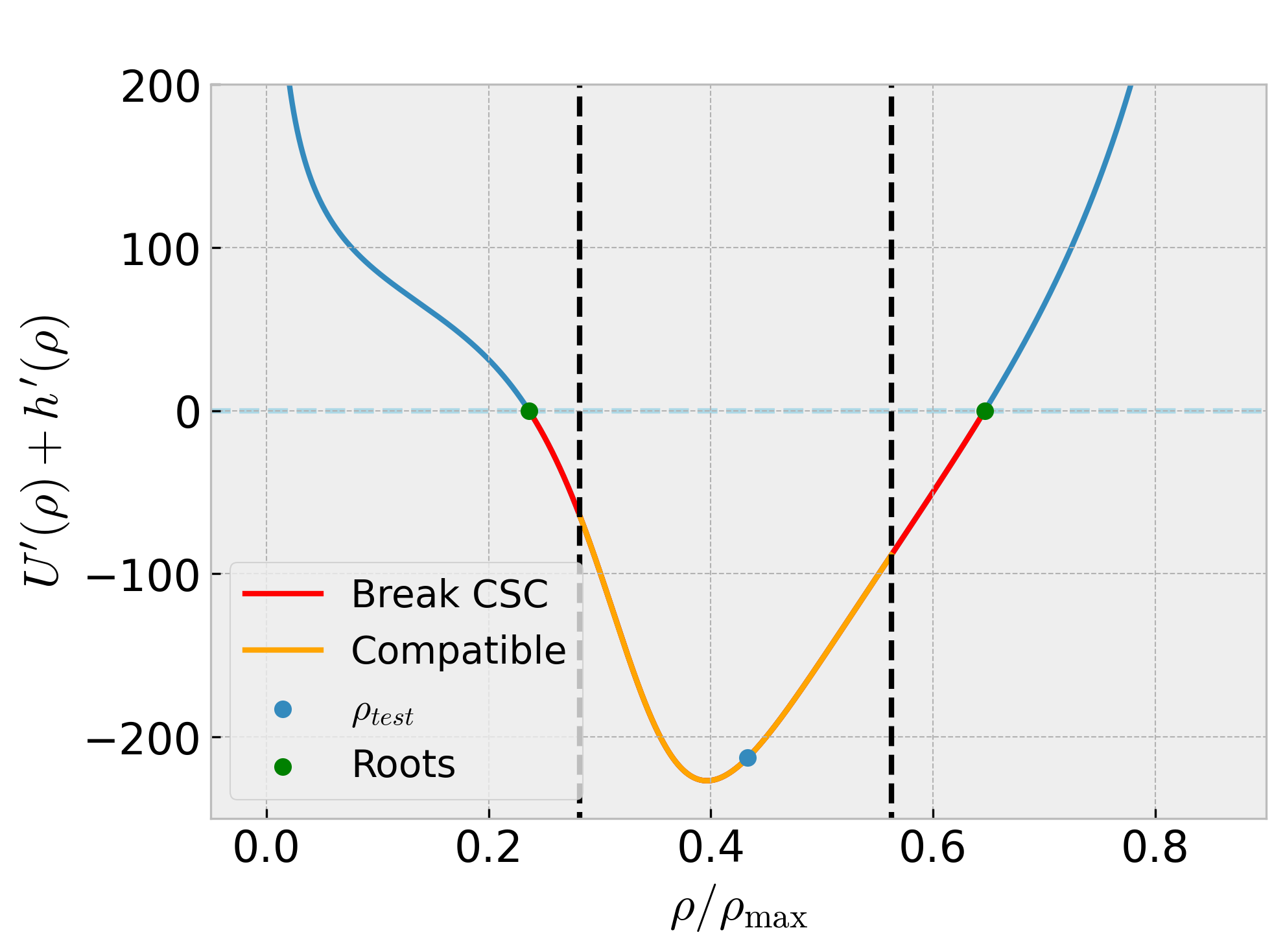}
   \caption{Compatible densities with a test jamiton. In red, possible densities that break the SCC condition.}%Densidades compatibles con jamiton de prueba.}
   \label{img:compatibles}
\end{figure}

As mentioned in Section \ref{cap:ARZ}, the collision of jamitons was observed in the formation of jamitinos on a circular route. In such a case, the collision of two jamitons results in a new jamiton with parameters possibly different from those of the original jamitons. Additionally, it was observed that for two jamitons to collide, \emph{they must share the value of $\rho_-$}, as they necessarily must be connected by the jump between $\rho_-$ and $\rho_+$ (or in Lagrangian variables, between $v_+$ and $v_-$).
Indeed, consider two jamiton-type solutions $\smash{Q_{\mathrm{jam}}^1}$ and $\smash{Q_{\mathrm{jam}}^2}$ with parameters~$\rho_s^1$, $\v_-^1$, $\v_+^1$, and  and $\rho_s^2$, $\v_-^2$, and $\v_+^2$, respectively. The idea is to obtain a necessary condition for $\smash{Q_{\mathrm{jam}}^1}$ to be able to join $\smash{Q_{\mathrm{jam}}^2}$. The smooth part of $\smash{Q_{\mathrm{jam}}^1}$ is obtained by integrating the ODE \eqref{eq:EDO_jam} 
between~$\v_+^1$ and~$\v_-^1$. However, ``the jump of $\smash{Q_{\mathrm{jam}}^2}$'' connects the values~$\v_-^2$ with~$\v_+^2$, so if one wants to join the smooth part of $\smash{Q_{\mathrm{jam}}^1}$ with 
 ``the jump of $\smash{Q_{\mathrm{jam}}^2}$,''  necessarily $\v_-^1 = \v_-^2$. This condition will be known as the {\em compatibility condition} between jamitons.
\begin{defn}[Compatibility]
    Two jamitons $\smash{Q_{\mathrm{jam}}^1}$ and $\smash{Q_{\mathrm{jam}}^2}$ with parameters $\v_-^1$ y $\v_-^2$, respectively, will be compatible if $\v_-^1 = \v_-^2$.
\end{defn}
It can be noted that the compatibility condition does not completely restrict the values of $\rho_s$ or $\v_+$, even though the latter depends on both $\rho_s$ (since $m=m(\rho_s)$ and $s = s(\rho_s)$) and $\v_-$. Since the family of jamiton solutions is parametrized by $\rho_s$, it is possible to study the compatibility of various jamitons of different sizes and lengths based on their sonic density for a fixed jamiton given by $\rho_s$. First, a $\rho_{\mathrm{test}}  $ will be chosen such that the selection interval for $\v_-$ is the largest among the densities that break the SCC, and $\v_-^{ \mathrm{test}}$ will be chosen as the midpoint of the interval found. Then, sonic densities will be chosen such that the maximal jamiton interval given by $[\v_R, \v_M]$ contains $\v_-^{ \mathrm{test}}$. The densities found are shown in Figure \ref{img:compatibles}. The test jamiton chosen has parameters $\rho_s \approx 0.4333$ and $\v_-^{ \mathrm{test}} \approx 26.602$.

This ensures that the study of collisions considers a variable amount of different sizes and lengths, and that the jamitons to collide are compatible, as $\v_-^{ \mathrm{test}}$ can be chosen within the jamitons to collide. However, the first step of choosing $\rho_s^{ \mathrm{test}}$ can be done for any $\rho_s$ that breaks the SCC. The particular choice made was to maximize the number of possible collisions with a single test jamiton. 

\subsection{Example 5: collision between two compatible  jamitons}

\begin{figure}[t]
   \centering
   (a) \\ 
   \includegraphics[width=0.9\linewidth]{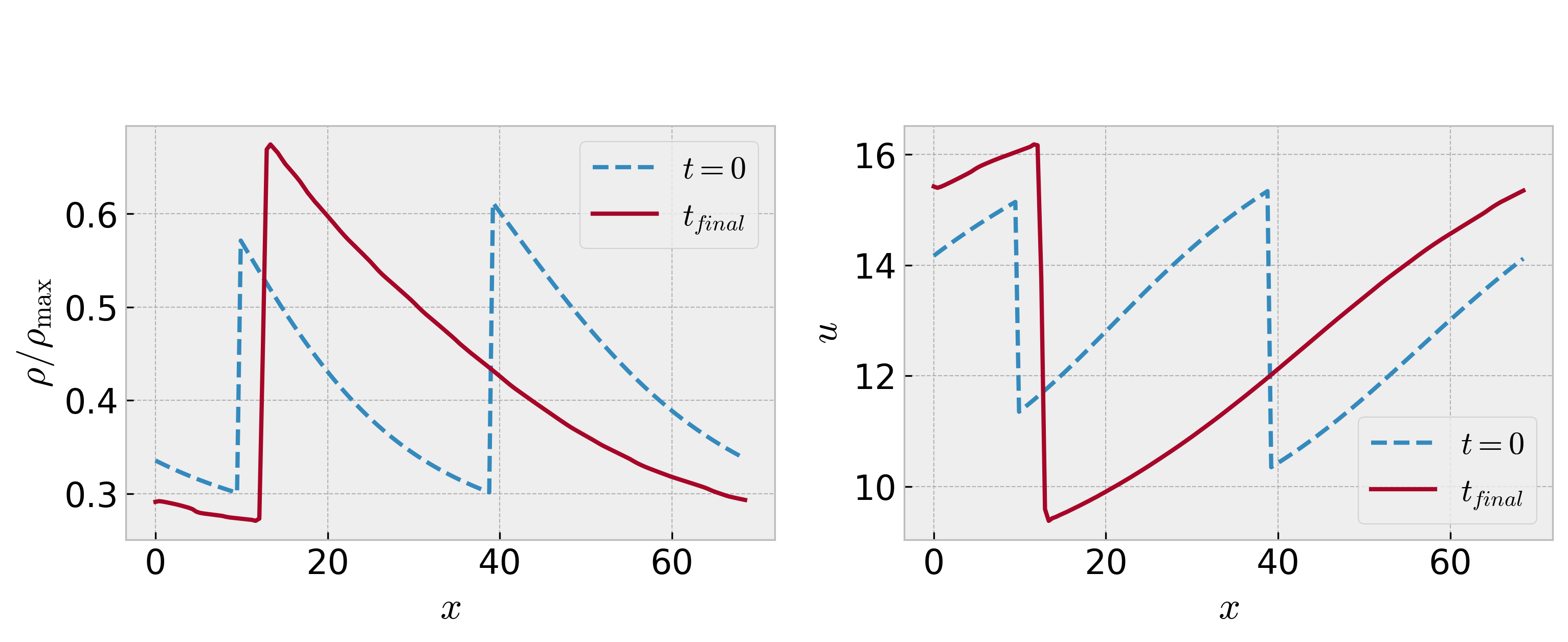} \\
   \begin{tabular}{cc} 
   (b) & (c) \\ 
   \includegraphics[width=.43\linewidth]{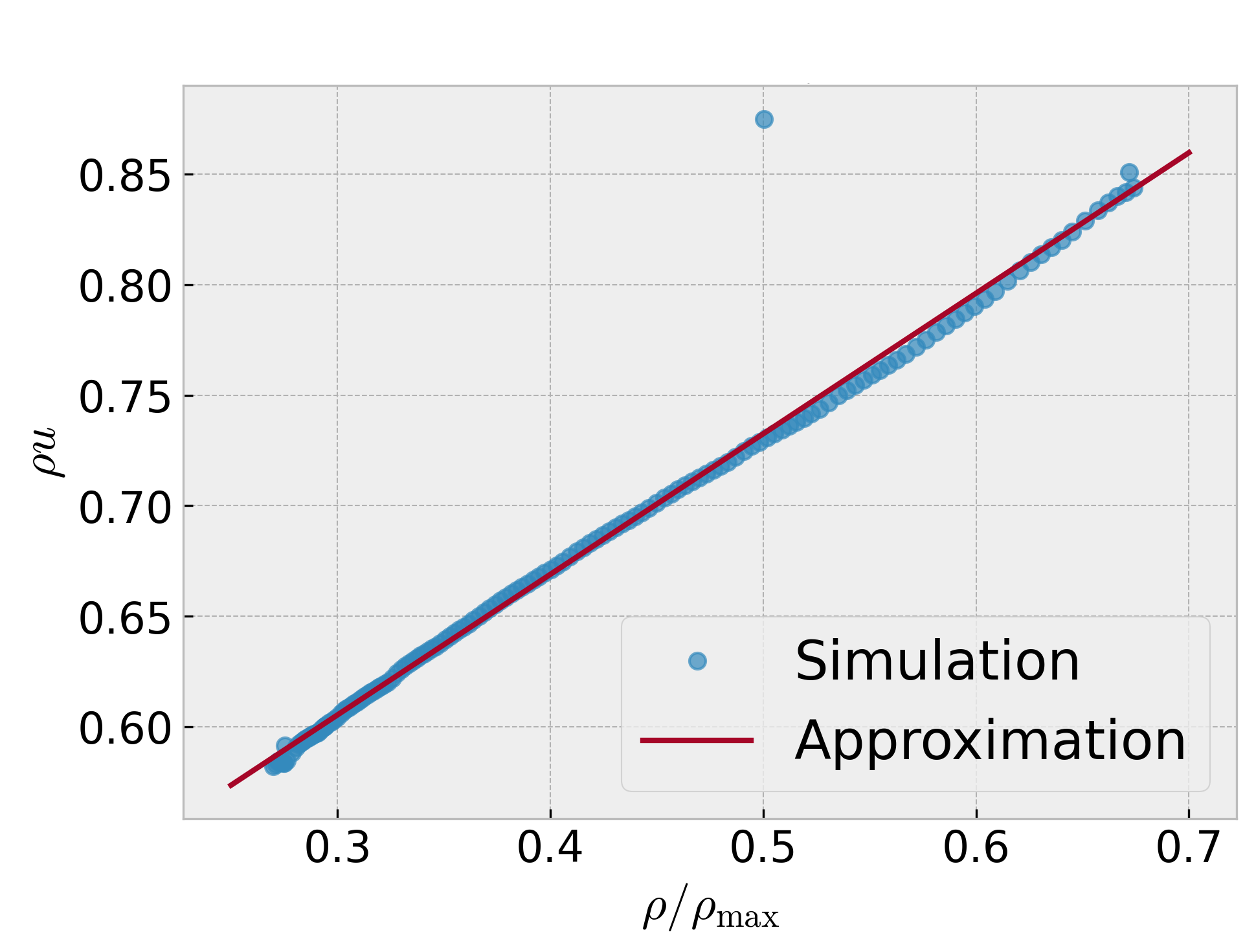} & 
   \includegraphics[width=.49\linewidth]{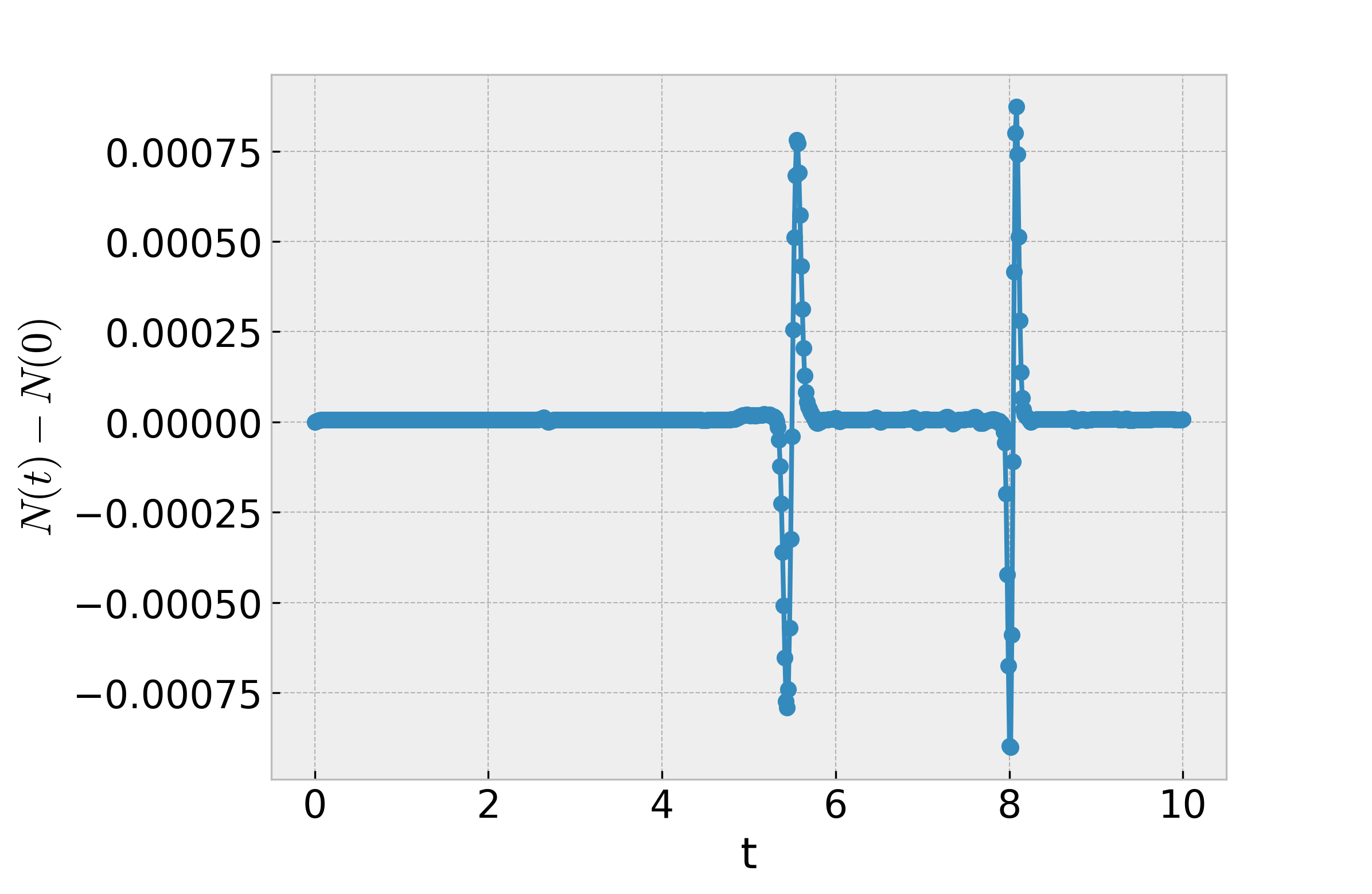} 
   \end{tabular} 
   \caption{Example 5 (collision between two compatible  jamitons): (a) collision between a jamiton with $\rho_s = 0.425\rho_{\max}$ and another with $\rho_s = 0.443\rho_{\max}$, both with $\v_- = 25$
   shown in the variables (left)~$\rho$ and (right)~$u$, (b) 
   segment in the $(\rho,q)$-plane approximating the jamiton obtained from the collision in~(a)
    (supporting the appearace of a post-collision jamiton), (c) illustration that
    the collision algorithm maintains the conservation law with little variation in the total number of vehicles.}
   \label{img:colision}
\end{figure}

As mentioned in Section \ref{cap:ARZ}, in the simulation of jamitinos, various collisions between jamitinos with different velocities are observed, as well as some being absorbed by the larger initial jamiton. This leads to the following conjecture, whose validation at the moment is purely numerical and can be observed in Figure~\ref{img:colision}.
\begin{conj}
The collision of two compatible jamitons generates a unique new jamiton whose parameters depend on the original jamitons.
\end{conj}
This conjecture is crucial as the simulations are conducted with the anticipation of the emergence of a new jamiton. The selected jamitons are merged with the test jamiton, which is then imposed as the initial condition in the simulation. The code is then allowed to run until a sufficiently long time has elapsed for the resulting jamiton from the collision to form. Since this jamiton is obtained purely numerically, to obtain its properties, such as $m$ and $s$, the method described in Section \ref{sec:aprox_m_s} is used based on the behavior of a jamiton in the fundamental diagram. Figure \ref{img:colision}(a)  shows an example of how two jamitons collide and the resulting jamiton after a time $t_{\mathrm{final}}$. In  Figure~\ref{img:colision}(b), the segment approximating the post-collision jamiton in the $(\rho, q)$ plane is observed.

\subsection{Example 6: simulation of multiple collisions}

\begin{figure}[t] 
\centering
\begin{tabular}{cc} 
(a) & (b) \\ 
 \includegraphics[width=.46\linewidth]{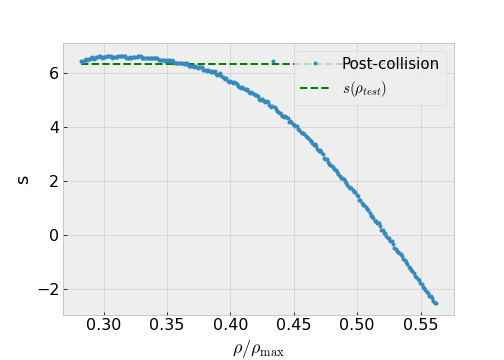}  & \includegraphics[width=.46\linewidth]{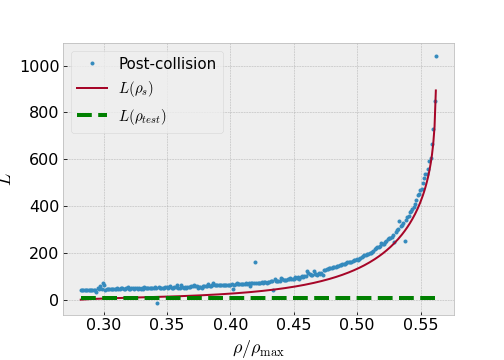} \\ 
 (c) & (d) \\ 
  \includegraphics[width=.46\linewidth]{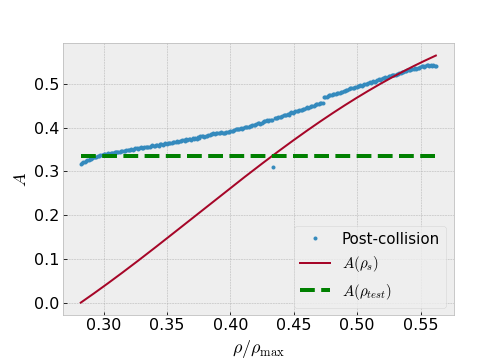} &   \includegraphics[width=.46\linewidth]{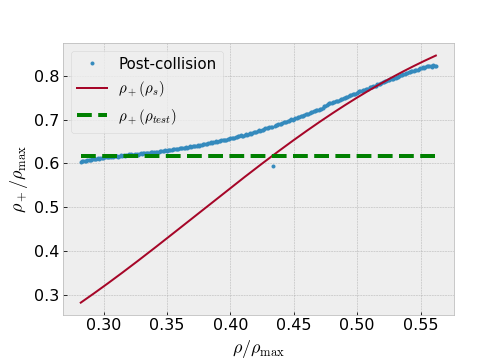} 
 \end{tabular} \\
 (e) \\
  \includegraphics[width=.46\linewidth]{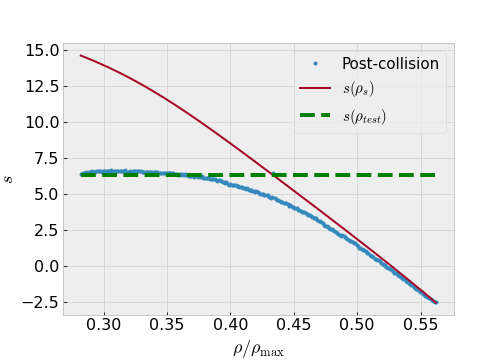}
\caption{Example 6 (simulation of multiple collisions): (a) plot of post-collision exit velocities (the data at $\rho/\rho_{\max} \approx 0.44$ corresponds to jamitons that did not collide because they had very close values of $\rho_s$;  note the small initial bump that appears in the post-collision exit velocities), (b) lengths obtained after each collision, (c) amplitudes obtained after each collision, (d) 
 values of $\rho_+$ after each collision, (e)  comparison with function $s(\rho_s)$. \label{img:colisiones_vel}}
\end{figure}  

With the proposed algorithm, 274 different collisions are performed between the test jamiton and the compatible jamitons, 
 where we use  $\tau=5$ and  $\mathcal{N}=160$ as the total number of points in the grid, and $\Delta t$ such that the CFL condition $s_{\max} \Delta t/  \Delta x  \leq 1/2$ 
is satisfied. The simulations of the various collisions were run in parallel on 20 HPC clusters. The resulting jamitons from the collisions have their propagation velocities calculated, as seen in Figure~\ref{img:colisiones_vel}(a)  as a function of the densities $\rho_s$ compatible with the test jamiton.
 Other properties can be calculated from the post-collision jamitons, such as length and amplitude given by equations \eqref{eq:Largo} and \eqref{eq:Amplitud}. 
 
 Figures~\ref{img:colisiones_vel}(b) and~(c) show  the lengths and amplitudes of the resulting jamitons as a function of the compatible $\rho_s$ values and compared to those of the jamitons before colliding. The dashed line represents the test jamiton. Figure~\ref{img:colisiones_vel}(d)   presents the values of $\rho_+$ after each collision. Similarly, the exit velocities can be compared along with the function $s(\rho_s)$, resulting in the graph in Figure~\ref{img:colisiones_vel}(e).

An interesting result can be observed in Figure~\ref{img:colisiones_vel}(a) for $\rho_s/\rho_{\max} \in [0.26, 0.35]$. Note that for some collisions, the exit velocities increased compared to the test jamiton. This could indicate that a jamiton can accelerate if it collides with smaller-sized jamitons, and if the relationship between velocity and size meets certain properties, it could even decrease the size of a jamiton. This provides an important tool in practical applications, as traffic could potentially be disentangled by colliding with chains of smaller-sized jamitons. On the other hand, in Figure~\ref{img:colisiones_vel}(e) it can be observed that the exit velocities correspond to a smoothing between the constant $\rho_ \mathrm{test}$ and the function $s(\rho_s)$. This indicates that the predominant jamitons after the collision are those of larger size, meaning that the smaller jamitons in length and amplitude are absorbed by the larger ones. However, in the collisions of jamitons with $\rho_s/\rho_{\max} \in [0.40, 0.45]$, there is a certain additivity in their amplitudes, as the amplitude obtained after the collision turns out to be greater than the amplitudes before the collision. A contrary effect, observable in Figure~\ref{img:colisiones_vel}(e), occurs with the exit velocities, which decrease compared to the input jamitons for the same density range. This phenomenon can be summarized in the following conjecture:
\begin{conj}[Post-collision inequalities]%Cotas post-colisión]
Let $Q_{\mathrm{jam}}$ be a jamiton with sonic density $\rho_s$. Let $\smash{\{Q_{\mathrm{jam}}(\rho_s^{ \mathrm{jam}})\}}$ be a corresponding non-empty set of compatible jamitons 	parametrized by their sonic densities $\rho_s^{ \mathrm{jam}}$. Let $Q_{\mathrm{post}}$ be the jamiton obtained from the collision of  $Q_{\mathrm{jam}}$ with $\smash{Q_{\mathrm{jam}}(\rho_s^{ \mathrm{jam}})}$ and with sonic density $\rho_s^{\mathrm{post}}$. Then the following hold:
    \begin{itemize}
       \item There exists a $\rho_s^0$ such that $\smash{s(\rho_s^{\mathrm{post}}) \leq s(\rho_s)}$ for all $\rho_s^\mathrm{jam} \geq \rho_s^0$ and 
        $s(\rho_s^{\mathrm{post}}) \leq s(\rho_s^ \mathrm{jam})$ for all~$\rho_s^ \mathrm{jam}$.
        \item There exists a $\rho_s^1$ such that $A(\rho_s^{\mathrm{post}}) \geq A(\rho_s^ \mathrm{jam})$ for all $\rho_s^\mathrm{jam} \leq \rho_s^1$.

        \item $L(\rho_s^\mathrm{post}) \geq L(\rho_s)$ and $L(\rho_s^\mathrm{post}) \geq L(\rho_s^{\mathrm{jam}})$ for all $\rho_s^{ \mathrm{jam}}$.    
     \end{itemize}
\end{conj}

\subsection{Example~7: effect of the relaxation time}

\begin{figure}[t]
   \centering
   \begin{tabular}{cc} 
   (a) & (b) \\ 
   \includegraphics[width=.48\linewidth]{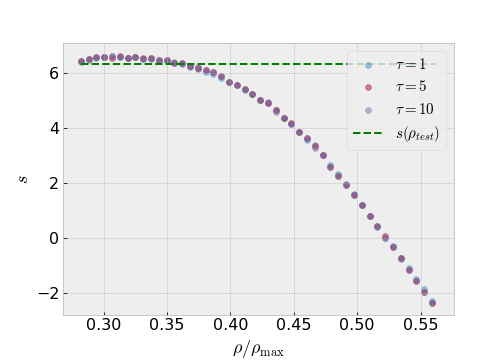} & \includegraphics[width=.48\linewidth]{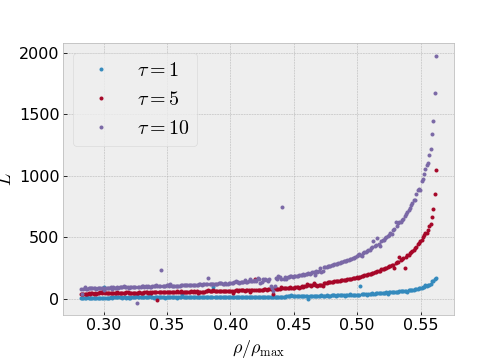} \\ 
   (c) & (d) \\
    \includegraphics[width=.48\linewidth]{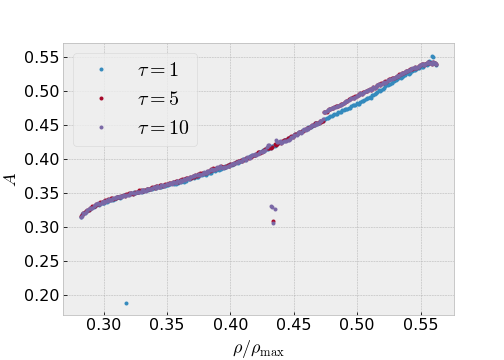} & 
    \includegraphics[width=.48\linewidth]{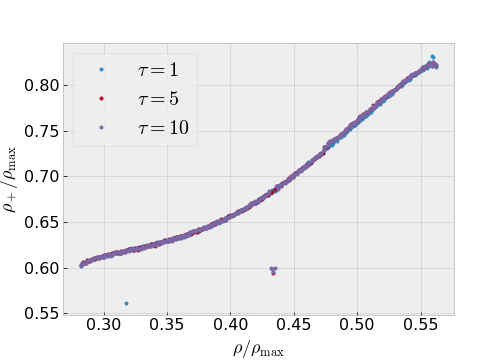}
   \end{tabular} 
   \caption{Example~7 (effect of the relaxation time): (a) plot of post-collision exit velocities , (b) lengths and (c) amplitudes 
    of the resulting jamiton, (d) values of~$\rho_+$, all for various 
     values of~$\tau$.} 
   \label{img:colisiones_taus}
\end{figure}

\begin{figure}[t]
   \centering
   \begin{tabular}{cc} 
   (a) & (b) \\ 
   \includegraphics[width=.48\linewidth]{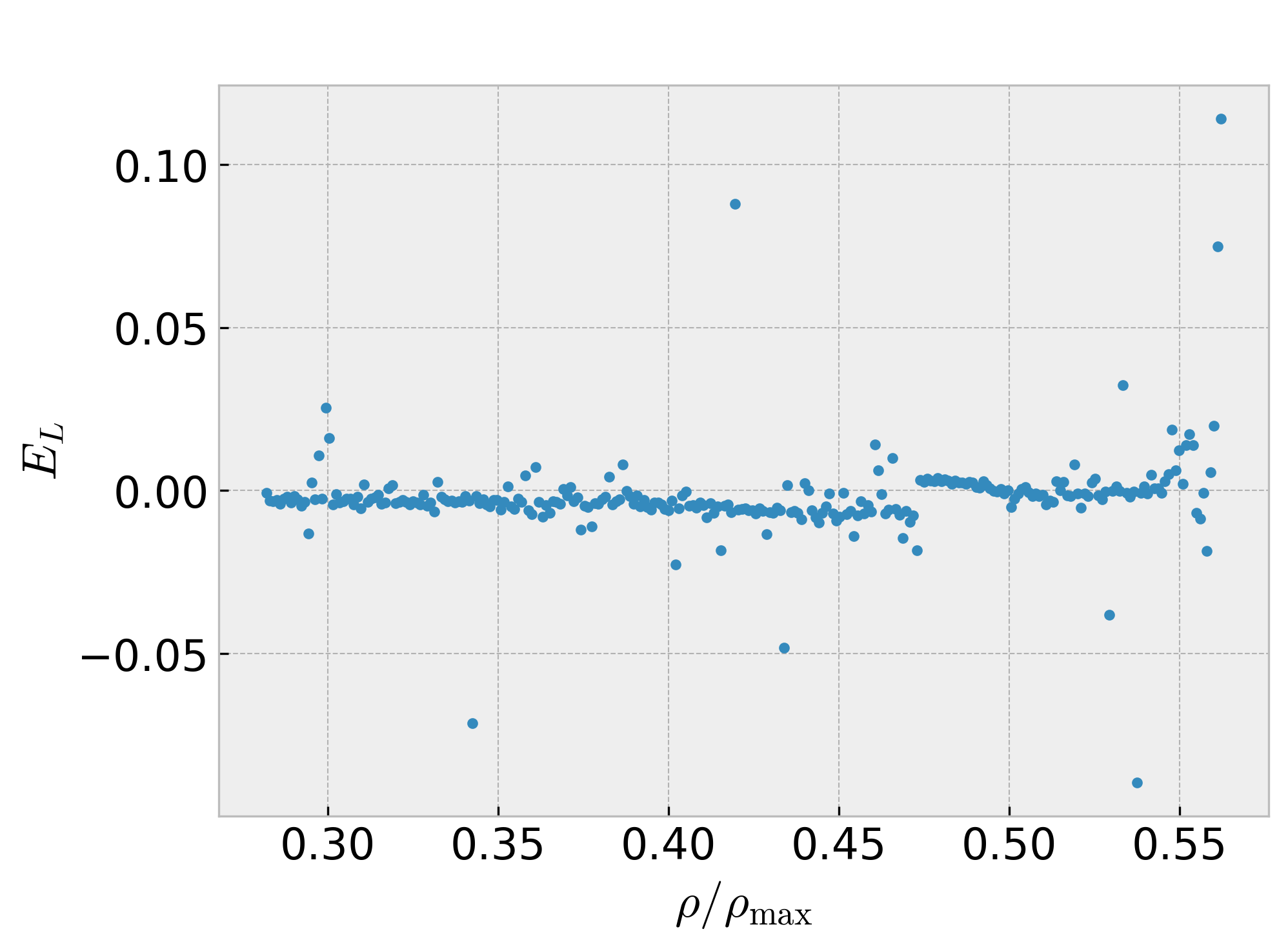} & 
    \includegraphics[width=.48\linewidth]{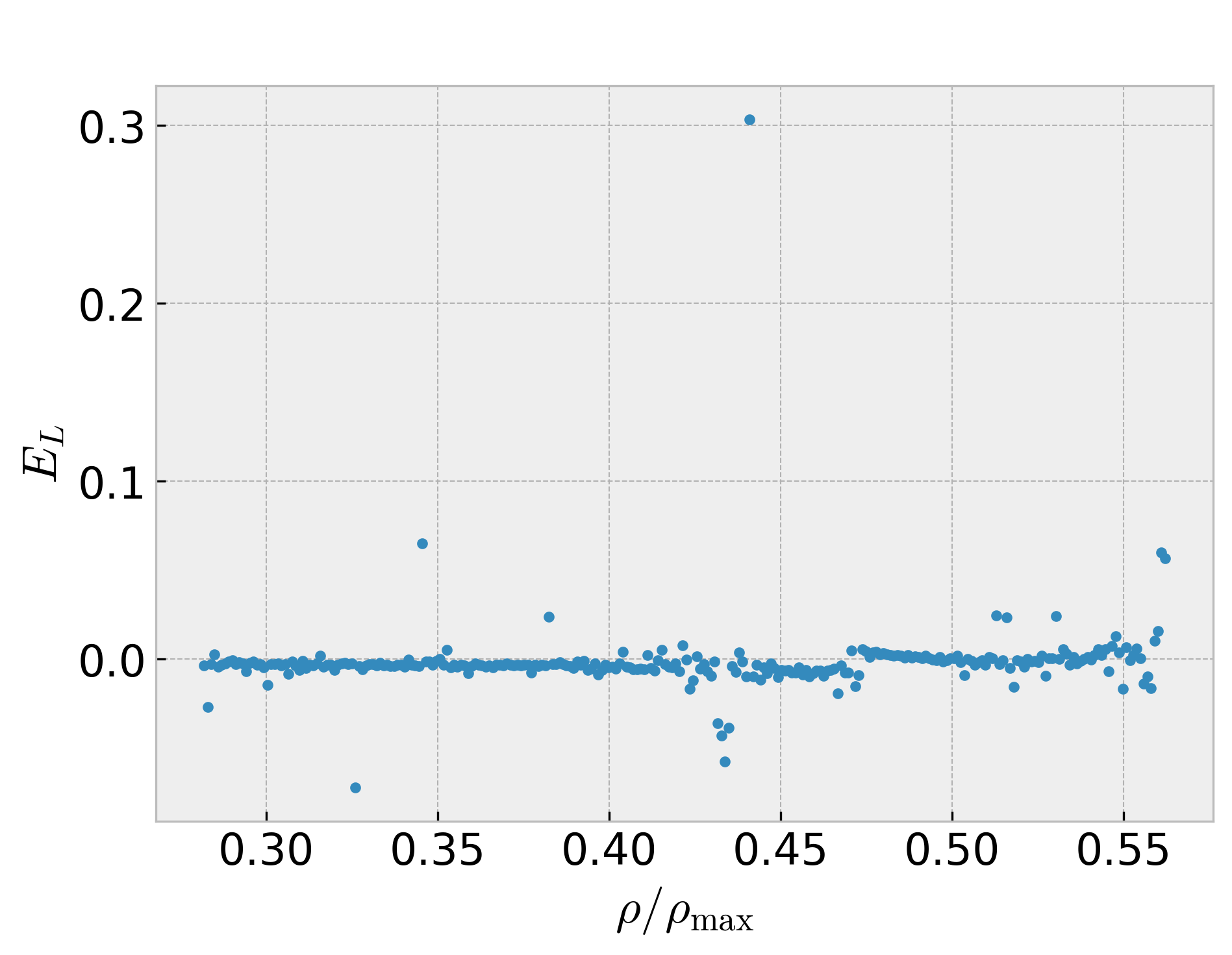}  
    \end{tabular} 
    (c) \\ 
    \includegraphics[width=.48\linewidth]{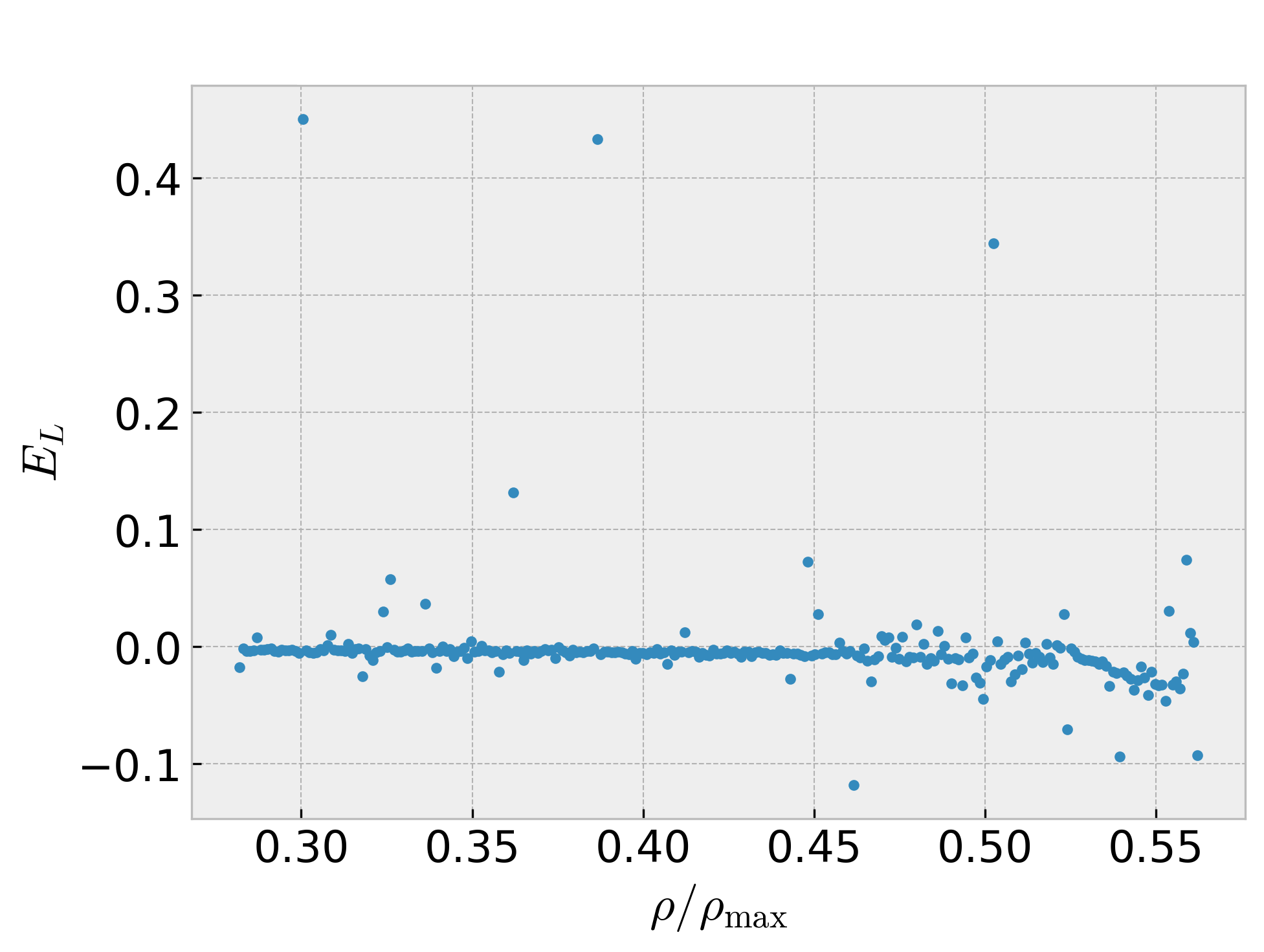}
   \caption{Example~7 (effect of the relaxation time): normalized length error~$E_L$ for (a)  $\tau=5$, 
    (b)  $\tau=10$ and (c) $\tau=1$.}%Largos y amplitudes para distintos $\tau$'s. Left: Error largos para $\tau=5$. Right: Error largos para $\tau=10$.}
   \label{img:dif_largo_taus}
\end{figure}

An interesting question   is whether the same phenomenon persists for different values of the relaxation time~$\tau$, 
 where we recall  that this parameter appears recurrently in the theoretical study of jamitons. To answer this question, the same experiment conducted previously is repeated, but with two additional different values of $\tau$ ($\tau=1$ and $\tau=10$). The outcome of the collisions is observed in Figure~\ref{img:colisiones_taus}(a).
The behavior of the exit velocities is consistent across different values of $\tau$. This makes sense because both $m$ and $s$ depend only on the sonic density, not on the value of $\tau$. This leads to the following conjecture:
\begin{conj}
    The sonic velocity of the jamiton resulting from a collision is independent of $\tau$.%La velocidad sónica del jamiton resultante de una colisión es independiente de $\tau$.
\end{conj}
The other properties of a jamiton, such as the amplitude, length, and $\rho_+$, can also be studied by varying $\tau$. In 
  Figures~\ref{img:colisiones_taus}(b)  to~(d), the amplitudes, lengths, and $\rho_+$ for $\tau=1$, $5$, and~$10$ are obtained.
It can be observed that such properties also do not depend on~$\tau$. This is also expected because both $\rho_+$ and $\rho_-$ depend only on $\rho_s$, which in turn is independent of $\tau$. This adds more invariants  to the collision of jamitons:

\begin{conj}
   The amplitude and value of $\rho_+$ of the jamiton resulting from a collision are independent of $\tau$.
\end{conj}
However, in Figure~\ref{img:colisiones_taus}(b), it is observed that the length of the post-collision jamiton changes, being greater as $\tau$ increases. This effect is related to the formula for length given in \eqref{eq:Largo}, which depends on $\tau$. This leads to the following conjecture:
\begin{conj}
    The length of the jamiton resulting from a collision depends on $\tau$.
\end{conj}
It is interesting to study whether there is additivity in the lengths of the jamitons, that is, whether the length of the jamiton post-collision corresponds to the sum of the lengths of the jamitons before colliding. For this, the lengths of the jamitons before collision are calculated, summed up, and compared to the lengths post-collision. The following normalized error is defined:
\begin{equation}\label{E_L}
    E_L \coloneqq  \frac{L_{\text{col}} - (L + L_{ \mathrm{test}})}{\displaystyle\max_{\rho}{\left(L + L_{ \mathrm{test}}\right)}},
\end{equation}
where $ L_{\text{col}} $ corresponds to the length of the jamiton post-collision, $ L_{ \mathrm{test}} $ to the length of the test jamiton, and $ L $ to the length of the jamiton collided with the test jamiton. Figure~\ref{img:dif_largo_taus}  shows the results obtained for $E_L$ in \eqref{E_L}. It is observed that indeed the error is small, especially for the case $\tau=5$. Additionally, the most noticeable differences could perfectly be caused by numerical errors in the simulation. This leads to the following conjecture:
\begin{conj}
	The length of the jamiton after colliding corresponds to the sum of the lengths of the colliding jamitons.
\end{conj}

\section{Conclusions}\label{cap:conclusion}
To conclude this work, a summary of the content and results in each section is provided, along with potential future lines of development.
Section \ref{cap:macro} introduced properties of the most used traffic models, being the inhomogeneous ARZ one of the key models studied here. In Section \ref{cap:jam}, the theoretical derivation of jamitons and their properties was observed. In Section \ref{cap:ARZ}, the numerical scheme was validated under numerical errors and comparisons with theoretical jamitons. Finally, in Section \ref{cap:colision}, several results in the collision of jamitons were defined and numerically tested: the collision of jamitons originates a new jamiton, the exit velocities correspond to the smoothing of two functions, the exit lengths increase, as well as the amplitudes (in a density range), and properties dependent and independent of $\tau$ post-collision were observed..

The collision of jamitons leaves several open questions to be answered, especially in the theoretical realm, as the results of this paper are only numerical but can provide direction for eventual proofs. Therefore, theoretical explanations or mathematically robust validations are tasks for future work. On the other hand, further experimentation with numerical simulations can be carried out to explore beyond double collisions. For example, simulating triple, quadruple collisions, or even collisions of $N$ jamitons. Of course, the results at this point must be carefully analyzed as they correspond to cases of greater complexity and whose final outcome is highly non-trivial. Similar cases have been studied in other types of systems, such as the Schrödinger equation where there is an explicit solution for the interaction of $N$ solitons \cite{hirota}. Following this line, one can also experiment with 'bombardments' using chains of jamitons on a larger-sized jamiton, hoping to decrease its amplitude and, in simpler terms, alleviate traffic congestion. The veracity of this effect is unknown, but Figure \ref{img:colisiones_vel} suggests that at least the speed of the larger jamiton can be increased.

Another experiment that can be carried out is to collide using other functions $h$, $U$, or $Q$ determined from another more localized fundamental diagram. For example, specific routes in a city or country. In this way, numerical simulations can be adjusted to different realities, as well as being able to apply them in real life and solve everyday problems.

Finally, a future work that can be carried out, although it deviates slightly from the theme of the paper, is to implement control models for traffic jams through numerical simulations of jamitons. As mentioned earlier, empirical control experiments have been conducted, but there is no mathematical rigor regarding them, so a good future work would be to experiment with these cases under numerical simulations and derive theoretical results.

\appendix

\bibliographystyle{natnumurl}
\bibliography{library_RB}

\begin{thebibliography}{10}
  \providecommand\backrefparscanfalse{}
  \providecommand\backrefparscantrue{}
  \providecommand\backrefprint{}

\backrefparscanfalse
\bibitem{presion}
Flynn, M.~R., Kasimov, A.~R., Nave, J.-C., Rosales, R.~R., and Seibold, B.,
  ``Self-sustained nonlinear waves in traffic flow'', { Phys. Rev. E (3)},
  vol.~79, no.~5, pp.~056113, 13, 2009,
  \href{https://dx.doi.org/10.1103/PhysRevE.79.056113}{\nolinkurl{doi:\detokenize{10.1103/PhysRevE.79.056113}}}.
\backrefprint\backrefparscantrue

\backrefparscanfalse
\bibitem{jamitones}
Seibold, B., Flynn, M.~R., Kasimov, A.~R., and Rosales, R.~R., ``Constructing
  set-valued fundamental diagrams from jamiton solutions in second order
  traffic models'', { Netw. Heterog. Media}, vol.~8, no.~3, pp.~745--772, 2013,
  \href{https://dx.doi.org/10.3934/nhm.2013.8.745}{\nolinkurl{doi:\detokenize{10.3934/nhm.2013.8.745}}}.
\backrefprint\backrefparscantrue

\backrefparscanfalse
\bibitem{Garavello-Piccoli}
Garavello, M., Han, K., and Piccoli, B., { Models for vehicular traffic on
  networks}, vol.~9 de { AIMS Series on Applied Mathematics}.
\newblock American Institute of Mathematical Sciences (AIMS), Springfield, MO,
  2016.
\backrefprint\backrefparscantrue

\backrefparscanfalse
\bibitem{Treiber-Kesting}
Treiber, M. and Kesting, A., { Traffic flow dynamics}.
\newblock Springer, Heidelberg, 2013,
  \href{https://dx.doi.org/10.1007/978-3-642-32460-4}{\nolinkurl{doi:\detokenize{10.1007/978-3-642-32460-4}}}.
\newblock Data, models and simulation, Translated by Treiber and Christian
  Thiemann.
\backrefprint\backrefparscantrue

\backrefparscanfalse
\bibitem{resurrection}
Aw, A. and Rascle, M., ``Resurrection of ``second order'' models of traffic
  flow'', { SIAM J. Appl. Math.}, vol.~60, no.~3, pp.~916--938, 2000,
  \href{https://dx.doi.org/10.1137/S0036139997332099}{\nolinkurl{doi:\detokenize{10.1137/S0036139997332099}}}.
\backrefprint\backrefparscantrue

\backrefparscanfalse
\bibitem{Z}
Zhang, H.~M., ``{A non-equilibrium traffic model devoid of gas-like
  behavior}'', { Transp. Res. Part B}, vol.~36, pp.~275--290, 2002,
  \href{https://dx.doi.org/10.1016/S0191-2615(00)00050-3}{\nolinkurl{doi:\detokenize{10.1016/S0191-2615(00)00050-3}}}.
\backrefprint\backrefparscantrue

\backrefparscanfalse
\bibitem{Na-Sch}
Nagel, K. and Schreckenberg, M., ``{A cellular automaton model for freeway
  traffic}'', { Journal de Physique I}, vol.~2, p.~2221, 1992,
  \href{https://dx.doi.org/10.1051/jp1:1992277}{\nolinkurl{doi:\detokenize{10.1051/jp1:1992277}}}.
\backrefprint\backrefparscantrue

\backrefparscanfalse
\bibitem{LW}
Lighthill, M.~J. and Whitham, G.~B., ``On kinematic waves. {II}. {A} theory of
  traffic flow on long crowded roads'', { Proc. Roy. Soc. London Ser. A},
  vol.~229, pp.~317--345, 1955,
  \href{https://dx.doi.org/10.1098/rspa.1955.0089}{\nolinkurl{doi:\detokenize{10.1098/rspa.1955.0089}}}.
\backrefprint\backrefparscantrue

\backrefparscanfalse
\bibitem{R}
Richards, P.~I., ``Shock waves on the highway'', { Operations Res.}, vol.~4,
  pp.~42--51, 1956,
  \href{https://dx.doi.org/10.1287/opre.4.1.42}{\nolinkurl{doi:\detokenize{10.1287/opre.4.1.42}}}.
\backrefprint\backrefparscantrue

\backrefparscanfalse
\bibitem{hw}
Hilliges, M. and Weidlich, W., ``A phenomenological model for dynamic traffic
  flow in networks'', { Transp. Res. Part B}, vol.~29, pp.~407--431, 1995,
  \href{https://dx.doi.org/10.1016/0191-2615(95)00018-9}{\nolinkurl{doi:\detokenize{10.1016/0191-2615(95)00018-9}}}.
\backrefprint\backrefparscantrue

\backrefparscanfalse
\bibitem{family}
B\"urger, R., Garc\'ia, A., Karlsen, K.~H., and Towers, J.~D., ``A family of
  numerical schemes for kinematic flows with discontinuous flux'', { J. Engrg.
  Math.}, vol.~60, no.~3-4, pp.~387--425, 2008,
  \href{https://dx.doi.org/10.1007/s10665-007-9148-4}{\nolinkurl{doi:\detokenize{10.1007/s10665-007-9148-4}}}.
\backrefprint\backrefparscantrue

\backrefparscanfalse
\bibitem{bcv23}
B\"urger, R., Contreras, H.~D., and Villada, L.~M., ``A
  {H}illiges-{W}eidlich-type scheme for a one-dimensional scalar conservation
  law with nonlocal flux'', { Netw. Heterog. Media}, vol.~18, no.~2,
  pp.~664--693, 2023,
  \href{https://dx.doi.org/10.3934/nhm.2023029}{\nolinkurl{doi:\detokenize{10.3934/nhm.2023029}}}.
\backrefprint\backrefparscantrue

\backrefparscanfalse
\bibitem{Jinbook2021}
Jin, W.-L., { Introduction to network traffic flow theory}.
\newblock Elsevier, Amsterdam, 2021,
  \href{https://dx.doi.org/10.1016/C2017-0-03540-1}{\nolinkurl{doi:\detokenize{10.1016/C2017-0-03540-1}}}.
\newblock Principles, concepts, models, and methods.
\backrefprint\backrefparscantrue

\backrefparscanfalse
\bibitem{micro-macro}
Aw, A., Klar, A., Materne, T., and Rascle, M., ``Derivation of continuum
  traffic flow models from microscopic follow-the-leader models'', { SIAM J.
  Appl. Math.}, vol.~63, no.~1, pp.~259--278, 2002,
  \href{https://dx.doi.org/10.1137/S0036139900380955}{\nolinkurl{doi:\detokenize{10.1137/S0036139900380955}}}.
\backrefprint\backrefparscantrue

\backrefparscanfalse
\bibitem{cellular-macro}
Alperovich, T. and Sopasakis, A., ``Stochastic description of traffic flow'', {
  J. Stat. Phys.}, vol.~133, no.~6, pp.~1083--1105, 2008,
  \href{https://dx.doi.org/10.1007/s10955-008-9652-6}{\nolinkurl{doi:\detokenize{10.1007/s10955-008-9652-6}}}.
\backrefprint\backrefparscantrue

\backrefparscanfalse
\bibitem{gps_macro}
Work, D., Tossavainen, O.-P., Blandin, S., Bayen, A., Iwuchukwu, T., and
  Tracton, K., ``{An Ensemble Kalman Filtering approach to highway traffic
  estimation using GPS enabled mobile devices}'', { Proceedings of the IEEE
  Conference on Decision and Control}, pp.~5062--5068, 2008,
  \href{https://dx.doi.org/10.1109/CDC.2008.4739016}{\nolinkurl{doi:\detokenize{10.1109/CDC.2008.4739016}}}.
\backrefprint\backrefparscantrue

\backrefparscanfalse
\bibitem{gps_macro_2}
Wang, Y. and Papageorgiou, M., ``{Real-time freeway traffic state estimation
  based on extended Kalman filter: a general approach}'', { Transp. Res. Part
  B}, vol.~39, no.~2, pp.~141--167, 2005,
  \href{https://dx.doi.org/https://doi.org/10.1016/j.trb.2004.03.003}{\nolinkurl{doi:\detokenize{https://doi.org/10.1016/j.trb.2004.03.003}}}.
\backrefprint\backrefparscantrue

\backrefparscanfalse
\bibitem{control}
Papageorgiou, M., ``{Some remarks on macroscopic traffic flow modelling}'', {
  Transp. Res. Part A}, vol.~32, no.~5, pp.~323--329, 1998,
  \href{https://dx.doi.org/https://doi.org/10.1016/S0965-8564(97)00048-7}{\nolinkurl{doi:\detokenize{https://doi.org/10.1016/S0965-8564(97)00048-7}}}.
\backrefprint\backrefparscantrue

\backrefparscanfalse
\bibitem{shocks}
LeVeque, R.~J., { Numerical methods for conservation laws}.
\newblock Lectures in Mathematics ETH Z\"urich, Birkh\"auser Verlag, Basel,
  second~ed., 1992,
  \href{https://dx.doi.org/10.1007/978-3-0348-8629-1}{\nolinkurl{doi:\detokenize{10.1007/978-3-0348-8629-1}}}.
\backrefprint\backrefparscantrue

\backrefparscanfalse
\bibitem{volumen_finito}
LeVeque, R.~J., { Finite volume methods for hyperbolic problems}.
\newblock Cambridge Texts in Applied Mathematics, Cambridge University Press,
  Cambridge, 2002,
  \href{https://dx.doi.org/10.1017/CBO9780511791253}{\nolinkurl{doi:\detokenize{10.1017/CBO9780511791253}}}.
\backrefprint\backrefparscantrue

\backrefparscanfalse
\bibitem{toro2009}
Toro, E.~F., { Riemann solvers and numerical methods for fluid dynamics}.
\newblock Springer-Verlag, Berlin, third~ed., 2009,
  \href{https://dx.doi.org/10.1007/b79761}{\nolinkurl{doi:\detokenize{10.1007/b79761}}}.
\newblock A practical introduction.
\backrefprint\backrefparscantrue

\backrefparscanfalse
\bibitem{hestbook}
Hesthaven, J.~S., { Numerical methods for conservation laws}, vol.~18 de {
  Computational Science \& Engineering}.
\newblock Society for Industrial and Applied Mathematics (SIAM), Philadelphia,
  PA, 2018,
  \href{https://dx.doi.org/10.1137/1.9781611975109}{\nolinkurl{doi:\detokenize{10.1137/1.9781611975109}}}.
\newblock From analysis to algorithms.
\backrefprint\backrefparscantrue

\backrefparscanfalse
\bibitem{kuzminbook}
Kuzmin, D. and Hajduk, H., { Property-preserving numerical schemes for
  conservation laws}.
\newblock World Scientific Publishing Co. Pte. Ltd., Hackensack, NJ, [2024]
  \copyright 2024.
\backrefprint\backrefparscantrue

\backrefparscanfalse
\bibitem{ablo11}
Ablowitz, M.~J., { Nonlinear dispersive waves}.
\newblock Cambridge Texts in Applied Mathematics, Cambridge University Press,
  New York, 2011,
  \href{https://dx.doi.org/10.1017/CBO9780511998324}{\nolinkurl{doi:\detokenize{10.1017/CBO9780511998324}}}.
\newblock Asymptotic analysis and solitons.
\backrefprint\backrefparscantrue

\backrefparscanfalse
\bibitem{schneid17}
Schneider, G. and Uecker, H., { Nonlinear {PDE}s}, vol.~182 de { Graduate
  Studies in Mathematics}.
\newblock American Mathematical Society, Providence, RI, 2017,
  \href{https://dx.doi.org/10.1090/gsm/182}{\nolinkurl{doi:\detokenize{10.1090/gsm/182}}}.
\newblock A dynamical systems approach.
\backrefprint\backrefparscantrue

\backrefparscanfalse
\bibitem{phi4}
Goodman, R.~H. and Haberman, R., ``Kink-antikink collisions in the {$\phi^4$}
  equation: the {$n$}-bounce resonance and the separatrix map'', { SIAM J.
  Appl. Dyn. Syst.}, vol.~4, no.~4, pp.~1195--1228, 2005,
  \href{https://dx.doi.org/10.1137/050632981}{\nolinkurl{doi:\detokenize{10.1137/050632981}}}.
\backrefprint\backrefparscantrue

\backrefparscanfalse
\bibitem{phi_4_1}
Campbell, D.~K., Peyrard, M., and Sodano, P., ``Kink-antikink interactions in
  the double sine-{G}ordon equation'', { Phys. D}, vol.~19, no.~2,
  pp.~165--205, 1986,
  \href{https://dx.doi.org/10.1016/0167-2789(86)90019-9}{\nolinkurl{doi:\detokenize{10.1016/0167-2789(86)90019-9}}}.
\backrefprint\backrefparscantrue

\backrefparscanfalse
\bibitem{phi_4_2}
Campbell, D.~K. and Peyrard, M., ``Solitary wave collisions revisited'', {
  Phys. D}, vol.~18, no.~1-3, pp.~47--53, 1986,
  \href{https://dx.doi.org/10.1016/0167-2789(86)90161-2}{\nolinkurl{doi:\detokenize{10.1016/0167-2789(86)90161-2}}}.
\newblock Solitons and coherent structures (Santa Barbara, Calif., 1985).
\backrefprint\backrefparscantrue

\backrefparscanfalse
\bibitem{phi_4_3}
Campbell, D.~K., Schonfeld, J.~F., and Wingate, C.~A., ``{Resonance structure
  in kink-antikink interactions in $\varphi^4$ theory}'', { Phys. D}, vol.~9,
  no.~1, pp.~1--32, 1983,
  \href{https://dx.doi.org/https://doi.org/10.1016/0167-2789(83)90289-0}{\nolinkurl{doi:\detokenize{https://doi.org/10.1016/0167-2789(83)90289-0}}}.
\backrefprint\backrefparscantrue

\backrefparscanfalse
\bibitem{phi_4_4}
Anninos, P., Oliveira, S., and Matzner, R.~A., ``{Fractal structure in the
  scalar
  $\ensuremath{\lambda}{({\ensuremath{\varphi}}^{2}\ensuremath{-}1)}^{2}$
  theory}'', { Phys. Rev. D}, vol.~44, pp.~1147--1160, 1991,
  \href{https://dx.doi.org/10.1103/PhysRevD.44.1147}{\nolinkurl{doi:\detokenize{10.1103/PhysRevD.44.1147}}}.
\backrefprint\backrefparscantrue

\backrefparscanfalse
\bibitem{phi_4_5}
Peyrard, M. and Campbell, D.~K., ``{Kink-antikink interactions in a modified
  sine-{Gordon} model}'', { Physica D: Nonlinear Phenomena}, vol.~9, no.~1,
  pp.~33--51, 1983,
  \href{https://dx.doi.org/https://doi.org/10.1016/0167-2789(83)90290-7}{\nolinkurl{doi:\detokenize{https://doi.org/10.1016/0167-2789(83)90290-7}}}.
\backrefprint\backrefparscantrue

\backrefparscanfalse
\bibitem{Bellouquid2012}
Bellouquid, A., De~Angelis, E., and Fermo, L., ``Towards the modeling of
  vehicular traffic as a complex system: a kinetic theory approach'', { Math.
  Models Methods Appl. Sci.}, vol.~22, pp.~1140003, 35, 2012,
  \href{https://dx.doi.org/10.1142/S0218202511400033}{\nolinkurl{doi:\detokenize{10.1142/S0218202511400033}}}.
\backrefprint\backrefparscantrue

\backrefparscanfalse
\bibitem{Albi2019}
Albi, G., Bellomo, N., Fermo, L., Ha, S.-Y., Kim, J., Pareschi, L., Poyato, D.,
  and Soler, J., ``Vehicular traffic, crowds, and swarms: from kinetic theory
  and multiscale methods to applications and research perspectives'', { Math.
  Models Methods Appl. Sci.}, vol.~29, no.~10, pp.~1901--2005, 2019,
  \href{https://dx.doi.org/10.1142/S0218202519500374}{\nolinkurl{doi:\detokenize{10.1142/S0218202519500374}}}.
\backrefprint\backrefparscantrue

\backrefparscanfalse
\bibitem{Delitala2017}
Delitala, M. and Tosin, A., ``Mathematical modeling of vehicular traffic: a
  discrete kinetic theory approach'', { Math. Models Methods Appl. Sci.},
  vol.~17, no.~6, pp.~901--932, 2007,
  \href{https://dx.doi.org/10.1142/S0218202507002157}{\nolinkurl{doi:\detokenize{10.1142/S0218202507002157}}}.
\backrefprint\backrefparscantrue

\backrefparscanfalse
\bibitem{P}
Payne, H.~J., ``Models of freeway traffic and control'', en { Mathematical
  models of public systems} (Bekey, G.~A., ed.), vol.~No. 1 de { Simulation
  Councils Proceedings Series, Vol. 1}, pp.~51--61, Simulation Councils, Inc.,
  1971.
\newblock National Invitational Seminar on Advanced Simulation, held in San
  Diego, Calif., September 1970.
\backrefprint\backrefparscantrue

\backrefparscanfalse
\bibitem{W}
Whitham, G.~B., { Linear and nonlinear waves}.
\newblock Pure and Applied Mathematics, Wiley-Interscience [John Wiley \&
  Sons], New York-London-Sydney, 1974.
\backrefprint\backrefparscantrue

\backrefparscanfalse
\bibitem{dataset}
{Minnesota Department of Transportation}, ``{Mn/DOT Traffic Data}'', 2017,
  \href{http://data.dot.state.mn.us/datatools/}{\nolinkurl{\detokenize{http://data.dot.state.mn.us/datatools/}}}.
\backrefprint\backrefparscantrue

\backrefparscanfalse
\bibitem{primer_DF}
Greenshields, B.~D., Bibbins, J.~R., Channing, W., and Miller, H.~H., ``{A
  study of traffic capacity}'', { Highway Research Board}, vol.~14,
  pp.~448--477, 1935,
  \href{https://api.semanticscholar.org/CorpusID:107546777}{\nolinkurl{\detokenize{https://api.semanticscholar.org/CorpusID:107546777}}}.
\backrefprint\backrefparscantrue

\backrefparscanfalse
\bibitem{Newell}
Newell, G., ``{A simplified theory of kinematic waves in highway traffic, part
  II: Queueing at freeway bottlenecks}'', { Transp. Res. Part B}, vol.~27,
  no.~4, pp.~289--303, 1993,
  \href{https://dx.doi.org/https://doi.org/10.1016/0191-2615(93)90039-D}{\nolinkurl{doi:\detokenize{https://doi.org/10.1016/0191-2615(93)90039-D}}}.
\backrefprint\backrefparscantrue

\backrefparscanfalse
\bibitem{Daganzo}
Daganzo, C.~F., ``{The cell transmission model: A dynamic representation of
  highway traffic consistent with the hydrodynamic theory}'', { Transp. Res.
  Part B}, vol.~28, no.~4, pp.~269--287, 1994,
  \href{https://dx.doi.org/https://doi.org/10.1016/0191-2615(94)90002-7}{\nolinkurl{doi:\detokenize{https://doi.org/10.1016/0191-2615(94)90002-7}}}.
\backrefprint\backrefparscantrue

\backrefparscanfalse
\bibitem{problema_riemann}
Mammar, S., Lebacque, J., and Salem, H.~H., ``Riemann problem resolution and
  {G}odunov scheme for the {A}w-{R}ascle-{Z}hang model'', { Transp. Sci.},
  vol.~43, no.~4, pp.~531--545, 2009,
  \href{http://www.jstor.org/stable/25769472}{\nolinkurl{\detokenize{http://www.jstor.org/stable/25769472}}}.
\backrefprint\backrefparscantrue

\backrefparscanfalse
\bibitem{df_teoria}
Kerner, B., { Introduction to Modern Traffic Flow Theory and Control: The Long
  Road to Three-Phase Traffic Theory}.
\newblock Springer Berlin, Heidelberg, 2009,
  \href{https://dx.doi.org/10.1007/978-3-642-02605-8}{\nolinkurl{doi:\detokenize{10.1007/978-3-642-02605-8}}}.
\backrefprint\backrefparscantrue

\backrefparscanfalse
\bibitem{parametros}
Fan, S. and Seibold, B., ``Data-fitted first-order traffic models and their
  second-order generalizations: Comparison by trajectory and sensor data'', {
  Transp. Res. Record}, vol.~2391, no.~1, pp.~32--43, 2013,
  \href{https://dx.doi.org/10.3141/2391-04}{\nolinkurl{doi:\detokenize{10.3141/2391-04}}}.
\backrefprint\backrefparscantrue

\backrefparscanfalse
\bibitem{taco_fantasma}
Sugiyama, Y., Fukui, M., Kikuchi, M., Hasebe, K., Nakayama, A., Nishinari, K.,
  Tadaki, S., and Yukawa, S., ``{Traffic jams without
  bottlenecks—experimental evidence for the physical mechanism of the
  formation of a jam}'', { New J. Phys.}, vol.~10, p.~033001, 2008,
  \href{https://dx.doi.org/10.1088/1367-2630/10/3/033001}{\nolinkurl{doi:\detokenize{10.1088/1367-2630/10/3/033001}}}.
\backrefprint\backrefparscantrue

\backrefparscanfalse
\bibitem{Requiem}
Daganzo, C.~F., ``{Requiem for second-order fluid approximations of traffic
  flow}'', { Transp. Res. Part B}, vol.~29, no.~4, pp.~277--286, 1995,
  \href{https://dx.doi.org/https://doi.org/10.1016/0191-2615(95)00007-Z}{\nolinkurl{doi:\detokenize{https://doi.org/10.1016/0191-2615(95)00007-Z}}}.
\backrefprint\backrefparscantrue

\backrefparscanfalse
\bibitem{lagrangeano_1}
Greenberg, J.~M., ``Congestion redux'', { SIAM J. Appl. Math.}, vol.~64, no.~4,
  pp.~1175--1185, 2004,
  \href{https://dx.doi.org/10.1137/S0036139903431737}{\nolinkurl{doi:\detokenize{10.1137/S0036139903431737}}}.
\backrefprint\backrefparscantrue

\backrefparscanfalse
\bibitem{lagrangeano_2}
Courant, R. and Friedrichs, K.~O., { Supersonic {F}low and {S}hock {W}aves}.
\newblock Interscience Publishers, Inc., New York, 1948.
\backrefprint\backrefparscantrue

\backrefparscanfalse
\bibitem{SPH-1}
Gingold, R.~A. and Monaghan, J.~J., ``{Smoothed particle hydrodynamics: theory
  and application to non-spherical stars}'', { Mon. Not. R. Astron. Soc.},
  vol.~181, pp.~375--389, 1977,
  \href{https://dx.doi.org/10.1093/mnras/181.3.375}{\nolinkurl{doi:\detokenize{10.1093/mnras/181.3.375}}}.
\backrefprint\backrefparscantrue

\backrefparscanfalse
\bibitem{SPH-2}
Lucy, L.~B., ``{Numerical approach to the testing of the fission hypothesis}'',
  { Astron. J. (United States)}, vol.~82:12, 1977,
  \href{https://dx.doi.org/10.1086/112164}{\nolinkurl{doi:\detokenize{10.1086/112164}}}.
\backrefprint\backrefparscantrue

\backrefparscanfalse
\bibitem{CSC_demo}
Whitham, G.~B., ``Some comments on wave propagation and shock wave structure
  with application to magnetohydrodynamics'', { Comm. Pure Appl. Math.},
  vol.~12, pp.~113--158, 1959,
  \href{https://dx.doi.org/10.1002/cpa.3160120107}{\nolinkurl{doi:\detokenize{10.1002/cpa.3160120107}}}.
\backrefprint\backrefparscantrue

\backrefparscanfalse
\bibitem{stabilidad_jam}
Ramadan, R., Rosales, R.~R., and Seibold, B., ``Structural properties of the
  stability of jamitons'', en { Mathematical descriptions of traffic flow:
  micro, macro and kinetic models}, vol.~12 de { ICIAM 2019 SEMA SIMAI Springer
  Ser.}, pp.~35--62, Springer, Cham, [2021] \copyright 2021,
  \href{https://dx.doi.org/10.1007/978-3-030-66560-9\_3}{\nolinkurl{doi:\detokenize{10.1007/978-3-030-66560-9\_3}}}.
\backrefprint\backrefparscantrue

\backrefparscanfalse
\bibitem{whitham_1}
Chen, G.~Q., Levermore, C.~D., and Liu, T.-P., ``Hyperbolic conservation laws
  with stiff relaxation terms and entropy'', { Comm. Pure Appl. Math.},
  vol.~47, no.~6, pp.~787--830, 1994,
  \href{https://dx.doi.org/10.1002/cpa.3160470602}{\nolinkurl{doi:\detokenize{10.1002/cpa.3160470602}}}.
\backrefprint\backrefparscantrue

\backrefparscanfalse
\bibitem{whitham_2}
Liu, T.-P., ``Hyperbolic conservation laws with relaxation'', { Comm. Math.
  Phys.}, vol.~108, no.~1, pp.~153--175, 1987,
  \href{http://projecteuclid.org/euclid.cmp/1104116362}{\nolinkurl{\detokenize{http://projecteuclid.org/euclid.cmp/1104116362}}}.
\backrefprint\backrefparscantrue

\backrefparscanfalse
\bibitem{whitham_3}
Li, T. and Liu, H., ``Stability of a traffic flow model with nonconvex
  relaxation'', { Commun. Math. Sci.}, vol.~3, no.~2, pp.~101--118, 2005,
  \href{http://projecteuclid.org/euclid.cms/1118778270}{\nolinkurl{\detokenize{http://projecteuclid.org/euclid.cms/1118778270}}}.
\backrefprint\backrefparscantrue

\backrefparscanfalse
\bibitem{deg_1}
Li, T., ``Global solutions and zero relaxation limit for a traffic flow
  model'', { SIAM J. Appl. Math.}, vol.~61, no.~3, pp.~1042--1061, 2000,
  \href{https://dx.doi.org/10.1137/S0036139999356788}{\nolinkurl{doi:\detokenize{10.1137/S0036139999356788}}}.
\backrefprint\backrefparscantrue

\backrefparscanfalse
\bibitem{deg_2}
Li, T. and Liu, H., ``Critical thresholds in a relaxation system with resonance
  of characteristic speeds'', { Discrete Contin. Dyn. Syst.}, vol.~24, no.~2,
  pp.~511--521, 2009,
  \href{https://dx.doi.org/10.3934/dcds.2009.24.511}{\nolinkurl{doi:\detokenize{10.3934/dcds.2009.24.511}}}.
\backrefprint\backrefparscantrue

\backrefparscanfalse
\bibitem{Ramadan}
Ramadan, R., ``Non-equilibrium dynamics of second order traffic models'', 2020,
  \href{http://hdl.handle.net/20.500.12613/2078}{\nolinkurl{\detokenize{http://hdl.handle.net/20.500.12613/2078}}}.
\backrefprint\backrefparscantrue

\backrefparscanfalse
\bibitem{ZND}
Fickett, W. and Davis, W.~C., { Detonation: Theory and Experiment}.
\newblock Dover Books on Physics, Dover Publications, 2012,
  \href{https://books.google.cl/books?id=QaejAQAAQBAJ}{\nolinkurl{\detokenize{https://books.google.cl/books?id=QaejAQAAQBAJ}}}.
\backrefprint\backrefparscantrue

\backrefparscanfalse
\bibitem{hllpaper}
Harten, A., Lax, P.~D., and van Leer, B., ``On upstream differencing and
  {G}odunov-type schemes for hyperbolic conservation laws'', { SIAM Rev.},
  vol.~25, no.~1, pp.~35--61, 1983,
  \href{https://dx.doi.org/10.1137/1025002}{\nolinkurl{doi:\detokenize{10.1137/1025002}}}.
\backrefprint\backrefparscantrue

\backrefparscanfalse
\bibitem{entropia}
Lax, P.~D., { Hyperbolic systems of conservation laws and the mathematical
  theory of shock waves}, vol.~No. 11 de { Conference Board of the Mathematical
  Sciences Regional Conference Series in Applied Mathematics}.
\newblock Society for Industrial and Applied Mathematics, Philadelphia, PA,
  1973.
\backrefprint\backrefparscantrue

\backrefparscanfalse
\bibitem{MIT_jam}
{Massachusetts Institute of Technology (MIT)}, ``{Traffic Modeling - Phantom
  Traffic Jams and Traveling Jamitons}''.,
  \href{https://math.mit.edu/traffic/}{\nolinkurl{\detokenize{https://math.mit.edu/traffic/}}}.
\backrefprint\backrefparscantrue

\backrefparscanfalse
\bibitem{hirota}
Hirota, R., ``Exact envelope-soliton solutions of a nonlinear wave equation'',
  { J. Mathematical Phys.}, vol.~14, pp.~805--809, 1973,
  \href{https://dx.doi.org/10.1063/1.1666399}{\nolinkurl{doi:\detokenize{10.1063/1.1666399}}}.
\backrefprint\backrefparscantrue

\end{thebibliography}

\end{document}